\documentclass[12pt]{amsart}
\font\emailfont=cmtt10

\usepackage{amsmath,amsthm,amsfonts,amscd,flafter,epsf}

\newcommand\divis{\mathfrak d}
\newcommand\ord{\mathrm{ord}}
\newcommand\commentable[1]{#1}
\newcommand\Id{\mathrm{Id}}

\newcommand{\Tors}{\mathrm{Tors}}
\newcommand{\rk}{\mathrm{rk}}
\newcommand{\HF}{HF}

\newtheorem{theorem}{Theorem}[section]
\newtheorem{prop}[theorem]{Proposition}
\newtheorem{cor}[theorem]{Corollary}
\newtheorem{conj}[theorem]{Conjecture}
\newtheorem{lemma}[theorem]{Lemma}

\newtheorem{defn}[theorem]{Definition}

\newtheorem{remark}[theorem]{Remark}

\def\endproof{\relax\ifmmode\expandafter\endproofmath\else
  \unskip\nobreak\hfil\penalty50\hskip.75em\hbox{}\nobreak\hfil\bull
  {\parfillskip=0pt \finalhyphendemerits=0 \bigbreak}\fi}
\def\endproofmath$${\eqno\bull$$\bigbreak}
\def\bull{\vbox{\hrule\hbox{\vrule\kern3pt\vbox{\kern6pt}\kern3pt\vrule}\hrule}}

\newcounter{bean}
\newcommand{\Q}{\mathbb{Q}}
\newcommand{\R}{\mathbb{R}}
\newcommand{\T}{\mathbb{T}}

\newcommand{\C}{\mathbb{C}}

\newcommand{\Z}{\mathbb{Z}}

\newcommand{\OneHalf}{\frac{1}{2}}

\newcommand{\OneQuarter}{\frac{1}{4}}

\newcommand{\Zmod}[1]{\Z/{#1}\Z}

\newcommand{\Ker}{\mathrm{Ker}}

\newcommand{\Image}{\mathrm{Im}}

\newcommand{\grad}{\vec\nabla}

\newcommand{\cm}{\cdot}

\newcommand{\Nbd}[1]{{\mathrm{nd}}(#1)}
\newcommand{\nbd}[1]{\Nbd{#1}}
\newcommand{\CDisk}{D}

\newcommand{\ModSWfour}{\mathcal{M}}
\newcommand{\ModFlow}{\ModSWfour}

\newcommand{\upm}{{\widehat {\ModFlow}}}
\newcommand{\HFswred}{HF^{SW} _{{\rm red}}}

\newcommand{\SpinC}{{\mathrm{Spin}}^c}

\newcommand{\goesto}{\mapsto}

\newcommand\Wedge{\Lambda}

\newcommand\Hom{\mathrm{Hom}}

\newcommand\abuts\Rightarrow
\newcommand\Sym{\mathrm{Sym}}

\newcommand\Ann{\mathbb A}
\newcommand\AnnOp{{\mathbb A}^{\circ}}

\def\sqr#1#2{{\vcenter{\vbox{\hrule height.#2pt
	\hbox{\vrule width.#2pt height#1pt \kern#1pt
	\vrule width.#2pt}
	\hrule height.#2pt}}}}

\newcommand\HFaCoh{{\widehat{HF}}^{*}}
\newcommand\HFpCoh{HF_{+}^{*}}
\newcommand\HFmCoh{HF_{-}^{*}}

\newcommand{\sj}{\mathfrak j}
\newcommand{\MT}{t}
\newcommand{\Tor}{\mathrm{Tor}}

\newcommand\HFmRed{\HFm_{\red}}

\newcommand\sRelSpinC{\underline\spinc}
\newcommand\RelSpinC{\underline{\SpinC}}
\newcommand\relspinc{\underline{\spinc}}

\newcommand\Filt{\mathcal F}
\newcommand\HFinfty{\HFinf}
\newcommand\CFinfty{\CFinf}
\newcommand\Tai{{\mathbb T}_{\alpha}^i}
\newcommand\Tbj{{\mathbb T}_{\beta}^j}
\newcommand\RightFp{R^+}
\newcommand\LeftFp{L^+}
\newcommand\RightFinf{R^\infty}
\newcommand\LeftFinf{L^\infty}
\newcommand\Area{\mathcal A}
\newcommand\PhiIn{\phi^{\mathrm{in}}}
\newcommand\PhiOut{\phi^{\mathrm{out}}}

\newcommand\x{\mathbf x}
\newcommand\w{\mathbf w}

\newcommand\p{\mathbf p}
\newcommand\q{\mathbf q}
\newcommand\y{\mathbf y}

\newcommand\ModSphere{\ModFlow\left({\mathbb S}\longrightarrow 
\Sym^{g-1}(\Sigma_{1})\times \Sym^2(\Sigma_{2})\right)}
\newcommand\ModSpheres\ModSphere
\newcommand\CF{CF}

\newcommand\CFa{\widehat{CF}}
\newcommand\CFp{\CFb}
\newcommand\CFm{\CF^-}
\newcommand\CFleq{\CF^{\leq 0}}
\newcommand\HFleq{\HF^{\leq 0}}

\newcommand\HFred{\HF_{\rm red}}
\newcommand{\red}{{\rm red}}

\newcommand\HFp{\HFb}
\newcommand\HFpm{HF^{\pm}}
\newcommand\HFm{\HF^-}
\newcommand\CFinf{CF^\infty}
\newcommand\HFinf{HF^\infty}
\newcommand\CFb{CF^+}
\newcommand\HFa{\widehat{HF}}
\newcommand\HFb{HF^+}
\newcommand\gr{\mathrm{gr}}
\newcommand\Mas{\mu}
\newcommand\UnparModSp{\widehat \ModSp}
\newcommand\UnparModFlow\UnparModSp
\newcommand\Mod\ModSp
\newcommand{\cpl}{{\mathcal C}^+}
\newcommand{\cmi}{{\mathcal C}^-}
\newcommand{\cplm}{{\mathcal C}^\pm}
\newcommand{\cald}{{\mathcal D}}

\newcommand\PD{\mathrm{PD}}

\newcommand{\spinc}{\mathfrak s}

\newcommand{\spinct}{\mathfrak t}

\newcommand\ModMaps{\mathcal M}
\newcommand\ModSp\ModMaps

\newcommand\Ta{{\mathbb T}_{\alpha}}
\newcommand\Tb{{\mathbb T}_{\beta}}
\newcommand\Tc{{\mathbb T}_{\gamma}}
\newcommand\Td{{\mathbb T}_{\delta}}

\newcommand\Strip{\mathbb{D}}

\newcommand\alphas{\mbox{\boldmath$\alpha$}}
\newcommand\xis{\mbox{\boldmath$\xi$}}
\newcommand\etas{\mbox{\boldmath$\eta$}}
\newcommand\betas{\mbox{\boldmath$\beta$}}
\newcommand\gammas{\mbox{\boldmath$\gamma$}}
\newcommand\deltas{\mbox{\boldmath$\delta$}}
\newcommand\HFto{HF _{{\rm to}} ^{SW}}
\newcommand\CFto{CF _{{\rm to}} ^{SW}}
\newcommand\HFfrom{HF_{{\rm from}} ^{SW}}
\newcommand\CFfrom{CF_{{\rm from}} ^{SW}}

\newcommand\PerDom{\mathcal P}

\newcommand\intPerDom{\mathrm{int}\PerDom}
\newcommand\csum{*}
\newcommand\CurveComp{\Sigma-\alpha_{1}-\ldots-\alpha_{g}-\beta_{1}-\ldots-\beta_{g}}
\newcommand\EmbSurf{Z}

\newcommand\uCF{\underline\CF}
\newcommand\uCFinf{\underline\CFinf}
\newcommand\uHF{\underline{\HF}}
\newcommand\uHFred{\underline{\HFred}}
\newcommand\uHFp{\underline{\HFp}}
\newcommand\uHFpm{\underline{\HFpm}}
\newcommand\uCFp{\underline{\CFp}}
\newcommand\uHFm{\underline{\HFm}}
\newcommand\uHFa{\underline{\HFa}}
\newcommand\uCFa{\underline{\CFa}}
\newcommand\uDel{\underline{\partial}}
\newcommand\uHFinf{\underline{\HFinf}}

\newcommand\ufinf{\underline{f^\infty_{\alpha,\beta,\gamma}}}
\newcommand\uHFleq{\underline{\HFleq}}
\newcommand\uCFleq{\underline{\CFleq}}
\newcommand\ufp[1]{\underline{f^+_{#1}}}
\newcommand\ufleq[1]{\underline{f^{\leq 0}_{#1}}}
\newcommand\ufa{\underline{{\widehat f}_{\alpha,\beta,\gamma}}}
\newcommand\uFa[1]{\underline{{\widehat F}_{#1}}}
\newcommand\uFm[1]{\underline{F^-_{#1}}}

\newcommand\uFp[1]{\underline{F^+_{#1}}}
\newcommand\uFinf[1]{\underline{F^\infty_{#1}}}
\newcommand\uFleq[1]{\underline{F^{\leq 0}_{#1}}}
\newcommand\uFstar[1]{\underline{F^{*}_{#1}}}

\newcommand\Fp[1]{F^{+}_{#1}}

\newcommand\Tx{{\mathbb T}_\xi}
\newcommand\Ty{{\mathbb T}_\eta}

\newcommand\uLeftFinf{{\underline L}^\infty}
\newcommand\uRightFinf{{\underline R}^\infty}
\newcommand\orient{\mathfrak o}
\newcommand\Field{\mathbb F}

\newcommand\TbPr{\Tb'}
\newcommand\TxPr{\Tx'}
\newcommand\SumMap{\Gamma}

\title[{Holomorphic disks and three-manifold invariants}]
{Holomorphic disks and three-manifold invariants: properties and applications}

\author[Peter Ozsv{\'a}th]{Peter Ozsv\'ath}
\address{Department of
Mathematics,  Princeton University, New Jersey 08540 \newline
\indent{\emailfont{petero@math.princeton.edu}}}\thanks{PSO was supported by NSF grant number DMS 9971950 and a Sloan 
Research Fellowship}

\author[Zolt{\'a}n Szab{\'o}]{Zolt{\'a}n Szab{\'o}} 
\address{Department of
Mathematics,  Princeton University, New Jersey 08540 \newline
\indent{\emailfont{szabo@math.princeton.edu}}}
\thanks{ZSz was supported by a Sloan 
Research Fellowship and a Packard Fellowship}

\begin{document}

\begin{abstract} 
    In \cite{HolDisk}, we introduced Floer homology theories
    $\HFm(Y,\spinc)$, $\HFinf(Y,\spinc)$, $\HFp(Y,\spinct)$,
    $\HFa(Y,\spinc) $,and
    $\HFred(Y,\spinc)$ associated to closed, oriented three-manifolds
    $Y$ equipped with a $\SpinC$ structures $\spinc \in \SpinC(Y)$. In
    the present paper, we give calculations and study the properties
    of these invariants. The calculations suggest a conjectured
    relationship with Seiberg-Witten theory.  The properties include a
    relationship between the Euler characteristics of $\HFpm$ and
    Turaev's torsion, a relationship with the minimal genus problem
    (Thurston norm), and surgery exact sequences. We also include some
    applications of these techniques to three-manifold topology.
\end{abstract}

\maketitle
\section{Introduction}
\label{sec:Introduction}

The present paper is a continuation of~\cite{HolDisk}, where we
defined topological invariants for closed oriented, three-manifolds
$Y$, equipped with a $\SpinC$ structure $\spinc\in\SpinC(Y)$. These
invariants are a collection of Floer homology groups $\HFm(Y,\spinc)$,
$\HFinf(Y,\spinc)$, $\HFp(Y,\spinc)$, and $\HFa(Y,\spinc)$. Our goal
here is to study these invariants: calculate them in several examples,
establish their fundamental properties, and give some applications.

We begin in Section~\ref{sec:BasProp} with some of the properties of
the groups, including their behaviour under orientation reversal of $Y$ and
conjugation of its $\SpinC$ structures. Moreover, we show that for any
three-manifold $Y$, there are at most finitely many $\SpinC$
structures $\spinc\in\SpinC(Y)$ with the property that
$\HFp(Y,\spinc)$ is non-trivial.\footnote{Throughout this
introduction, and indeed through most of this paper, we will suppress
the orientation system $\orient$ used in the definition. This is
justified in part by the fact that our statements typically hold for
all possible orientation systems on $Y$ (and if not, then it is
easy to supply necessary quantifiers). A more compelling justification
is given by the fact that in Section~\ref{sec:HFinfty}, we show how to
equip an arbitrary oriented three-manifold with $b_1(Y)>0$ with a
canonical orientation system. And finally, of course, orientation
systems become irrelevant if we were to work with coefficients in $\Zmod{2}$}

In Section~\ref{sec:Examples}, we illustrate the Floer homology
theories by computing the invariants for certain rational homology
three-spheres. These calculations are done by explicitly identifying
the relevant moduli spaces of flow-lines. In Section~\ref{sec:SW} we
compare them with invariants with corresponding ``equivariant
Seiberg-Witten-Floer homologies''$\HFto$, $\HFfrom$, and $\HFswred$
for the three-manifolds studied in Section~\ref{sec:Examples},
compare~\cite{Marcolli},
\cite{KMIrvine}.

These calculations support the following conjecture:

\begin{conj}
\label{Conj}
    Let $Y$ be an oriented rational homology three-sphere. Then for
all $\SpinC$ structures  $\spinc\in\SpinC(Y)$ 
    there are isomorphisms\footnote{This manuscript was written before the appearance of~\cite{Manolescu} and~\cite{KManolescu}. 
In the second paper, Kronheimer and Manolescu propose alternate
Seiberg-Witten constructions, and indeed give one which they
conjecture to agree with our $\HFa$, see also~\cite{MarcolliWangNew}.}
    $$\HFto(Y,\spinc)\cong \HFp(Y,\spinc),\ \ 
 	\HFfrom(Y,\spinc)\cong \HFm(Y,\spinc),\ \
	\HFswred(Y,\spinc)\cong \HFred(Y,\spinc).$$
\end{conj}

After the specific calculations, we turn back to general
properties. In Section~\ref{sec:EulerCharacteristic}, we consider the
Euler characteristics of the theories. The Euler characteristic of
$\HFa(Y,\spinc)$ turns out to depend only on homological information
of $Y$, but the Euler characteristic of $\HFp$ has richer structure:
indeed, when $b_1(Y)>0$, we establish a relationship between it and
Turaev's torsion function (c.f. Theorem~\ref{thm:EulerOne} in the
case where $b_1(Y)=1$ and Theorem~\ref{thm:Euler} when $b_1(Y)>1$):

\begin{theorem}
\label{intro:Euler}
Let $Y$ be a three-manifold with $b_1(Y)>0$, and $\spinc$ be a
non-torsion $\SpinC$ structure, then
$$\chi(\HFp(Y,\spinc))=\pm\tau(Y,\spinc),$$
where $\tau \colon \SpinC(Y)\longrightarrow \Z$ is Turaev's torsion
function. In the case where $b_1(Y)=1$, $\tau(\spinc)$ is calculated
in the ``chamber'' containing $c_1(\spinc)$.
\end{theorem}

For zero-surgery on a knot, there is a well-known formula for the
Turaev torsion in terms of the Alexander polynomial,
see~\cite{Turaev}. With this, the above theorem has the following
corollary (a more precise version of which is given in
Proposition~\ref{prop:PreciseChi}, where the signs are clarified):

\begin{cor}
\label{cor:AlexPoly}
Let $Y_0$ be the three-manifold obtained by zero-surgery on a knot
$K\subset S^3$, and write its symmetrized Alexander polynomial as
$$\Delta_K=a_0+\sum_{i=1}^d a_i(T^i+T^{-i}).$$
Then, for each $i\neq 0$,
$$\chi(\HFp(Y_0,\spinc_0+iH))=\pm \sum_{j=1}^d j a_{|i|+j},$$
where $\spinc_0$ is the $\SpinC$ structure with trivial first Chern
class, and $H$ is a generator for $H^2(Y_0;\Z)$.
\end{cor}

Indeed, a variant of Theorem~\ref{intro:Euler} also holds in the case
where the first Chern class is torsion, except that in this case, the
homology must be appropriately truncated to obtain a finite Euler
characteristic (see Theorem~\ref{thm:TruncEuler}). Also, a similar
result holds for $\HFm(Y,\spinc)$, see Section~\ref{subsec:EulerHFm}.

As one might expect, these homology theories contain more information
than Turaev's torsion. This  can be seen, for instance, from their behaviour
under connected sums, which is studied in Section~\ref{sec:ConnectedSums}.
(Recall that if $Y_1$ and $Y_2$ are a pair of three-manifolds both with positive
first Betti number, then the Turaev torsion of their connected sum vanishes.)

We have the following result:

\begin{theorem}
\label{intro:ConnectedSumNonTriv}
Let $Y_1$ and $Y_2$ be a pair of oriented three-manifolds, and let $Y_1\#
Y_2$ denote their connected sum. A $\SpinC$ structure over $Y_1\# Y_2$
has non-trivial $\HFp$ if and only if it splits as a sum
$\spinc_1\#\spinc_2$ with $\SpinC$ structures $\spinc_i$ over $Y_i$ for $i=1,2$,
with the property that both groups $\HFp(Y_i,\spinc_i)$ are non-trivial.
\end{theorem}

More concretely, we have the following K\"unneth principle concerning
the behaviour of the invariants under connected sums.

\begin{theorem}
\label{intro:ConnectedSum}
Let $Y_1$ and $Y_2$ be a pair of three-manifolds, equipped with $\SpinC$ structures
$\spinc_1$ and $\spinc_2$ respectively. Then, we have identifications
\begin{eqnarray*}
\HFa(Y_1\#Y_2,\spinc_1\#\spinc_2)&\cong& H_*(\CFa(Y_1,\spinc_1)\otimes_\Z
\CFa(Y_2,\spinc_2)) \\
\HFm(Y_1\#Y_2,\spinc_1\#\spinc_2)&\cong& H_*(\CFm(Y_1,\spinc_1)\otimes_{\Z[U]}
\CFm(Y_2,\spinc_2)),
\end{eqnarray*}
where the chain complexes $\CFa(Y_i,\spinc_i)$ and
$\CFm(Y_i,\spinc_i)$  represent
$\HFa(Y_i,\spinc_i)$ and $\HFm(Y_i,\spinc_i)$ respectively.
\end{theorem}

In Section~\ref{sec:Adjunction}, we turn to a property which
underscores the close connection of the invariants with the the
minimal genus problem in three dimensions (which could alternatively
be stated in terms of Thurston's semi-norm,
c.f. Section~\ref{sec:Adjunction}):

\begin{theorem}
    \label{intro:Adjunction} Let $\EmbSurf \subset Y$ be an oriented,
    connected, embedded surface of genus $g(\EmbSurf)>0$ in an
    oriented three-manifold with $b_{1}(Y)>0$.  If $\spinc$ is a $\SpinC$
    structure for which $\HFp(Y,\spinc)\neq 0$, then $$\big|\langle
    c_{1}(\spinc),[\EmbSurf]\rangle\big| \leq 2g(\EmbSurf)-2.$$
\end{theorem}

In Section~\ref{sec:TwistedCoeffs}, we give a technical interlude,
wherein we give a variant of Floer homologies with $b_1(Y)>0$ with
``twisted coefficients.'' Once again, these are Floer homology groups
associated to a closed, oriented three-manifold $Y$ equipped with $\spinc\in\SpinC(Y)$,
but now, we have one more piece of input: a module $M$ over the
group-ring $\Z[H^1(Y;\Z)]$. This construction gives a collection of
Floer homology modules $\uHFinf(Y,\spinc,M)$, $\uHFpm(Y,\spinc,M)$,
and $\uHFa(Y,\spinc,M)$ which are modules over the ring
$\Z[U]\otimes_\Z \Z[H^1(Y;\Z)]$. In the case where $M$ is the trivial
$\Z[H^1(Y;\Z)]$-module $\Z$, this construction gives back the usual
``untwisted'' homology groups from~\cite{HolDisk}.

In Section~\ref{sec:Surgeries}, we give a very useful calculational
device for studying how $\HFp(Y)$ (and $\HFa(Y)$) changes as the
three-manifold undergoes surgeries: the surgery long exact sequence.
There are several variants of this result. The first one we give is
the following: suppose $Y$ is an integral homology three-sphere,
$K\subset Y$ be a knot, and let $Y_p(K)$ denote the three-manifold
obtained by surgery on the knot with integral framing $p$.  When
$p=0$, we let $\HFp(Y_0)$ denote $$
\HFp(Y_0)=\bigoplus_{\spinc\in\SpinC(Y_0)}\HFp(Y_0,\spinc), $$ thought
of as a $\Z[U]$ module with a relative $\Zmod{2}$ grading.

\begin{theorem}
\label{intro:Surgeries}
If $Y$ is an integral homology three-sphere, then there is a
an exact sequence of $\Z[U]$-modules
$$\begin{CD} ... @>>> \HFp(Y)@>>>
\HFp(Y_0)@>>> \HFp(Y_1) @>>> ...
\end{CD},
$$
where all maps respect the relative $\Zmod{2}$-relative gradings.
\end{theorem}

A more general version of the above theorem is given in
Section~\ref{sec:Surgeries} which relates $\HFp$ for an oriented
three-manifold $Y$ and the three-manifolds obtained by surgery on a
knot $K\subset Y$ with framing $h$, $Y_{h}$, and the three-manifold
obtained by surgery along $K$ with framing given by $h+m$ (where $m$
is the meridian of $K$ and $h\cm m = 1$),
c.f. Theorem~\ref{thm:GeneralSurgery}.  Other generalizations include:
the case of $1/q$ surgeries (Subsection~\ref{subsec:FracSurg}), the
case of integer surgeries (Subsection~\ref{subsec:IntSurg}), a version
using twisted coefficients (Subsection~\ref{subsec:SurgeriesTwisted}), 
and an analogous results for $\HFa$
(Subsection~\ref{subsec:SurgeriesHFa}).

In Section~\ref{sec:HFinfty}, we study $\HFinfty(Y,\spinc)$.  We
prove that if $b_1(Y)=0$, then for any $\SpinC$ structure $\spinc$,
$\HFinf(Y,\spinc)\cong \Z[U,U^{-1}]$. More generally, if the Betti
number if $b_1(Y)\leq 2$, $\HFinf$ is determined by $H_1(Y;\Z)$. This
is no longer the case when $b_1(Y)>2$ (see~\cite{HolDiskGraded}).
However, if we use totally twisted coefficients (i.e. twisting by
$\Z[H^1(Y;\Z)]$, thought of as a trivial $\Z[H^1(Y;\Z)]$-module), then
$\uHFinf(Y,\spinc)$ is always determined by $H_1(Y;\Z)$
(Theorem~\ref{thm:HFinfTwist}). This non-vanishing result allows us to
endow the Floer homologies with an absolute $\Zmod{2}$ grading, and
also a canonical isomorphism class of coherent orientation system.

We conclude with two applications.

\subsection{First application: complexity of three-manifolds and 
surgeries}

As described in~\cite{HolDisk}, there is a finite-dimensional theory
which can be extracted from $\HFp(Y)$, given by $$
\HFred(Y)=\HFp(Y)/\Image U^d, $$ where $d$ is any sufficiently large
integer. This can be used to define a numerical complexity of an
integral homology three-sphere $Y$: $$ N(Y)=\rk
\HFred(Y). $$ An easy calculation shows that $N(S^3)=0$ (c.f.
Proposition~\ref{prop:Lensspaces}).

Correspondingly, we define a  complexity of the symmetrized Alexander polynomial of a knot
$$\Delta_K(T)=a_0+\sum_{i=1}^d a_i(T^i+T^{-i})$$
by the following formula:
$$\|\Delta_K\|_{\circ}=\max(0,-\MT_0(K))+2 \sum_{i=1}^d \big| \MT_i(K) \big|, $$
where 
$$\MT_i(K)=  \sum_{j=1}^d j a_{|i|+j}.$$

As an application of the theory outlined above, we have the following:

\begin{theorem}
\label{thm:Complexity}
Let $K\subset Y$ be a knot in an integer homology three-sphere, and
$n>0$ be an integer, then 
$$n\cm \big\|\Delta_K\big\|_{\circ} \leq
N(Y)+N(Y_{1/n}), $$ where $\Delta_K$ is the Alexander polynomial
of the knot, and $Y_{1/n}$ is the three-manifold obtained by $1/n$
surgery on $Y$ along $K$.
\end{theorem}

This has the following immediate consequences:

\begin{cor}
\label{cor:ComplexitySurgeries}
If $N(Y)=0$ (for example, if $Y\cong S^3$), and the symmetrized
Alexander polynomial of $K$ has degree greater than one, then
$N(Y_{1/n})>0$; in particular, $Y_{1/n}$ is not homeomorphic to $S^3$.
\end{cor}

And also:

\begin{cor}
\label{cor:ComplexitySurgeriesII}
Let $Y$ and $Y'$ be a pair of integer homology three-spheres. Then
there is a constant $C=C(Y,Y')$ with the property that if $Y'$ can be obtained
from $Y$ by $\pm 1/n$-surgery on a knot $K\subset Y$ with $n>0$, then
$\|\Delta_K\|_{\circ}\leq C/n$.
\end{cor}

It is interesting to compare these results to analogous results
obtained using Casson's invariant. Apart from the case where the
degree of $\Delta_K$ is one, Corollary~\ref{cor:ComplexitySurgeries}
applies to a wider class of knots. On the other hand, at present,
$N(Y)$ does not give information about the fundamental group of $Y$.
There are generalizations of Theorem~\ref{thm:Complexity}
(and its corollaries) using an absolute grading on the
homology theories given in~\cite{HolDiskGraded}.

Corollary~\ref{cor:ComplexitySurgeries} should be compared with the
result of Gordon and Luecke which states that no non-trivial surgery
on a non-trivial knot in the three-sphere can give back the
three-sphere, see~\cite{GorLueckI}, \cite{GorLueckII}, see also
~\cite{CGLS}.

\subsection{Second application: bounding the number of
gradient trajectories}
\label{subsec:TrajectoryNumbers}
We give another application, to Morse theory over homology
three-spheres. 

Consider the following question. Fix an integral homology three-sphere
$Y$. Equip $Y$ with a self-indexing Morse function $f\colon
Y\longrightarrow \R$ with only one index zero critical point and one
index three critical point, and $g$ index one and two critical points.
Endowing $Y$ with a generic metric $\mu$, we then obtain a gradient
flow equation over $Y$, for which all the gradient flow-lines
connecting index one and two critical points are isolated. Let
$m(f,\mu)$ denote the number of $g$-tuples of disjoint gradient
flowlines connecting the index one and two critical points (note that
this is {\em not} a signed count).  Let $M(Y)$ denote the minimum of
$m(f,\mu)$, as $f$ varies over all such Morse functions and $\mu$
varies over all such (generic) Riemannian metrics. Of course, $M(Y)$
has an interpretation in terms of Heegaard diagrams: $M(Y)$ is the
minimum number of intersection points between the tori $\Ta$ and $\Tb$
for any Heegaard diagram $(\Sigma, \alphas,\betas)$ where the
attaching circles are in general position or, more concretely, the
minimum (again, over all Heegaard diagrams) of the quantity
$$m(\Sigma,\alphas,\betas)=\sum_{\sigma\in S_g}\left(\prod_{i=1}^g
\Big|\alpha_i\cap
\beta_{\sigma(i)}\Big|\right),$$ where $S_g$ is the symmetric group on $g$
letters and $|\alpha\cap\beta|$ is the number of intersection points
between curves $\alpha$ and $\beta$ in $\Sigma$.

We call this
quantity the {\em simultaneous trajectory number} of $Y$.  It is easy to see that
if $M(Y)=1$, then $Y$ is the three-sphere.  It is natural to consider
the following
\vskip.3cm
\noindent{\bf Problem:} if $Y$ is a three-manifold, find $M(Y)$.
\vskip.3cm

Since the complex $\CFa(Y)$ calculating $\HFa(Y)$ is generated by intersection points between $\Ta$ and $\Tb$,
it is easy to see that we have the
following:

\begin{theorem}
\label{thm:BoundIntPts}
If $Y$ is an integral homology three-sphere, then
$$\rk \HFa(Y) \leq M(Y).$$
\end{theorem}

Using this, the relationship between $\HFp(Y)$ and $\HFa(Y)$
(Proposition~\ref{prop:NonVanishHFa}), and a surgery sequence for
$\HFa$ analogous to Theorem~\ref{intro:Surgeries}
(Theorem~\ref{thm:GeneralSurgeryHFa}), we obtain the following result,
whose proof is given in Section~\ref{sec:Applications}:

\begin{theorem}
\label{thm:BoundIntPtsSurg}
Let $K\subset S^3$ be a knot, and let $Y_{1/n}$ be the three-manifold
obtained by $+1/n$-surgery on $K$, then $$M(Y)\geq 4 k +1,$$ where $k$ is
the number of positive integers $i$ for which $\MT_i(K)$ is non-zero.
\end{theorem}

\subsection{Relationship with gauge theory}

The close relationship between this theory and Seiberg-Witten theory
should be apparent.

For example, Conjecture~\ref{Conj} is closely related to the Atiyah-Floer
conjecture (see~\cite{AtiyahFloer}, see also \cite{AFSalamon},
\cite{DostSal}), a loose statement of which is the following. A
Heegaard decomposition of an integral homology three-sphere
$Y=U_{0}\cup_{\Sigma}U_{1}$ gives rise to a space $M$, the space of
$SU(2)$-representations of $\pi_{1}(\Sigma)$ modulo conjugation, and a
pair of half-dimensional subspaces $L_{0}$ and $L_{1}$ corresponding
to those representations of the fundamental group which extend over
$U_{0}$ and $U_{1}$ respectively. Away from the singularities of $M$
(corresponding to the Abelian representations), $M$ admits a natural
symplectic structure for which $L_{0}$ and $L_{1}$ are Lagrangian. The
Atiyah-Floer conjecture states that there is an isomorphism between
the associated Lagrangian Floer homology
$HF^{\mathrm{Lag}}(M;L_{0},L_1)$ and the instanton Floer homology
$HF^{\mathrm{Inst}}(Y)$ for the three-manifold $Y$,
$$
          HF^{\mathrm{Inst}}(Y)\cong HF^{\mathrm{Lag}}(M;L_{0},L_{1}). 
$$ Thus,
Conjecture~\ref{Conj} could be interpreted as an analogue of the
Atiyah-Floer conjecture for Seiberg-Witten-Floer homology.

Of course, this is only a conjecture. But aside from the calculations
of Sections~\ref{sec:Examples} and \ref{sec:SW}, the close connection
is also illustrated by several of the theorems, including the Euler
characteristic calculation, which has its natural analogue in
Seiberg-Witten theory (see~\cite{MengTaubes},
\cite{TuraevTwo}), and the adjunction inequalities, which exist in
both worlds (compare~\cite{Auckly} and~\cite{KMthurston}).

Two additional results presented in this paper -- the surgery exact
sequence and the algebraic structure of the Floer homology groups
which follow from the $\HFinf$ calculations -- have analogues in
Floer's instanton homology, and conjectural analogues in
Seiberg-Witten theory, with some partial results already established.
For instance, a surgery exact sequence (analogous to
Theorem~\ref{intro:Surgeries}) was established for instanton homology,
see~\cite{Floer}, \cite{BraamDonaldson}. Also, the algebraic structure
of ``Seiberg-Witten-Floer'' homology for three-manifolds with positive
first Betti number is still largely conjectural, but expected to match
with the structure of $\HFp$ in large degrees
(compare~\cite{KMIrvine}, \cite{Marcolli}, \cite{Embeddings}); see
also~\cite{AusBra} for some corresponding results in instanton
homology.

However, the geometric content of these homology theories, which gives
rise to bounds on the number of gradient trajectories
(Theorem~\ref{thm:BoundIntPts} and
Theorem~\ref{thm:BoundIntPtsSurg}) has, at present, no direct analogue
in Seiberg-Witten theory; but it is interesting to compare it with
Taubes' results connecting Seiberg-Witten theory over four-manifolds
with the theory of pseudo-holomorphic curves, see~\cite{TauSWGromov}.
For discussions on $S^1$-valued Morse theory and
Seiberg-Witten invariants, see~\cite{TaubesGeoSW} and
\cite{HutchingsLee}.

Gauge-theoretic invariants in three dimensions are closely related to
smooth four-manifold topology: Floer's instanton homology is linked to
Donaldson invariants, Seiberg-Witten-Floer homology should be the
counterpart to Seiberg-Witten invariants for four-manifolds. In fact,
there are four-manifold invariants related to the constructions
studied here. Manifestations of this four-dimensional picture can
already be found in the discussion on holomorphic triangles
(c.f. Section~\ref{HolDiskOne:sec:HolTriangles} of~\cite{HolDisk} and
Section~\ref{sec:Surgeries} of the present paper).  
These four-manifold invariants are presented
in~\cite{HolDiskThree}.

Although the link with Seiberg-Witten theory was our primary
motivation for finding the invariants, we emphasize that the
invariants studied here require no gauge theory to define and
calculate, only pseudo-holomorphic disks in the symmetric
product. Indeed, in many cases, such disks boil down to holomorphic
maps between domains in Riemann surfaces. Thus, we hope that these
invariants are accessible to a wider audience.

\section{Basic properties}
\label{sec:BasProp}

We collect here some properties of $\HFa$, $\HFp$, $\HFm$, and
$\HFinf$ which follow easily from the definitions.

\subsection{Finiteness properties}

Note that $\HFa$ and $\HFp$ distinguish certain $\SpinC$ structures on
$Y$ -- those for which the groups do not vanish. 

\begin{prop}
\label{prop:NonVanishHFa} 
For an oriented three-manifold $Y$ with $\SpinC$ structure $\spinc$,
$\HFa(Y,\spinc)$ is non-trivial if and only if $\HFp(Y,\spinc)$ is
non-trivial (for the same orientation system).
\end{prop}

\begin{proof}
This follows from the natural long exact sequence:
$$
\begin{CD}
... @>>>  \HFa(Y,\spinc) @>>> \HFp(Y,\spinc) @>{U}>> \HFp(Y,\spinc) @>>> ...
\end{CD}
$$
induced from the short exact sequence of chain complexes
$$
\begin{CD}
0@>>> \CFa(Y,\spinc) @>>> \CFp(Y,\spinc) @>{U}>> \CFp(Y,\spinc)@>>> 0.
\end{CD}
$$

Now, observe that $U$ is an isomorphism on $\HFp(Y,\spinc)$ if and only if
the latter group is trivial, since each element of $\HFp(Y,\spinc)$ is
annihilated by a sufficiently large power of $U$.
\end{proof}

\begin{remark}
Indeed, the above proposition holds when we use an arbitrary coefficient ring. In particular, 
the rank of $\HFp(Y,\spinc)$ is non-zero if and only if the rank of $\HFa(Y,\spinc)$ is non-zero.
\end{remark}

Moreover, there are finitely many such $\SpinC$ structures:

\begin{theorem}
\label{thm:FinitelyMany}
There are finitely many $\SpinC$ structures $\spinc$ for which
$\HF^+(Y,\spinc)$ is non-zero. The same holds for $\HFa(Y,\spinc)$.
\end{theorem}

\begin{proof}
Consider a Heegaard diagram which is weakly $\spinc$-admissible for
all $\SpinC$ structures (i.e. a diagram which is $\spinc_0$-admissible
Heegaard diagram, where $\spinc_0$ is any torsion $\SpinC$ structure,
c.f. Remark~\ref{HolDiskOne:rmk:AdmissibleAtOnce} and, of course,
Lemma~\ref{HolDiskOne:lemma:StronglyAdmissible} of~\cite{HolDisk}). This
diagram can be used to calculate $\HFp$ and $\HFa$ for all
$\SpinC$-structures simultaneously. But the tori $\Ta$ and $\Tb$ have
only finitely many intersection points, so there are only finitely
many $\SpinC$ structures for which the chain complexes
$\CFp(Y,\spinc)$ and $\CFa(Y,\spinc)$ are non-zero.
\end{proof}

\subsection{Conjugation and orientation reversal}

Recall that the set of $\SpinC$ structures comes equipped with a
natural involution, which we denote $\spinc\mapsto {\overline
\spinc}$: if $v$ is a non-vanishing vector field which represents 
$\spinc$, then $-v$ represents represents ${\overline \spinc}$.  
The homology groups are symmetric under
this involution:

\begin{theorem}
\label{thm:ConjInvar}
There are $\Z[U]\otimes_{\Z}\Wedge^*H_1(Y;\Z)/\Tors$-module 
isomorphisms identifications
$$
\HFpm(Y,\spinc)\cong \HFpm(Y,{\overline\spinc}),{\hskip.3in}
\HFinf(Y,\spinc)\cong \HFinf(Y,{\overline\spinc}),{\hskip.3in}
\HFa(Y,\spinc)\cong \HFa(Y,{\overline\spinc}),
$$
%	\begin{eqnarray*}
%	\HFm(Y,\spinc)&{\cong}&\HFm(Y,{\overline\spinc}) \\
%	\HFinf(Y,\spinc)&{\cong}&\HFinf(Y,{\overline\spinc}) \\
%	\HFp(Y,\spinc)&{\cong}&\HFp(Y,{\overline\spinc}) \\
%	\HFa(Y,\spinc)&{\cong}&\HFa(Y,{\overline\spinc}) \\
%	\end{eqnarray*}
\end{theorem}

\begin{proof}
Let $(\Sigma,\alphas,\betas,z)$ be a strongly $\spinc$-admissible
pointed Heegaard diagram for $Y$.  If we switch the roles of $\alphas$
and $\betas$, and reverse the orientation of $\Sigma$, then this
leaves the orientation of $Y$ unchanged. Of course, the set of
intersection points $\Ta\cap\Tb$ is unchanged, and indeed to each pair
of intersection points $\x,\y\in\Ta\cap \Tb$, for each
$\phi\in\pi_2(\x,\y)$, the moduli spaces of holomorphic disks
connecting $\x$ and $\y$ are identical for both sets of data.
However, switching the roles of the $\alphas$ and $\betas$ changes the
map from intersection points to $\SpinC$ structures.  If $f$ is a
Morse function compatible with the original data
$(\Sigma,\alphas,\betas,z)$, then $-f$ is compatible with the new data
$(-\Sigma,\betas,\alphas,z)$; thus, if $s_{z}(\x)$ is the $\SpinC$
structure associated to an intersection point $\x\in\Ta\cap
\Tb$ with respect to the original data, then ${\overline {s_{z}(\x)}}$ 
is the $\SpinC$ structure associated to the new data. (Note also that
the new Heegaard diagram is strongly ${\overline\spinc}$-admissible.)
This proves the result.
\end{proof}

Of course, the Floer complexes give rise to cohomology theories as
well.  To draw attention to the distinction between the cohomology and
the homology, it is convenient to adopt conventions from algebraic
topology, letting $\HFa_{*}$, $\HFp_{*}$, and $\HFm_{*}$ denote the
Floer homologies defined before, and $\HFaCoh(Y,\spinc)$,
$\HFpCoh(Y,\spinc)$, and $\HFmCoh(Y,\spinc)$ denote the homologies of
the dual complexes $\Hom(\CFa(Y,\spinc),\Z)$,
$\Hom(\CFp(Y,\spinc),\Z)$ and $\Hom(\CFm(Y,\spinc),\Z)$ respectively.

\begin{prop}
\label{prop:Duality}
Let $Y$ be a three-manifold with and $\spinc$ be a torsion $\SpinC$
structure. Then, 
there are natural isomorphisms:
\begin{eqnarray*}
    \HFaCoh(Y,\spinc)\cong \HFa_{*}(-Y,\spinc)
&{\text{and}}&
HF^*_\pm(Y,\spinc)\cong HF_*^\mp (-Y,\spinc),
\end{eqnarray*}
where $-Y$ denotes $Y$ with the opposite orientation.
\end{prop}

\begin{proof}
Changing the orientation of $Y$ is equivalent to reversing the
orientation of $\Sigma$. 
Thus, for each $\x,\y\in\Ta\cap\Tb$, and each class 
$\phi\in\pi_{2}(\x,\y)$, there is a natural identification
$$\ModFlow_{J_s}(\phi)\cong \ModFlow_{-J_s}(\phi'),$$
where $\phi'\in\pi_2(\y,\x)$ is the class with
$n_z(\phi')=n_z(\phi)$, obtained by pre-composing each holomorphic 
map by complex conjugation. This induces the stated isomorphisms in 
an obvious manner.
\end{proof}

\section{Sample calculations}
\label{sec:Examples}

In this section, we give some calculations for $\HFa$, $\HF^\pm$, and 
$\HFred$ for several families of three-manifolds.

\subsection{Genus one examples}
\label{subsec:GenusOne}

First we consider an easy case, where $Y$ is the lens space $L(p,q)$.
(Of course, this includes the case where $Y$ is a sphere).

We will introduce some shorthand. Let ${\mathcal T}^{\infty}$ denote
$\Z[U,U^{-1}]$, thought of as a graded $\Z[U]$-module, where the
grading of the element $U^{d}$ is $-2d$. We let ${\mathcal T}^-$
denote the submodule generated by all elements with grading $\leq -2$
(i.e. this is a free $\Z[U]$-module), and ${\mathcal T}^+$ denote the
quotient, given its naturally induced grading.

\begin{prop}
\label{prop:Lensspaces}
If $Y=L(p,q)$, then for each $\SpinC$ structure $\spinc$,
$$\HFa(Y,\spinc)=\Z,\hskip.3in
\HFm(Y,\spinc)\cong {\mathcal T}^-,
\hskip.3in
\HFinf(Y,\spinc)\cong {\mathcal T}^\infty,
\hskip.3in
\HFp(Y,\spinc)\cong {\mathcal T}^+.$$
Furthermore, $\HFred(Y,\spinc)=0$.
\end{prop}

\begin{proof}
Consider the genus one Heegaard splitting of $Y$. Here we can arrange
for $\alpha$ to meet $\beta$ in precisely $p$ points. Each
intersection point corresponds to a different $\SpinC$ structure, and,
of course, all boundary maps are trivial.
\end{proof}

Next, we turn to $S^1\times S^2$.  Consider the torus $\Sigma$ with a
homotopically non-trivial embedded curve $\alpha$, and an isotopic
translate $\beta$. The data $(\Sigma,\alpha,\beta)$ gives a Heegaard
diagram for $S^1\times S^2$.

We can choose the curves disjoint, dividing $\Sigma$ into a pair of
annuli. If the basepoint $z$ lies in one annulus, the other annulus
${\PerDom}$ is a periodic domain. Since there are no intersection
points, one might be tempted to think that the homology groups are
trivial; but this is not the case, as the Heegaard diagram is not
weakly admissible for $\spinc_0$, and also not strongly admissible for
any $\SpinC$ structure.

To make the diagram weakly admissible for the torsion $\SpinC$
structure $\spinc_0$, the periodic domain must have coefficients with
both signs. This can be arranged by introducing canceling pairs of
intersection points between $\alpha$ in $\beta$ (compare
Subsection~\ref{HolDiskOne:subsec:STwoTimesSOne}
of~\cite{HolDisk}). The simplest such case occurs when there is only
one pair of intersection points $x^+$ and $x^-$. There is now a pair
of (non-homotopic) holomorphic disks connecting $x^+$ and $x^-$ (both
with Maslov index one), showing at once that $$\begin{array}{cc}
\HFa(S^1\times S^2,\spinc_0)\cong H_*(S^1), & 
\HFinf(S^1\times S^2,\spinc_0)\cong
H_*(S^1)\otimes_\Z {\mathcal T}^\infty, \\
\HFp(S^1\times S^2,\spinc_0)\cong
H_*(S^1)\otimes_\Z {\mathcal T}^+
&
\HFm(S^1\times S^2,\spinc_0)\cong
H_*(S^1)\otimes_\Z {\mathcal T}^-.
\end{array}
$$ (We are free to choose here the orientation system so that the two disks
algebraically cancel;
but there are in fact two equivalence classes
orientation systems giving two different Floer homology groups, just as
there are two locally constant $\Z$ coefficient systems over $S^1$
giving two possible  homology groups.) Since the described
Heegaard decomposition is weakly admissible for all $\SpinC$
structures, and both intersection points represent $\spinc_0$, it
follows that $$\HFa(S^1\times S^2,
\spinc)=\HFp(S^1\times S^2,\spinc)=0$$ if $\spinc\neq \spinc_0$. 

To calculate the other homologies in non-torsion $\SpinC$ structures,
we must wind transverse to $\alpha$, and then push the basepoint $z$
across $\alpha$ some number of times, to achieve strong
admissibility. Indeed, it is straightforward to verify that if $h\in
H^2(S^1\times S^2)$ is a generator, then for $\spinc=\spinc_0+n\cm h$
with $n>0$,
$$\partial^\infty[x^+,i]=[x^-,i]-[x^-,i-n];$$
in particular, 
$$\HFm(S^2\times S^1,\spinc_0+nh)\cong 
\HFinfty(S^2\times S^1,\spinc_0+nh)\cong \Z[U]/(U^n - 1).$$

\subsection{Surgeries on the trefoil}

Next, we consider the three-manifold $Y$ which is obtained by $+n$
surgery on the left-handed trefoil, i.e. the $(2,3)$ torus knot, with $n>6$.

\begin{prop}
\label{prop:Trefoil}
Let $Y= Y_{1,n}$ denote the three-manifold obtained by $+n$ surgery on a
$(2,3)$ torus knot. Then, if $n>6$, there is a unique $\SpinC$
structure $\spinc_0$, with the following properties:
\begin{enumerate}
\item For all  $\spinc\neq \spinc_0$, the Floer theories are trivial, i.e.
$\HFa(Y,\spinc)\cong \Z$, $\HFp(Y,\spinc)\cong {\mathcal T}^+$, 
$\HFm(Y,\spinc)\cong {\mathcal T}^-$,
and 
$\HFred(Y,\spinc)=0$.
\item  $\HFa(Y,\spinc_0)$ is freely 
generated by three elements $a,b,c$ where $\gr(b,a)=\gr(b,c)=1$. 
\item $\HFp(Y,\spinc_0)$ is freely 
generated by elements $y$, and $x_i$ for $i\geq 0$, with 
$\gr(x_i,y)=2i$, $U_+(x_i)=x_{i-1}$, $U_+(x_0)=0$. 
\item  $\HFm(Y,\spinc_0)$ is freely 
generated by elements $y$, and $x_i$ for  $i < 0$, with 
$\gr(x_i, y)=2i+1$, $U_-(x_i)=x_{i-1}$.
\item $\HFred(Y,\spinc_0)\cong \Z$. 
\end{enumerate}
\end{prop}

\begin{figure}
\mbox{\vbox{\epsfbox{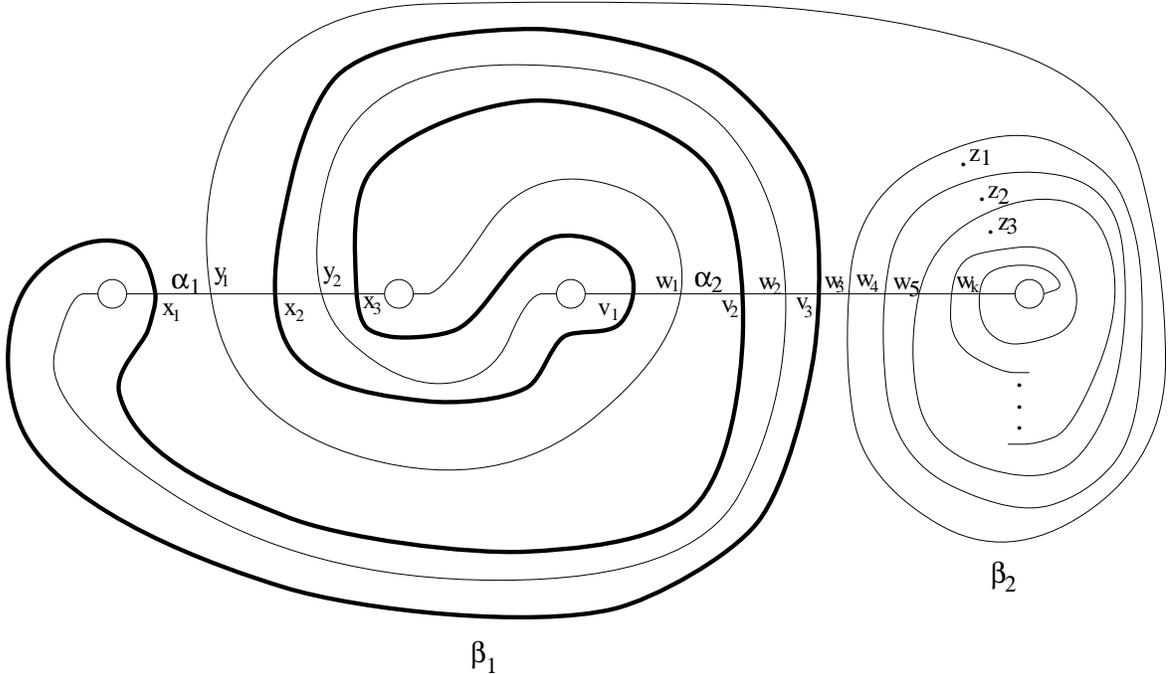}}}
\caption{\label{fig:TTwoThree}
Surgeries on the $(2,3)$ torus knot.}
\end{figure}

Before proving this proposition, we introduce some notation and 
several lemmas. 
For $Y$ we exhibit a genus 2 Heegaard decomposition and attaching
circles (see Figure~\ref{fig:TTwoThree}), where $k=n+6$, and  where the
spiral on the right hand side of the picture meets the horizontal
circle $k-2$ times. For a general discussion on constructing  Heegaard
decompositions from link diagrams see~\cite{GompfStipsicz}.

The picture is to be interpreted as follows. Attach a one-handle
connecting the two little circles on the left, and another one handle
connecting the two little circles on the right, to obtain a genus two
surface. Extend the horizontal arcs (labeled $\alpha_1$ and
$\alpha_2$) to go through the one-handles, to obtain the attaching
circles. Also extend $\beta_2$ to go through both of these one-handles
(without introducing new intersection points between $\beta_2$ and
$\alpha_i$). Note that here $\alpha _1$, $\alpha _2$, $\beta _1$ 
correspond to the left-handed trefoil: if we take the genus 2 handlebody determined
by $\alpha _1$, $\alpha _2$, and add a two-handle along $\beta _1$ then we get the complement
of the left-handed trefoil in $S^3$. Now varying $\beta _2$ corresponds to different surgeries
along the trefoil.

We have labeled  $\alpha_1\cap \beta_1=\{x_1,x_2,x_3\}$,
$\alpha_2\cap \beta_1=\{v_1, v_2, v_3\}$, $\alpha_1\cap
\beta_2=\{y_1,y_2\}$, and $\alpha_2\cap \beta_2=\{w_1,...,w_k\}$. 
Let us also fix basepoints $z_1,...,z_{k-2}$ labeled from outside to inside
in the spiral at the right side of the picture.
Since $H_1(Y_n;\Z)\cong \Z/n\Z$, the intersection points 
$\{x_i,w_j\}$, $\{v_i,y_j\}$ of $\Ta\cap \Tb$
can be partitioned into $n$ equivalence classes, c.f.
Subsection~\ref{HolDiskOne:subsec:SpinCStructures} of~\cite{HolDisk}. 
As $n$ increases
by 1 the number of intersection points in $\Ta\cap \Tb$ increases
by 3. We will use the following:

\begin{lemma}
For $n>6$ the points $\{x_1,w_9\}$, $\{x_2,w_8\}$, and $\{x_3,w_7\}$ are in the
same equivalence class, and all other intersection points are in 
different equivalence classes. By varying the base point $z$ among the
$\{z_5,...,z_{k-2}\}$,  we get the Floer homologies in all $\SpinC$
structures. 
\end{lemma}

\begin{proof}
From the picture, it is clear that (for some appropriate orientation
of $\{\alpha_1,\alpha_2\}$ and $\{\beta_1,\beta_2\}$) we have:
\begin{eqnarray*}
{[\alpha_1]}\cm {[\beta_1]} &=& -1 \\
{[\alpha_2]}\cm {[\beta_1]} &=& -1 \\
{[\alpha_1]}\cm {[\beta_2]} &=& 2 \\
{[\alpha_2]}\cm {[\beta_2]} &=& n+2.
\end{eqnarray*}
Thus, if $\{[\alpha_1],B_1,[\alpha_2],B_2\}$ is a standard symplectic 
basis for $H_1(\Sigma_2)$, then
\begin{eqnarray*}
{[\beta_1]}&\equiv & -B_1-B_2 \\
{[\beta_2]}& \equiv & 2B_1+(n+2)B_2
\end{eqnarray*}
in $H_1(\Sigma)/\langle [\alpha_1], [\alpha_2] \rangle$. It follows that $H_1(Y_n)\cong \Zmod{n}$ is generated by
$B_1=-B_2=h$. 

We can calculate, for example, $\epsilon(\{x_1,w_i\}, \{x_2,w_i\})$
as follows. We find a closed loop in $\Sigma_2$ which is composed of
one arc $a\subset \alpha_1$, and another in $b\subset \beta_1$
both of which connect $x_1$ and $x_2$. We then calculate the
intersection number $(a-b)\cap \alpha_1 = 0$, $(a-b)\cap \alpha_2 =
-1$. It follows that $a-b = h$ in $H_1(Y)$. So, $\epsilon(\{x_1,w_i\},
\{x_2,w_i\})=h$.

Proceeding in a similar manner, we calculate:
\begin{eqnarray*}
\epsilon(\{x_2,w_i\},\{x_3,w_i\})&=& h \\
\epsilon(\{y_1,v_i\},\{y_2,v_i\})&=& 3h \\
\epsilon(\{y_i,v_1\},\{y_i,v_2\})&=& -h \\
\epsilon(\{y_i,v_2\},\{y_i,v_3\})&=& -h \\
\epsilon(\{x_i,w_1\},\{x_i,w_2\})&=& h \\
\epsilon(\{x_i,w_2\},\{x_i,w_3\})&=& -2h \\
\epsilon(\{x_i,w_j\},\{x_i,w_{j+1}\})&=& h \\
\end{eqnarray*}
for $j=3,...,k-1$.
Finally,
$\epsilon(\{y_1,v_3\}, \{x_1,w_3\})=0$, as these intersections can be
connected by a square. 

It follows from this that the equivalence class containing
$\{x_1,w_9\}$ contains three intersection points:
$\{x_1,w_9\}$,$\{x_2,w_8\}$, and $\{x_3,w_7\}$.

Finally, note that
$s_{z_{i+1}}({\bf x})-s_{z_{i}}({\bf x})=\epsilon
\beta_2^*$, for some fixed $\epsilon=\pm 1$, according to
Lemma~\ref{HolDiskOne:lemma:VarySpinC} of~\cite{HolDisk}, 
and $\beta_2^*$ generates $H^2(Y;\Z)$,
according to the intersection numbers between the $\alpha_i$ and
$\beta_j$ calculated above.
\end{proof}

%  Export at 60%

We can identify certain flows as follows:

\begin{lemma}
\label{lemma:Flows}
For all $3\leq i\leq k-2 $ there is a
$\phi\in\pi_2(\{x_3,w_i\},\{x_2,w_{i+1}\})$ and  a\\ 
$\psi\in\pi_2(\{x_1,w_{i+2}\},\{x_2,w_{i+1}\})$
with $\Mas(\phi)=1=\Mas(\psi)$. Moreover, 
$$\#\UnparModFlow(\phi)=\#\UnparModFlow(\psi)=\pm 1.$$ Furthermore,
$n_{z_r}(\phi)=0$ for $r<i-2$, and $n_{z_r}(\phi)=1$ for
$r\geq i-2$. 
Also,
$n_{z_r}(\psi)=1$ for $r\leq i-2$, and $n_{z_r}(\psi)=0$ for
$r>i-2$. 
\end{lemma}

\begin{proof}
We draw the domains $\cald(\phi)$ and $\cald(\psi)$ belonging to
$\phi$ and $\psi$ in Figures~\ref{fig:phi} and
\ref{fig:psi} respectively, where the coefficients 
are equal to 1 in the shaded regions and 0 otherwise. Let $\delta _1$,
$\delta _2$ denote the part of $\alpha _2$, $\beta _2$ that lies in
the shaded region of ${\cald}(\phi)$.  Once again, we consider the constant almost-complex structure structure $J_s\equiv\Sym^{2}(\sj)$.

Suppose that $f$ is a holomorphic representative of $\phi$, i.e. $f\in 
\ModFlow(\phi)$, and
let  $\pi: F\longrightarrow\Strip$ 
denote the corresponding 2-fold branched covering of the disk (see Lemma~\ref{HolDiskOne:lemma:Correspondence} of~\cite{HolDisk}). Also let
${\widehat f}: F\longrightarrow \Sigma$ denote the corresponding holomorphic map to $\Sigma$.
Since $\cald (\phi)$ has only 0 and 1 coefficients, it follows that
$F$ is holomorphically identified with its image, 
which is topologically an annulus. 
This annulus is obtained by first choosing $\ell=1$ or $2$ and 
then cutting the shaded region along an interval
$I\subset \delta_\ell$ starting at $w_{i+1}$. 
Let $c\in [0,1)$ denote the length of this cut.
Note that by uniformization, we can identify the interior of
$F$ with a 
standard open annulus $\AnnOp(r)=\{z\in\C\big| r<|z|<1\}$ 
for some $0<r<1$ (where, of course, $r$ 
depends on the cut-length $c$ and direction $\ell=1$ or $2$).

In fact, given any $\ell=1,2$ and $c\in [0,1)$, 
we can consider the annular region $F$
obtained by cutting the region 
corresponding to $\phi$ in the direction $\delta_{\ell}$ with length 
$c$. Once again, 
we have a conformal identification $\Phi$ of the region $F\subset \Sigma$
with 
some standard annulus 
$\AnnOp(r)$,
whose inverse extends to the boundary to give
a map
$\Psi\colon \Ann(r) \longrightarrow \Sigma$.
For a given $\ell$ and $c$  let $a_1$, $a_2$, $b_1$, $b_2$ denote the arcs
in the boundary of the annulus which 
map to $\alpha _1$, $\alpha _2$, $\beta _1$, $\beta_2$ 
respectively, and
let $\angle(a_j)$, $\angle(b_j)$ denote angle spanned by these
arcs in the standard annulus $\Ann(r)$.
A branched covering over $\Strip$ as above 
corresponds to an involution 
$\tau: F\longrightarrow F$ which permutes the arcs:
$\tau (a_1)=a_2$, $\tau(b_1)=b_2$. Such  an involution
exists if and only if $\angle (a_1)=\angle(a_2)$ in which case
it is unique (see Lemma~\ref{HolDiskOne:lemma:Annuli} of~\cite{HolDisk}).
According to the generic perturbation theorem, if the curves are in generic position then 
these solutions are transversally cut out. It follows that
$\Mas(\phi)=1$. 

We argue that for $\ell=1$ and $c\rightarrow 1$ the angles converge to
$\angle (a_1)\rightarrow 0$, $\angle (a_2)\rightarrow 2\pi $.  To see
this, consider a map $\Theta\colon \CDisk\longrightarrow \Sigma$,
which induces a conformal identification between the interior of the
disk and the contractible region in $\Sigma$ corresponding to $\ell=1$
and $c=1$. One can see that the continuous extension of the composite
$\Phi_{c}\circ \Theta$, as a map from the disk to itself converges to
a constant map, for some constant on the boundary. (It is easy to
verify that the limit map carries the unit circle into
the unit circle, and has winding number zero about the origin, so it
must be constant.)
Thus, as $c\goesto 1$, both curves $a_{1}$
and $b_{2}$ converge to a point on the boundary of the disk, proving
the above claim.  In a similar way, for $\ell=2$ and $c\rightarrow 1$
the angles converge to $\angle (a_1)\rightarrow 2\pi $, $\angle
(a_2)\rightarrow 0$.

Now suppose that for $c=0$ we have $\angle(a_1)<\angle(a_2)$. Then the
signed sum of solutions with $\ell=1$ cuts is equal to zero, and the
signed sum of solutions with $\ell=2$ cuts is equal to $\pm
1$. Similarly if for $c=0$ we have $\angle(a_2)<\arg(a_1)$, then the
signed sum of solutions with $\ell=1$ cuts is equal to $\pm 1$, and
the signed sum of solutions with $\ell=2$ cuts is equal to zero. This
finishes the proof for $\phi$, and the case of $\psi$ is completely
analogous.

Although the domains $\phi$ and $\psi$ do not satisfy the
boundary-injectivity hypothesis in
Proposition~\ref{HolDiskOne:prop:MoveToriTransversality} of~\cite{HolDisk},
transversality can still be achieved by the same argument as in that
proposition.  For example, consider $\phi$, and suppose we cut along
$\ell=1$, so that the map $f$ induced by some holomorphic disk $u$ is
two-to-one along part of its boundary mapping to $\alpha_2$. Then, it
must map injectively to the $\beta$-curves so, for generic position of
those curves, the holomorphic map $u$ is cut out transversally.
Arguing similarly for the $\ell=1$ cut, we can arrange that the moduli
space $\ModFlow(\phi)$ is smooth. The same considerations ensure
transversality for $\psi$.

Note also that we have counted points in $\UnparModFlow(\phi)$ and
$\UnparModFlow(\psi)$, for the family $J_s\equiv \Sym^{2}(\sj)$, but
it follows easily that the same point-counts must hold for small
perturbations of this constant family.

%  Export at 60%

\begin{figure}
\mbox{\vbox{\epsfbox{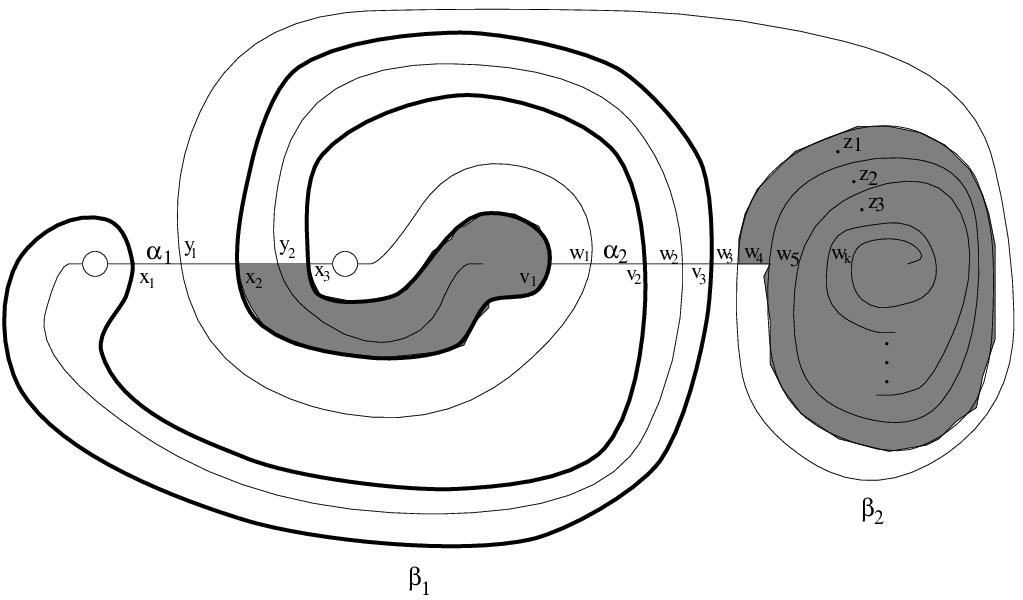}}}
\caption{\label{fig:phi}
Domain belonging to $\phi$ and $i=3$.}
\end{figure}

\begin{figure}
\mbox{\vbox{\epsfbox{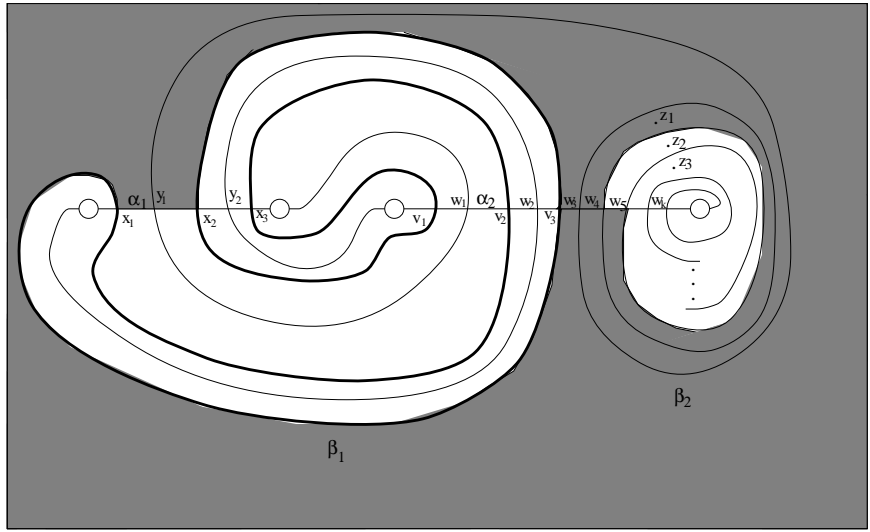}}}
\caption{\label{fig:psi}
Domain belonging to $\psi$ and $i=3$.}
\end{figure}

\end{proof}

\vskip.3cm
\noindent{\bf{Proof of Proposition~\ref{prop:Trefoil}.}}
Consider the equivalence class containing the elements
$\{x_1,w_9\}$, $\{x_2,w_8\}$, and $\{x_3,w_7\}$, denoted $a$, $b$, and
$c$ respectively. Let $\spinc _0$ denote the $\SpinC$ structure 
corresponding to this equivalence class and the basepoint $z_5$.
 According to Lemma~\ref{lemma:Flows}, in this $\SpinC$ structure
we have
$$\partial^\infty [a,j]=\pm [b,j-1], \ \ \ \partial^\infty [c,j]= 
\pm [b,j-1].$$
From the fact that $(\partial^\infty)^2=0$, it follows that
$\partial^\infty[b,j]=0$. The calculations for $\spinc_0$ follow.

Varying the basepoint $z_r$ with $r=6,...,k-2$, we capture all the
other $\SpinC$ structures. According to Lemma~\ref{lemma:Flows}, with
this choice, 
$$\partial^\infty [a,j]=\pm [b,j]  , \ \ \ \partial^\infty[c,j]= \pm[b,j-1]$$ 
This implies
the result for all the other $\SpinC$ structures.
\qed
\vskip0.3cm

More generally  let $Y_{m,n}$ denote the oriented 3-manifold
obtained by a $+n$ surgery along the torus knot $T_{2,2m+1}$. (Again we use the
left-handed versions of these knots, so for example $+1$ surgery would give the
Brieskorn sphere $\Sigma (2,2m+1, 4m+3)$).
In the following we 
will compute the Floer homologies of $Y_{m,n}$ for the case $n>6m$.

\begin{figure}
\mbox{\vbox{\epsfbox{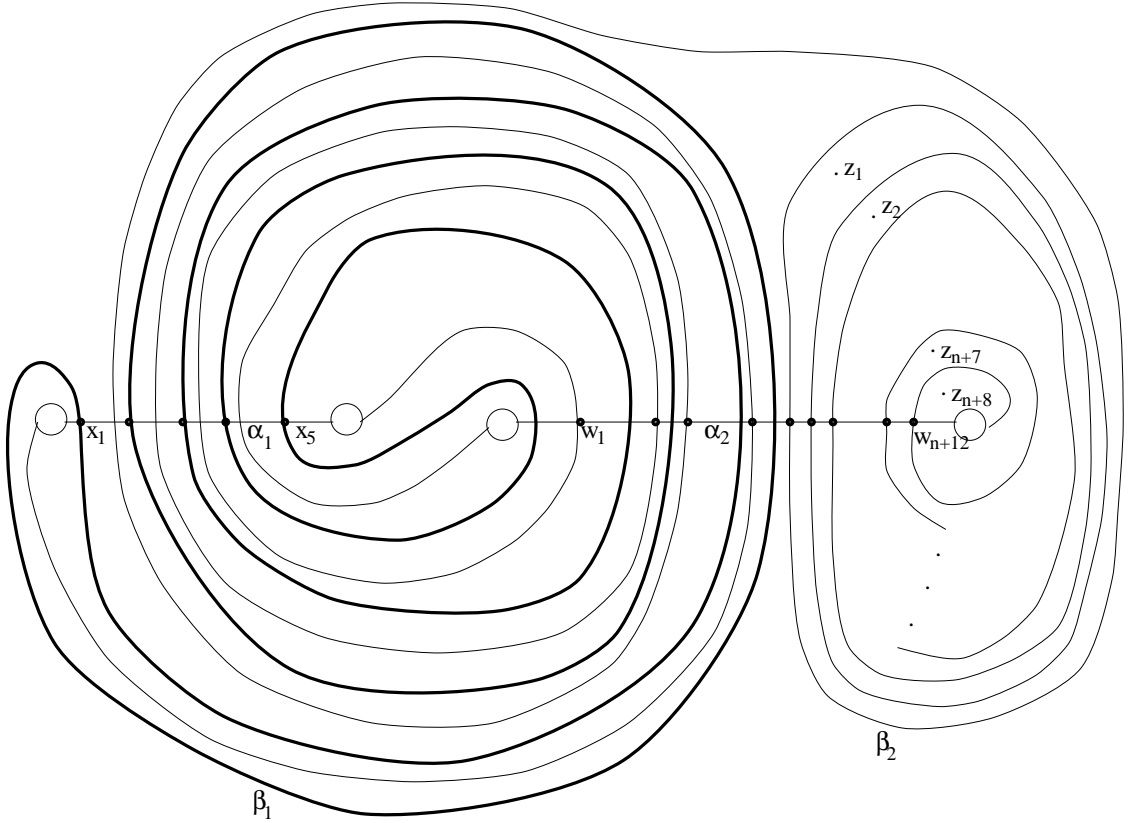}}}
\caption{\label{fig:kettoot}
$+n$ surgery on the $(2,5)$ torus knot.}
\end{figure}

First note that $Y_{m,n}$ admits a Heegaard decomposition of genus 2.
The corresponding picture is analogous to the $m=1$ case, except that
now $\beta _1$ and $\beta _2$ spiral more around $\alpha _1$,
$\alpha _2$, see  Figure \ref{fig:kettoot} for $m=2$. In general
 the $\beta _1$ curve hits both $\alpha _1$ and $\alpha _2$
in $2m+1$ points,  $\beta _2$ intersects $\alpha _1$ in $2m$ points and
$\alpha _2$ in   $n+6m$ points. Let $x_1,..., x_{2m+1}$ denote
the intersection points of $\alpha _1\cap \beta _1$, labeled from
left to right. Similarly let $w_1,...,w_{n+6m}$ denote
the intersection points of $\alpha _2\cap \beta _2$ labeled from
left to right. We also choose basepoints 
$z_1, z_2,...,z_{n+4m}$ in the spiral at the right hand side, labeled
from outside to inside.

\begin{lemma}
\label{lemma:BoundaryComp}
If $n>6m$, then there is an equivalence class
containing only the intersection
points  $a_i=\{x_i, w_{8m+2-i} \}$ for $i=1,...,2m+1$.
Furthermore if $\spinc _t$ denotes the $\SpinC$ structure determined
by this equivalence class and base point $z_{5m+t}$,
 for $ 1-m\leq t\leq  n-m$, then in this
$\SpinC$ structure we have
\begin{itemize}
\item  
$\partial _\infty [a_{2v+1},j]= \pm [a_{2v},j]\pm [a_{2v+2}, j-1]$,
for $t<m-2v$
\item 
$\partial _\infty [a_{2v+1},j]= \pm [a_{2v},j]\pm [a_{2v+2},j]$
for $t=m-2v$,
\item $\partial _\infty [a_{2v+1},j]= \pm [a_{2v},j-1]\pm [a_{2v+2}, j],$
for $t>m-2v$,
\end{itemize}
where $0\leq v\leq m$, and $a_0=a_{2m+2}=0$.
\end{lemma}

\begin{proof}
This is the same argument as in the $m=1$ case, together with the
observation that if $\phi\in \pi_2(a_{2v+1},a_{2\ell})$, and $\ell\neq
v$ or $v+1$, and $\Mas(\phi)=1$, then the domain $\cald(\phi)$
contains regions with negative coefficients (so the moduli space is
empty). Moreover, since $(\partial^\infty)^2=0$, it follows that $\partial^\infty([a_{2v},i])=0$.
\end{proof}

Note that $\spinc _{t+1}-\spinc _{t}\in H^2(Y_{m,n})$ is the Poincare dual of the meridian
of the knot. Since the meridian of the knot generates
$H_1(Y_{m,n})=\Z/n\Z$, it follows that $\{\spinc _t|\ 1-m\leq t\leq n-m\}= 
\SpinC (Y_{m,n})$, i.e we get all the $\SpinC$ structures this way.
 Now a straightforward computation gives the 
Floer homology groups of $Y_{m,n}$:

\begin{cor}
Let $Y=Y_{m,n}$ denote the three-manifold obtained by $+n$ surgery on the
$(2,2m+1)$ torus knot. Suppose that $n>6m$, and let $\spinc _t$ 
denote
the $\SpinC$ structures defined above.
For $m-1<t\leq n-m$ the Floer theories are trivial,
i.e  $\HFa (Y_{m,n}, \spinc _t)\cong \Z$, 
$\HFred (Y_{m,n}, \spinc _t)=0$, 
$\HFp (Y_{m,n}, \spinc _t)\cong {\mathcal T}^+$, and 
$\HFm (Y_{m,n}, \spinc _t)\cong {\mathcal T}^-$.
For $-m+1\leq t< 0$, the Floer homologies
of $\spinc _t$ are isomorphic to the corresponding Floer homologies
of $\spinc _{-t}$. 
Furthermore for $ 0\leq t\leq m-1$ we have
\begin{enumerate}
\item $\HFa (Y_{m,n},\spinc _t)$ is generated by $a, b, c$ with 
$\gr(b,a)= 1+2v_{m,t} +2t$, $\gr(b,c)=1+2v_{m,t}$.
\item $\HFp (Y_{m,n}, \spinc _t)$ is generated by $x_i$, $y_j$, for 
$0\leq i$, $0\leq j\leq v_{m,t}$,  $\gr(y_j, x_i)= 2(j-i+t)$ and
$U_+(x_i)= x_{i-1}$, $U_+(x_0)=0$, $U_+(y_i)=y_{i-1}$, $U_+(y_0)=0$.
\item $\HFm (Y_{m,n}, \spinc _t)$ is generated by $x_i$, $y_j$, for 
$i<0$, $0\leq j\leq v_{m,t}$, $\gr(y_j,x_i)= 2(j-i+t)-1$ and
$U_-(x_i)= x_{i-1}$, $U_-(y_i)=y_{i-1}$, $U_-(y_0)=0$.
\item $\HFred (Y_{m,n}, \spinc _t)$ is generated by  $y_j$, for 
 $0\leq j\leq v_{m,t}$,  $\gr(y_i, y_j)= 2i-2j$,
\end{enumerate}
where $v_{m,t}= \lfloor \frac{m-t-1}{2} \rfloor$, i.e. the greatest integer 
less than or equal to
$(m-t-1)/2$.
\end{cor}

\begin{remark}
\label{rmk:IdentifySpinStructure}
The symmetry of the Floer homology under the involution on the set of
$\SpinC$ structures ensures that $\spinc_0$ comes from a spin
structure. If $n$ is odd, there is a unique spin structure. With some
additional work one can show that, regardless of the parity of $n$,
$\spinc_0$ can be uniquely characterized as follows.  Let $X_{m,n}$ be
the four-manifold obtained by adding a two-handle to the four-ball
along the $(2,2m+1)$ torus knot with framing $+n$. Then, $\spinc_0$
extends to give a $\SpinC$ structure ${\mathfrak r}$ over $X_{m,n}$
with the property that $\langle c_1({\mathfrak r}), [S]\rangle = \pm
n$, where $S$ is a generator of $H_2(X_{m,n};\Z)$.  This calculation,
which is done in~\cite{HolDiskGraded}, follows easily from the
four-dimensional theory developed in~\cite{HolDiskThree}.
\end{remark}

In fact, Lemma \ref{lemma:BoundaryComp} can be used to prove that for
any $\SpinC$ structure on $Y_{m,n}$, $\HFinf(Y_{m,n},\spinc)\cong
{\mathcal T}^\infty$. Actually, it will be shown in Section~\ref{sec:HFinfty}
that for any rational homology three-sphere, $\HFinf(Y,\spinc)\cong {\mathcal T}^\infty$.

\section{Comparison with Seiberg-Witten theory}
\label{sec:SW}

\subsection{Equivariant Seiberg-Witten Floer homology}

We recall briefly the construction of equivariant Seiberg-Witten Floer
homologies $\HFto$, $\HFfrom$ and $\HFswred$. Our presentation here
follows the lectures of Kronheimer and Mrowka~\cite{KMIrvine}. For
more discussion, see~\cite{AusBra} for the instanton Floer homology
analogue, and also~\cite{Froyshov}, \cite{Marcolli}, \cite{Witten}.

Let $Y$ be an oriented rational homology 3-sphere, and $\spinc\in \SpinC(Y)$.
After fixing additional data (a Riemannian metric over $Y$ and some perturbation)
 the Seiberg-Witten equations over $Y$ in the $\SpinC$ structure $\spinc$
give a smooth moduli space consisting
of finitely many irreducible solutions $\gamma _1,..., \gamma _k$ 
and a smooth reducible solution $\theta$.

The chain-group $\CFto$ is freely generated by $\gamma _1,...,\gamma _k$
and $[\theta, i]$, for $i\geq 0$. Let $S$ denote this set of generators.
 The relative grading is given by 
$$\gr(\gamma _ j, [\theta,i])= {\rm dim}\left( 
\ModFlow(\gamma _j, \theta)\right) -2i,\ \ \gr(\gamma _j, \gamma_i )= {\rm dim}\left( 
\ModFlow(\gamma _j, \gamma _i)\right) $$
where $\ModFlow(\gamma _j,\theta)$ (resp. $\ModFlow(\gamma _j,\gamma _i)$)
 denotes the Seiberg-Witten
moduli space  of flows from $\gamma _j$ to $\theta$ (resp. $\gamma _j $ to $\gamma _i$).

\begin{defn}
For each $x,y \in S$ with $\gr(x,y)=1$ we define an incidence number
$c(x,y)\in \Z$, in the following way:
\begin{enumerate}
\item If $x=[\theta, i]$, then $c(x,y)=0$,
\item $c(\gamma _j, \gamma _i)= \#\upm(\gamma_j, \gamma _i)$,
\item  $c(\gamma _j, [\theta, 0])= \#\upm(\gamma_j, \theta)$
\item  $c(\gamma _j, [\theta, i])= \#(\upm(\gamma_j, \theta) \cap
\mu ({\rm pt})^i)$,
\end{enumerate}
where $\upm$ denotes the quotient of $\ModFlow$ by the $\R$ action of
translations, and $\cap \mu({\rm pt})^i$ denotes cutting down by a
geometric representative for $\mu({\rm pt})^i$ in a time-slice close
to $\theta$ (measured using the Chern-Simons-Dirac functional).  We
define the boundary map $\partial_{\rm to}$ on $\CFto$ by
$$\partial_{\rm to} (x)= \sum_{\{y\in S | \ \gr(x,y)=1\}} c(x,y)\cdot
y$$
\end{defn} 

It follows from the broken flowline compactification of two-dimensional
flows, modulo the $\R$ action, that $(\CFto, \partial _{\rm to})$ is a chain complex. Let $\HFto$ denote
the corresponding relative $\Z$ graded homology.

Similarly we can define the chain complex $(\CFfrom, \partial _{\rm from})$.
$\CFfrom$ is freely generated by $\gamma _1,...,\gamma _k$
and $[\theta, i]$, for $i \leq  0$. Let $S '$ denote this set of generators.
 The relative grading is determined 
by 
$$\gr([\theta,i], \gamma _j)= {\rm dim}\left( 
\ModFlow(\theta, \gamma _j)\right) +2i
,\ \ \gr(\gamma _j, \gamma_i )= {\rm dim}\left( 
\ModFlow(\gamma _j, \gamma _i)\right).
$$

\begin{defn}
For each $x,y \in S'$ with $\gr(x,y)=1$ we define an incidence number
$c'(x,y)\in \Z$, in the following way:
\begin{enumerate}
\item If $y=[\theta, i]$, then $c'(x,y)=0$,
\item $c'(\gamma _j, \gamma _i)= \#\upm(\gamma_j, \gamma _i)$,
\item  $c'( [\theta,0], \gamma _j)= \#\upm( \theta, \gamma _j)$
\item If $i < 0$, then  
 $c'([\theta, i], \gamma _j)= \#(\upm(\theta, \gamma _j) \cap
\mu ({\rm pt})^{-i})$.
\end{enumerate}
We define the boundary map $\partial_{\rm from}$ on $\CFfrom$ by
$$\partial_{\rm from} (x)= 
\sum _{\{y\in S' | \ \gr(x,y)=1\}} c'(x,y)\cdot y. $$
\end{defn} 
 
Again this gives a chain complex and we denote its homology by
$\HFfrom$. We also have a chain map 
$$f: \CFto\longrightarrow \CFfrom $$
given by $f(\gamma _j)= \gamma _j$, $f([\theta, i])=0$. Let $f_*$
denote the induced map between the Floer-homologies, and define
$$\HFswred= \HFto/({\rm Ker} f_*).$$

One reason to introduce these equivariant Floer-homologies
is that the irreducible Seiberg-Witten Floer homology (generated only by
$\gamma_1,...,\gamma _k$) is metric dependent. Analogy with
equivariant Morse theory suggests that the equivariant theories are
metric independent. Indeed the following was stated by Kronheimer and Mrowka,
\cite{KMIrvine}.

\begin{conj}
For oriented rational homology $3$-spheres $Y$ and $\SpinC$ structures
$\spinc \in \SpinC(Y)$ the 
equivariant Seiberg-Witten Floer homologies $\HFto(Y,\spinc)$,
$\HFfrom(Y,\spinc)$, and $\HFswred(Y, \spinc)$ are well-defined, 
i.e. they are independent of the
particular choice of metrics and perturbations.
\end{conj}

%\begin{remark}This conjecture is due to Kronheimer and Mrowka,
%who also outlined an argument to prove this. There are a
%few technical details like gluing to reducibles with obstructions.
%See also  \cite{Marcoli}. 
%\end{remark}

\subsection{Computations}

In this subsection we will compute $\HFto$, $\HFfrom$ and $\HFswred$
for the 3-manifolds studied in Section \ref{sec:Examples}, and for a
particular choice of perturbations of the Seiberg-Witten equations.
First, note that lens spaces all have trivial Seiberg-Witten Floer
homology, since they admit metrics with positive scalar curvature, in
particular, $\HFto(L(p,q),\spinc)$, $\HFfrom(L(p,q),\spinc)$ and
$\HFswred(L(p,q),\spinc)$ are isomorphic to ${\mathcal T}^+$, ${\mathcal T}^-$, and $0$
respectively. Note that all the 3-manifolds $Y=Y_{m,n}$ from
Section~\ref{sec:Examples} are Seifert-fibered so we can use
\cite{MOY} to compute their Seiberg-Witten Floer homology.

\begin{prop}
\label{prop:Comparison}
Let $Y=Y_{m,n}$ denote the oriented $3$-manifold obtained by $+n$ surgery
along the torus knot $T_{2, 2m+1}$. Suppose also that $n> 6m$.
Then for each $\spinc \in \SpinC(Y)$ we have
$$\HFto(Y,\spinc)\cong \HFp (Y,\spinc), \ \ \ \HFfrom(Y,\spinc)\cong \HFm(Y,\spinc), \ \ \ \HFswred (Y,\spinc)\cong 
\HFred (Y, \spinc),$$
where the isomorphisms are between relative $\Z$-graded Abelian groups,
and $\HFto (Y,\spinc)$, $\HFfrom (Y,\spinc)$, $\HFswred(Y, \spinc)$ are computed using a reducible connection
on the tangent bundle induced from the Seifert fibration of $Y$, and 
an additional  perturbation.
\end{prop}

\begin{proof}
First note  that
$Y_{m,n}$ is the boundary of the 4-manifold described by the plumbing diagram
in Figure 
\ref{fig:plumbing}, where the  number of $-2$ spheres in the right chain 
is $n+4m+1$. 
This gives a description of $Y_{m,n}$ as the
total space of an orbifold circle bundle over the sphere
with $3$ marked points with multiplicities  $2, 2m+1,k$ respectively,
where $k=n+4m+2$. The circle bundle $N$ has Seifert data
$$N=(-2,1,m+1,k-1).$$
and the canonical bundle is $K=(-2,1,2m,k-1)$.

\begin{figure}
\mbox{\vbox{\epsfbox{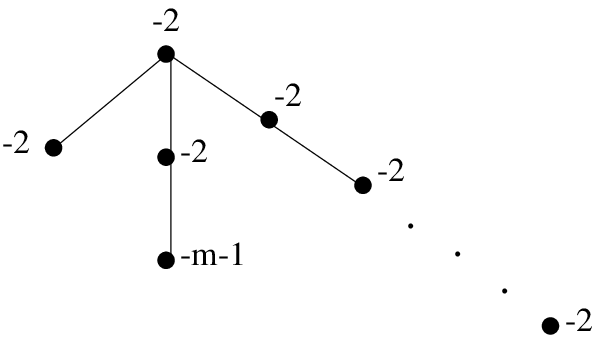}}}
{\caption{\label{fig:plumbing}}}
\end{figure}

Now we can apply \cite{MOY} to compute the irreducible solutions, relative gradings
and the boundary maps. 

Let us recall that for the unperturbed moduli space there is a 2 to 1 map from the set of irreducible solutions 
to the set of orbifold divisors $E$
with $E\geq 0$ and 
$$ {\rm deg}E <  {{\rm deg}(K) \over 2},$$
where the preimage consists of a holomorphic and an anti-holomorphic solution, that we denote
by $\cpl(E)$ and $\cmi(E)$ respectively. Note that $\cpl(E)$, $\cmi(E)$ lie in the $\SpinC$ structures determined
by the line-bundles $E$, $K\otimes E^{-1}$ respectively. 

In order to simplify the computation we will
use a certain perturbation of the Seiberg-Witten equation. Using the notation of \cite{HigherType} this
perturbation depends on a real parameter $u$, and corresponds  to adding a two-form
$iu(*d\eta)$ to the curvature equation, where $\eta$ is the connection form for $Y$ over the
orbifold. Now holomorphic solutions $\cpl(E)$ correspond to effective divisors
with
$$ {\rm deg}E <  {{\rm deg}(K) \over 2} -u{{\rm deg}(N)\over 2},$$ 
and anti-holomorphic solutions $\cmi(E)$ correspond to effective divisors
with
$$ {\rm deg}E <  {{\rm deg}(K) \over 2} +u{{\rm deg}(N)\over 2}.$$

 According to \cite{KMO}
the expected dimension of the moduli space between the reducible solution $\theta$ and
$\cplm(E)$ is computed by
$${\rm dim}\ModFlow(\theta, \cplm(E))= 1+2\left( \sum _{i\in I^\pm } \chi(E\otimes N^i)\right),$$
where $\chi (E\otimes N ^i)$ denotes the holomorphic Euler characteristic 
of the bundle $E\otimes N^i$, and $I ^\pm \subset \Z$ is given by the inequalities
$${\rm deg}E < {\rm deg}(E\otimes N^i)< {{\rm deg}(K)\over 2} \mp  u{{\rm deg}(N)\over 2} .$$

Returning to our examples let $E(a,b)$ denote the divisor $(0,0,a,b)$. It is easy to see that $\cmi (E(a,b))$
and $\cmi (E(a+1, b-2))$ are in the same $\SpinC$ structure. Also $\cmi(E(0,b))$ and
$\cpl (E(0, 2m-2-b))$ are in the same $\SpinC$ structure. From now
on let  $\spinc _0$ denote the $\SpinC$ structure given by  the line bundle
$E(0,m-1)$, and $\spinc _t$  corresponds to the
line-bundle $E(0,m-1+t)$. Clearly $\spinc _t\equiv \spinc _{t+n}$, since 
$H_1(Y,\Z)=\Z/n\Z$. 
 
Since 
$${\rm deg}E(a,b)={a\over {2m+1}}+ {b\over k},\ \ \ 
{\rm deg}K= {2m-1 \over 4m+2} -{1\over k},$$ 
for all $\spinc _t$ with  $n/4\leq |t| \leq n/2$ the unperturbed moduli space
(with $u=0$) have no irreducible solutions. It follows that
$\HFto(Y,\spinc _t)$ and $\HFfrom(Y, \spinc _t)$ are generated by $[\theta, i]$ and we have
the corresponding isomorphisms with ${\mathcal T}^+$, ${\mathcal T}^-$ respectively.

Clearly the $J$ action maps $\spinc _t$ to
$\spinc _{-t}$, so in the light of the $J$ symmetry in Seiberg-Witten theory,
it is enough to compute the equivariant Floer homologies for $0\leq t\leq n/4$. 
For these $\SpinC$ structures let us fix a perturbation with parameter $u$
satisfying 
$$ {\rm deg}(K) -u{\rm deg}(N) =-\epsilon,$$
where $\epsilon>0$ is sufficiently small. This perturbation eliminates all the
holomorphic solutions. It still remains to compute the anti-holomorphic solutions.

First let $0\leq t\leq m-1$. Since 
$${\rm deg}E(a,b)={a\over {2m+1}}+ {b\over k},\ \ \ 
{\rm deg}K= {2m-1 \over 4m+2} -{1\over k},$$ 
 the irreducible solutions in $\spinc _t$ are 
 $\delta  _r= \cmi (E(r, m-1-t-2r))$ for 
 $0\leq r\leq  {m-1-t \over 2}$. It is easy to see from 
\cite{MOY}, see also \cite{HigherType}, that the irreducible solutions and
$\theta$ are all transversally cut out
by the equations.

Computing the holomorphic Euler characteristic we get $\chi (E\otimes
N^{2i})= 1, $ for $0<2i\leq m-1-t -2r$, $\chi (E\otimes N^{2i+1})= -1,
$ for $m-1-t-2r<2i+1 \leq 2(m-r)-1$, and $\chi (E\otimes N^j)=0$ for
all other $j\in I ^-$, where $E=E(r, m-1-t-2r)$. The dimension formula
then gives $${\rm dim}\ModFlow(\theta , \delta _r)= -2t-2r-1.$$ As a
corollary we see that $\partial _{\rm from}$ is zero, since all these
moduli spaces have negative formal dimensions, and relative gradings
between the irreducible generators are even. In $\CFto$ the relative
gradings between all the generators are even, so $\partial _{\rm to}$
is trivial as well. Now the isomorphism between $\HFto(Y, \spinc _t)$
and $\HFp (Y, \spinc _t)$ corresponds to mapping $[\theta, i]$ to
$x_i$, and $\delta _r$ to $y_r$. Similarly the isomorphism between
$\HFfrom(Y, \spinc _t)$ and $\HFm (Y, \spinc _t)$ corresponds to
mapping $[\theta, i]$ to $x_{i-1}$, and $\delta _r$ to
$y_r$. Furthermore $\HFswred$ is freely generated by $\delta _r$ and
the map $\delta _r\rightarrow y_r$ gives the isomorphism with $\HFred$.

Now suppose that $m-1 <t \leq n/4$. Then there are no irreducible
solutions for the perturbed equation. So $\HFto$ and $\HFfrom$ are
generated by $[\theta, i]$ and we have the corresponding isomorphisms
with ${\mathcal T}^+$, ${\mathcal T}^-$ respectively.

For $ -n/4\leq t<0$ we get the analogous 
results by replacing $u$ with $-u$.
\end{proof}

\section{Euler characteristics}
\label{sec:EulerCharacteristic}

In this section, we analyze the Euler characteristics of the Floer
homology theories. In Subsection~\ref{subsec:EulerHFa}, we show that
the Euler characteristic of $\HFa$ is determined by $H_1(Y;\Z)$. After
that, we turn to the study of $\HFp$ for three-manifolds with $b_1>0$.

In~\cite{Turaev}, Turaev defines a torsion function
$$\tau_Y\colon \SpinC(Y)\longrightarrow \Z,$$ which is a
generalization of the Alexander polynomial.
This function can be calculated from a Heegaard diagram of $Y$
as follows.  Fix integers $i$ and $j$ between $1$ and $g$,
and  consider corresponding tori
\begin{eqnarray*}
\Ta^i = \alpha_1\times..\times {\widehat{\alpha_i}}\times...\times\alpha_g
&{\text{and}}&
\Tb^j=\beta_1\times..\times {\widehat{\beta_j}}\times...\times\beta_g
\end{eqnarray*}
in $\Sym^{g-1}(\Sigma)$
(where the hat denotes an omitted entry).  There is a map $\sigma$ from $\Ta^i\cap
\Tbj$ to $\SpinC(Y)$, which is given by thinking of each intersection
point as a $(g-1)$-tuple of connecting trajectories from index one to
index two critical points. Moreover, orienting $\alpha_i$, there is a
distinguished trajectory connecting the index zero critical point to
the index one critical point $a_i$ corresponding
to $\alpha_i$; similarly, orienting $\beta_j$,
there is a distinguished trajectory connecting the critical point
$b_j$ corresponding to the circle $\beta_j$ to the index index three
critical point in $Y$. This $(g+1)$-tuple of trajectories then gives
rise to a $\SpinC$ structure in the usual manner (modifying the
upward gradient flow in the neighborhoods of these
trajectories). Thus, we can define $$\Delta_{i,j}(\spinc)
=\pm \sum_{\{\x\in\Tai\cap\Tbj\big|
\sigma(\x)=\spinc\}} \epsilon(\x),$$
where $\epsilon(\x)$ is the local intersection number of $\Tai$ and
$\Tb^j$ at $\x$, and the overall sign depends on $i$, $j$ and $g$. (It
is straightforward to verify that this geometric interpretation is
equivalent to the more algebraic definition of $\Delta_{i,j}$ given
in~\cite{Turaev}, see for instance Section~\ref{Theta:sec:Alex}
from~\cite{Theta}.)

Choose $i$ and $j$ so that both $\alpha_i^*$ and $\beta_j^*$ have
non-zero image in $H^2(Y;\R)$.
When $b_1(Y)>1$, Turaev's torsion is characterized by the equation
\begin{equation}
\label{eq:TuraevsEquation}
\tau(\spinc)-\tau(\spinc+\alpha_i^*)-\tau(\spinc+\beta_j^*)+\tau(\spinc+\alpha_i^*+\beta_j^*)=
\Delta_{i,j}(\spinc),
\end{equation}
and the property that it has finite support.  (To define $\beta_j^*$
here, let $C$ be a curve in $\Sigma$ with $\beta_i\cap
C=\delta_{i,j}$, and let $\beta_j^*$ be Poincar\'e dual to the induced
homology class in $Y$.)  When $b_1(Y)=1$, we need a direction $t$ in
$H^2(Y;\R)$, which we think of as a component of $H^2(Y;\R)-0$.  Then,
$\tau_t$ is characterized by the above equation and the property that
$\tau_t$ has finite support amongst $\SpinC$ structures whose first
Chern class lies in the component of $t$.

For a three-manifold $Y$ with $\SpinC$ structure $\spinc$, the chain
complex $\CFp(Y,\spinc)$ can be viewed as a relatively
$\Zmod{2}$-graded complex (since the grading indeterminacy
$\divis(\spinc)$ is always even). Alternatively, this relative
$\Zmod{2}$ grading between $[\x,i]$ and $[\y,j]$ is calculated by
orienting $\Ta$ and $\Tb$, and letting the relative degree be given by
the product of the local intersection numbers of $\Ta$ and $\Tb$ at
$\x$ and $\y$. This relative $\Zmod{2}$-grading can be used to define
an Euler characteristic $\chi(\HFp(Y,\spinc))$ (when the homology
groups are finitely generated), which is well-defined up to an overall
sign.

In this section, we relate the Euler characteristics of
$\HFp(Y,\spinc)$ with Turaev's torsion function, when $c_1(\spinc)$ is
non-torsion. (The case where $c_1(\spinc)$ is torsion will be covered in
Subsection~\ref{subsec:TruncEuler}, after more is known about
$\HFinfty$; related results also hold for $\HFm$,
c.f. Subsection~\ref{subsec:EulerHFm}.)

The overall sign on $\chi(\HFp(Y,\spinc))$ will be pinned down once we define an absolute $\Zmod{2}$
grading on $\HFp(Y,\spinc)$ in Subsection~\ref{subsec:AbsoluteGradings}.

\subsection{Euler characteristic of $\HFa$}
\label{subsec:EulerHFa}

We first dispense with this simple object.

\begin{prop}
\label{prop:EulerHFa}
The Euler characteristic of $\HFa$ is given by
$$\chi(\HFa(Y,\spinc))
=\left\{\begin{array}{ll}
1 & {\text{if $b_1(Y)=0$}} \\
0 & {\text{if $b_1(Y)>0$}}
\end{array}\right..$$
\end{prop}

\begin{proof}
Both cases follow from the observation that $\chi(\HFa(Y,\spinc))$ is
independent of the $\SpinC$ structure $\spinc$. To see this, note that
for any $\beta_j$, we can wind normal to the $\alphas$ so that
$(\Sigma,\alphas,\betas,z)$ and $(\Sigma,\alphas,\betas,z')$ are both
weakly $\spinc$-admissible, where $z$ and $z'$ are two choices of
basepoint which can be connected by an arc which meets only
$\beta_j$. Now, both $\HFa(Y,\spinc)$ and
$\HFa(Y,\spinc+\PD[\beta_j^*])$ are calculated by the same equivalence
class of intersection points, using the basepoint $z$ in the first
case and $z'$ in the second. This changes only the boundary map, but
leaves the (finitely generated) chain groups unchanged, hence leaving
the Euler characteristic unchanged.

The result for $b_1(Y)>0$ then follows from this observation, together with
Theorem~\ref{thm:FinitelyMany}.

For the case where $b_1(Y)=0$, recall that the Heegaard decomposition
gives $Y$ a chain complex with $g$ one-dimensional generators
corresponding to the $\alphas$ (each of which is a cycle), and $g$
two-dimensional generators corresponding to the $\betas$. On the one
hand, the determinant of the boundary map is the order of the finite
group $H_1(Y;\Z)$ (which, in turn, is the number of distinct $\SpinC$
structures over $Y$); on the other hand, this determinant is easily
seen to agree with the intersection number
$\#(\Ta\cap\Tb)=\sum_{\spinc\in\SpinC(Y)}\chi(\HFa(Y,\spinc))$. The
result follows from this, together with $\spinc$-independence of
$\chi(\HFa(Y,\spinc))$.
\end{proof}

\subsection{$\chi(\HFp(Y,\spinc))$ when $b_1(Y)=1$ and $\spinc$ is non-torsion}
\label{subsec:BOneHPlus}

Our aim is to prove the following:

\begin{theorem}
\label{thm:EulerOne}
Suppose $b_1(Y)=1$. 
If $\spinc$ is a non-torsion $\SpinC$ structure, then $\HFp(Y,\spinc)$ is finitely generated, and indeed,
$$\chi(\HFp(Y,\spinc))=\pm \tau_{t}(Y,\spinc),$$
where $\tau_{t}$ is Turaev's torsion function, with respect to the
component $t$ of $H^2(Y;\R)-0$ containing $c_1(\spinc)$.
\end{theorem}

As usual, the Euler characteristic appearing above can be thought of
as the Euler characteristic of $\HFp(Y,\spinc)$ as a $\Z$-module; or,
alternatively, we could consider $\HFp(Y,\spinc,\Field)$ with
coefficients in an arbitrary field $\Field$.

The proof of Theorem~\ref{thm:EulerOne} occupies the rest of the present
subsection. 

Let $\spinc$ be a non-torsion $\SpinC$ structure on $Y$. Let $H$ be
the generator of $H_2(Y;\Z)$ with the property that $$ \langle
c_1(\spinc),H \rangle < 0. $$ After handleslides, we can arrange that
the periodic domain $\PerDom$ corresponding to $H$ contains $\alpha_1$
with multiplicity one in its boundary.

Choose a curve $\gamma$ transverse to $\alpha_1$ and disjoint from all
other $\alpha_i$ for $i>1$, oriented so that $\alpha_1\cap \gamma =
+1$. (Note that $\PD[\gamma]=\alpha_1^*$.) This curve has the property, then, that
$$\langle PD[\gamma], H\rangle =-1.$$
Let $\Tc=\gamma\times
\alpha_2\times...\times\alpha_g$. Winding $\alpha_1$ $n$ times along
$\gamma$, we obtain a new $\alpha$-torus, which we denote $\Ta(n)$. 
For each intersection point $\x\in\Tc\cap\Tb$
we obtain $2n$ intersection points in $\Ta(n)\cap\Tb$
$$\x^\pm_1,...,\x^{\pm}_n,$$ which we order with decreasing distance to $\gamma$, with a
sign $\pm$ indicating which side of $\gamma$ they lie on ($-$
indicates left, $+$ indicates right). We call the points in
$\Ta(n)\cap\Tb$ {\em $\gamma$-induced}: equivalently, a $\gamma$-induced
intersection point between $\Ta(n)$ and $\Tb$ is a $g$-tuple of points
in $\Sigma$, one of which lies in the winding region about $\gamma$.
It is easy to see that $\x^+_i$ and $\x^-_i$ lie in the same
equivalence class: indeed, there is a canonical
flow-line (with Maslov index $1$) connecting each $\x^+_i$ to
$\x^-_i$. 
Thus, (for any choice of base-point $z$), 
\begin{eqnarray*}
\spinc_z(\x^+_i)-\spinc_z(\x^+_j)&=&(i-j)\PD(\gamma), \\
\spinc_z(\x^+_i)&=&\spinc_z(\x^-_i).
\end{eqnarray*}

% 75 % .eps
\begin{figure}
	\mbox{\vbox{\epsfbox{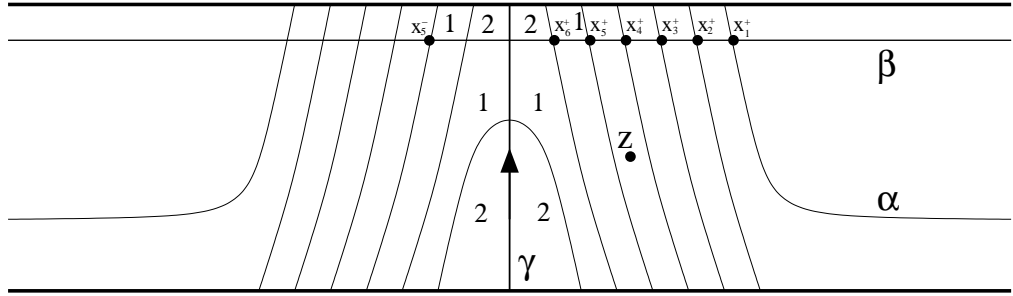}}} \caption{\label{fig:Wind2}
	{\bf Winding transverse to $\alpha$.}  We have pictured, once
	again, the cylindrical neighborhood of $\gamma$, and an
	$\alpha$-curve obtained by winding six times transverse to
	$\gamma$. The basepoint $z$ is placed in the third region, and
	intersection points corresponding to some $\beta$ are
	labeled. The multiplicities correspond to the domain of a flow
	connecting $x_5^+$ to $x_5^-$.}
\end{figure}

Our twisting will always be done in a ``sufficiently small'' area, so
that the area of each component of
$\Sigma-{\mathrm{nd}}(\gamma)-\alpha_1-\alpha_2-...-\alpha_g-\beta_1-...-\beta_g$ is
greater than $n$ times the area of ${\mathrm{nd}}(\gamma)$.

We will place our base-point $z$ to the right of $\gamma$, in the
$\left(\frac{n}{2}\right)^{th}$ subregion of the winding region about
$\gamma$. For this choice of basepoint, if $\x\in\Tc\cap\Tb$ then the
$\SpinC$ structure induced by $\x^\pm_{n/2}$ is independent of $n$.
Of course, the base-point is not uniquely determined by this
requirement: this region is divided into components by the
$\beta$-curves which intersect $\gamma$; but we fix any one such
region, for the time being.

\begin{lemma}
\label{lemma:PerClass}
If we wind $n$ times, and place the basepoint in the
$\left(\frac{n}{2}\right)^{th}$ subregion, and let $\PerDom_n$ denote the
corresponding periodic domain, then there is a constant $c$ with the property that we can find basepoints  
$w_1$ and $w_2$ (near $\gamma$ and away from $\gamma$ respectively), so that
\begin{eqnarray*}
n_{w_1}(\PerDom_n)\leq c-\frac{n}{2},
&{\text{and}}&
n_{w_2}(\PerDom_n)\geq c+\frac{n}{2}.
\end{eqnarray*}
\end{lemma}

\begin{lemma}
\label{lemma:InducedFromGamma}
Fix a $\SpinC$ structure $\spinc\in Y$. Then, 
if $n$ is sufficiently large, the $\gamma$-induced intersection points of 
$\Ta(n)\cap\Tb$ are the only ones which represent any of the $\SpinC$
structures of the form $\spinc+k\cm \PD[\gamma]$ for $k\geq 0$.
\end{lemma}

\begin{proof}
The intersection points between $\Ta(n)$ and $\Tb$ which are not
induced from $\gamma$ correspond to the intersection points
between the original $\Ta$ and $\Tb$. So, suppose that $\x$ is an
intersection point between $\Ta$ and $\Tb$ (there are,
of course, finitely many such intersection points), and let $z_0$ be
some basepoint outside the winding region.  As we wind $\alpha_1$ $n$
times, and place the new basepoint $z$ inside the winding region as
above (so as not to cross any additional $\beta$-curves), 
we see that $$s_z(\x)-s_{z_0}(\x)=-\frac{n}{2}\PD[\gamma],$$
where we think of $[\gamma]$ a one-dimension homology class in $Y$.
The lemma then follows.
\end{proof}

Let $(\Ta(n)\cap\Tb)^L\subset {\mathcal S}$ denote subset of $\gamma$-induced
intersection points where the $\alpha_1$ part lies to the ``left''
of $\gamma$, and $(\Ta(n)\cap
\Tb)^R$ denote subset of $\gamma$-induced intersection
points where the $\alpha_1$ part lies to the ``right'' of $\gamma$.  
(Note here that ${\mathcal S}$ denotes the subset of intersection
points which induce the given $\SpinC$ structure $\spinc$ over $Y$.)
There are corresponding subgroups
$\LeftFp$ and $\RightFp\subset \CFp(Y)$; similarly we have $\LeftFinf$ and $\RightFinf\subset\CFinf(Y)$. 

\begin{lemma}
\label{lemma:TwoHoClasses}
Fix $\spinc\in\SpinC(Y)$ and an integer $n$ sufficiently large (in comparison with
$\langle c_1(\spinc),\PerDom\rangle$). Then, for each $\gamma$-induced
pair $\x^+$ and $\y^-$ inducing $\spinc$, there are at most two
homotopy classes $\PhiIn, \PhiOut\in \pi_2(\x^+,\y^-)$ with Maslov
index one and with only non-negative multiplicities. Moreover, there
are no such classes in $\pi_2(\y^-,\x^+)$. 
\end{lemma}

\begin{proof}
Assume $\gr(\x^+,\y^-)$ is odd, and 
let $\PhiIn_{n}$ be the class with $\Mas(\PhiIn_n)=1$, and whose $\alpha_1$
boundary lies entirely inside the tubular neighborhood of $\gamma$. 
We claim that $\cald(\PhiIn_{n+2})$ is obtained from $\cald(\PhiIn_{n})$ by
winding only its $\alpha_1$-boundary (and hence leaving the domain
unchanged outside the winding region).  This follows from the fact
that the Maslov index is unchanged under totally real isotopies of the
boundary. 
It follows then that the multiplicities of $\PhiIn_n$ inside a
neighborhood of $\gamma$ grow like $n/2$. Recall that the
multiplicities of $\PerDom_n$ inside grow like $-n/2$, while outside
they grow like $n/2$.

Now, the set of all $\Mas=1$ homotopic classes connecting $\x^+$ to $\y^-$
is given by $$\PhiIn_n + k \left(\PerDom_n-
\frac{\langle c_1(\spinc), \PerDom
\rangle}{2} S\right).$$ If this class is to have non-negative
multiplicities, we must have that $k=0$ or $1$. 
This proves the assertion concerning classes from $\x^+$ to $\y^-$,
letting 
$\PhiOut_n=\PhiIn_n+\left(\PerDom_n-\frac{\langle
c_1(\spinc),\PerDom\rangle}{2} S\right)$.

Considering classes from $\y^-$ to $\x^+$, note that all $\Mas=1$
classes have the form $$(S-\PhiIn) + k\left(\PerDom_n - \frac{\langle
c_1(\spinc),\PerDom\rangle}{2} S\right).$$ When $k<0$, these classes have
negative multiplicities outside $\gamma$. When $k\geq 0$, these have negative
multiplicities inside the neighborhood of $\gamma$.
\end{proof}

\begin{prop}
\label{prop:LeftIsSubcomplex}
Given a $\SpinC$ structure $\spinc$ and an
$n$ sufficiently large, the subgroup $\LeftFinf\subset \CFinf(Y,\spinc)$ is a
subcomplex. 
\end{prop}

\begin{proof}
This follows immediately from the previous lemma. 
\end{proof}

Of course, the above proposition allows us to think of $\RightFinf$ as a chain
complex, as well, with differential induced from the quotient structure
$\CFinf/\LeftFinf$. 

There is a natural map
$$\delta \colon \RightFinf \longrightarrow \LeftFinf $$ given by taking the
$\LeftFinf$-component of the boundary of each element in
$\RightFinf$. This induces the connecting homomorphism for the long exact
sequence associated to the short exact sequence of complexes:
$$ \begin{CD}
0 @>>> \LeftFinf @>>> \CFinf @>>> \RightFinf @>>> 0.
\end{CD}
$$

To understand the homomorphism $\delta$, 
let 
$$f_1\colon \RightFinf \longrightarrow \LeftFinf$$ be the homomorphism induced by
$f_1([\x^+_i,j])=[\x^-_i,j-n_z(\phi)]$,
where $\phi$ the disk connecting $\x^+_i$ to $\x^-_i$ which is
supported in the tubular neighborhood of $\gamma$. 

We can define an ordering on the $\gamma$-induced intersection points
representing $\spinc$ as follows. Let $[\x,i], [\y,j]\in {\mathfrak
S}\times \Z$, then there is a unique $\phi\in\pi_2(\x,\y)$
with $n_z(\phi)=i-j$ and $\partial(\cald(\phi))\cap \alpha_1$
supported inside the tubular neighborhood of $\gamma$. We denote the
class $\phi$ by $\phi_{[\x,i], [\y,j]}$.  We then say that
$$[\x,i]>[\y,j]$$ if $$\Mas(\phi_{[\x,i],[\y,j]})>0$$ or if
$$\Mas(\phi_{[\x,i],[\y,j]})=0$$ and the area
$\Area(\cald(\phi_{[\x,i],[\y,j]}))>0$.  Note that an ordering gives
us a partial ordering for elements in $\CFinf(Y,\spinc)$: fix $\xi,
\eta\in \CFp(Y,\spinc)$, we say that $\xi<\eta$ if each
$[\x,i]\in{\mathcal S}\times \Z$ which appears with non-zero
multiplicity in the expansion of $\xi$ is smaller than each $[\y,j]\in
{\mathcal S}\times \Z$ which appears with non-zero
multiplicity in the expansion of $\eta$.

In the following lemma, it is crucial to work with {\em negative}
$\SpinC$ structures, i.e. those for which $\langle c_1(\spinc),
\PerDom \rangle < 0$.

\begin{lemma}
\label{lemma:DecomposeDelta}
If $\spinc$ is a negative $\SpinC$ structure, then 
the map $$\delta \colon \RightFinf \longrightarrow \LeftFinf$$ can be written
as $$\delta = f_1 + f_2,$$ so that $$f_2(g)<f_1(g)$$ for
each $g=[\x,i]\in \RightFinf$.
\end{lemma}

\begin{proof}
Consider a pair of generators $[\x^+,i]$ and $[\y^-,j]$, for which the
coefficient of $\delta$ is non-zero, i.e. that gives a homotopy class 
$\psi$ for which $\Mas(\psi)=1$ and $\cald(\psi)\geq 0$. Thus,
by Lemma~\ref{lemma:TwoHoClasses}, there are two possible cases, where
$\psi=\PhiIn$ or $\psi=\PhiOut$
(for $\x^+$ and $\y^-$). Note also that $\PhiIn=\phi_{[\x^+,i][\y^-,j]}$.

The case where $\psi=\PhiIn$, has two subcases, according to whether
or not $[\y^-,j]=f_1([\x^+,i])$. 
If 
$[\y^-,j]=f_1([\x^+,i])$,
$\psi=\phi_{[\x^+,i]f_1([\x^+,i])}$, and it follows easily that
$\#\Mod(\psi)=1$. Since the periodic domains have both positive and negative coefficients, 
the $[\y^-,j]$ coefficient of $f_2[\x^+,i]$
must vanish. If $[\y^-,j]\neq f_1([\x^+,i])$, then the domain of
$\phi_{f_1([\x^+,i]),[\y^-,j]}$ must include some region outside the
neighborhood of $\gamma$. Moreover, since
$$\phi_{[\x^+,i],f_1([\x^+,i])}+ \phi_{f_1([\x^+,i]),[\y^-,j]} =
\psi,$$ we have that  $\Mas(\phi_{f_1([\x^+,i]),[\y^-,j]})=0$; but
since the support of the twisting region is sufficiently small, it
follows that
$$\Area(\phi_{f_1([\x^+,i]),[\y^-,j]}) > 0;$$
i.e. $f_1([\x^+,i])>[\y^-,j]$.

When $\psi=\PhiOut$, it is easy to see that 
$$\phi_{[\x^+,i],[\y^-,j]}= \PhiOut - \PerDom.$$
It follows that
$\Mas(\phi_{[\x^+,i],[\y^-,j]})=1-\langle
c_1(\spinc),H(\PerDom)\rangle$. Moreover, 
$$\phi_{[\x^+,i],f_1([\x^+,i])}+\phi_{f_1([\x^+,i]),[\y^-,j]}
= \phi_{[x^+,i],[\y^-,j]},$$ 
so $\Mas(\phi_{f_1([\x^+,i]),[\y^-,j]})=-\langle
c_1(\spinc),H(\PerDom)\rangle>0$, by our hypothesis on $\spinc$, so that $f_1([\x^+,i])>[\y^-,j]$.
\end{proof}

\begin{prop}
\label{prop:DeltaSurj}
For negative $\SpinC$ structures $\spinc$, the map $\delta^+\colon
\RightFp\longrightarrow\LeftFp$ is surjective, and its kernel is
identified with the kernel of $f_1^+$ (as a $\Zmod{{\mathfrak
d}(\spinc)}$-graded groups).
\end{prop}

\begin{proof}
This is an algebraic consequence of Lemma~\ref{lemma:DecomposeDelta}.

We can define a right inverse to $f_1$,
$$P_1[\x_i^-, j] = [\x_i^+,j+n_z(\phi)],$$
where $\phi$ is the disk connecting $\x_i^+$ to $\x_i^-$. 
Then, we define a map
$$P= \sum_{N=0}^\infty P_1 \circ (-f_2\circ P_1)^{\circ N}.$$
Note that the right-hand-side makes sense, since the map $f_2\circ
P_1$ decreases the ordering (which is bounded below), so
for any fixed $\xi\in \RightFp$, there is some $N$ for which 
$$(-f_2\circ P_1)^{\circ N}(\xi)=0.$$
It is easy to verify that $P$ is a right inverse for $\delta^+$. 

The map sending $\xi \mapsto \xi-P\circ {\delta^+}(\xi)$ induces a
map from $\Ker f_1$ to $\Ker \delta^+$, which is injective, since for
any $\xi\in\Ker f_1$, we have that $$P\circ\delta^+(\xi)=P\circ 
f_2(\xi)<\xi.$$ Similarly, the map $\xi\mapsto \xi - P_1\circ f_1(\xi)$ supplies 
an injection $\Ker \delta^+ \longrightarrow \Ker f_1$. It follows 
that $\Ker f_1\cong \Ker \delta^+$. 
\end{proof}

\begin{prop}
\label{prop:IdentifyKerDelta}
For negative $\SpinC$ structures, the rank $\HFp(Y,\spinc)$ is finite. Moreover, we have that
$\chi(H_*(\ker \delta_\spinc^+))=\chi(\HFp(Y,\spinc))$.
\end{prop}

\begin{proof}
According to Proposition~\ref{prop:DeltaSurj} we have the short exact sequence
$$\begin{CD}
0@>>> \ker \delta^+ @>>> \RightFp @>{\delta^+}>> \LeftFp @>>> 0
\end{CD},$$
which we compare with the short exact sequence
$$\begin{CD}
0@>>>\LeftFp @>>> \CFp @>>> \RightFp @>>> 0
\end{CD}.$$ 
The result then follows by comparing the associated long exact
sequences, and observing that the connecting homomorphism for the second sequence agrees with the map on homology induced by $\delta^+$.
\end{proof}

\begin{prop} 
\label{prop:TuraevTorsion}
Let $\spinc$ be a negative $\SpinC$ structure, then 
$$\chi (\Ker f_1(\spinc)) = \pm\tau_{t}(\spinc),$$
where $t$ is the component of $H^2(Y,\Z)$ containing $c_1(\spinc)$.
\end{prop}

\begin{proof}
The map $f_1$ depends on a base-point and an equivalence class of
intersection point.  However, according to
Propositions~\ref{prop:DeltaSurj} and \ref{prop:IdentifyKerDelta},
$\chi(\Ker f_1^+(\spinc))$ depends on this data only through the
underlying $\SpinC$ structure $\spinc$ (when the latter is negative).
Let $\chi(\spinc)$ denote the Euler characteristic $\chi(\Ker
f_1|_{\spinc})$.  We fix a basepoint $z$ as before. We have a map
$$S_z\colon \Tc\cap\Tb \longrightarrow \SpinC(Y),$$ defined as
follows. Given $\x\in\Tc\cap \Tb$, we have
$$\spinc_z(\x_1^+)+(n_z(\phi)-1)\alpha_1^*,$$ where $\phi$ is the
canonical homotopy class connecting $\x_1^+$ and $\x_1^-$, and
$\alpha_1^*=\PD[\gamma]$.  (In fact, it is easy to see that the above
assignment is actually independent of the number of times we twist
$\alpha_1$ about $\gamma$.)  There is a naturally induced function
(depending on the basepoint) $$a_z\colon \SpinC(Y)\longrightarrow \Z$$
by $$a_z(\spinc)=\sum_{\{\x\in\Tc\cap\Tb\big|
S_z(\x)=\spinc\}}\epsilon(\x),$$ where $\epsilon(\x)$ is the local
intersection number of $\Tc\cap\Tb$ at $\x$.  It is clear that
$$\chi(\spinc)=\sum_{n=0}^\infty (n+1)\cdot a_z(\spinc+n \cm
\alpha_1^*).$$

It follows that 
\begin{equation}
\label{eq:FirstDifference}
\chi(\spinc)-\chi(\spinc+\alpha_1^*) = 
\sum_{n=0}^\infty a_z(\spinc+ n \cm \alpha_1^*).
\end{equation} 

We investigate the dependence of $a_z$ on the basepoint $z$. Note
first that there must be some curve $\beta_j$ which meets
$\gamma$ whose induced cohomology class $\beta_j^*$ is not a torsion
element in $H^2(Y;\Z)$: indeed, any $\beta_j$ appearing in the
expression $\partial \PerDom$ with non-zero multiplicity has this property.
Suppose that $z_1$ and $z_2$ are a pair of possible base-points
which can be connected by a path $z_t$ disjoint from all the attaching
circles except $\beta_j$, which it crosses transversally once, with $\#(\beta_j\cap z_t)=+1$.
We have a corresponding intersection point
$w\in \gamma\cap \beta_j$. We orient $\beta_j$ so that this
intersection number is negative (so that $\beta_j$ points in the same
direction as $\alpha_1$).

Now, we have two classes of intersection points
$\x\in\Tc\cap\Tb$: those which contain $w$ (each of these have the
form
$w\times \Ta^1\cap \Tb^j$), and those which do not. If $\x$ lies in
the first set, then 
$$S_{z_1}(\x)=S_{z_2}(\x)+\beta_j^*-\alpha_1^*;$$
if $\x$ lies in the second set, then 
$$S_{z_1}(\x)=S_{z_2}(\x)+\beta_j^*.$$
Note that there is an assignment:
$$\sigma'\colon \Ta^1\cap \Tb^j\longrightarrow \SpinC(Y)$$
obtained by restricting $S_{z_2}$ to 
$w\times (\Ta^1\cap\Tb^j)\subset \Tc\cap
\Tb$, and hence a
corresponding map
$$\Delta'\colon \SpinC(Y)\longrightarrow \Z.$$ We have the relation
that
\begin{equation}
\label{eq:ARelation}
a_{z_2}(\spinc)-a_{z_1}(\spinc+\beta_j^*)=
\Delta'(\spinc)-\Delta'(\spinc+\alpha_1^*).
\end{equation}

It follows from Equations~\eqref{eq:FirstDifference} 
and \eqref{eq:ARelation} that 
\begin{eqnarray*}
\chi(\spinc)-\chi(\spinc+\alpha_1^*)-\chi(\spinc+\beta_j^*)+\chi(\spinc+\alpha_1^*+\beta_j^*)
&=&
\sum_{n=0}^\infty
a_{z_2}(\spinc+n\alpha_1^*)-
a_{z_1}(\spinc+n\alpha_1^*+\beta_j^*) \\
&=&
\sum_{n=0}^\infty
\Delta'(\spinc+n\alpha_1^*)-\Delta'(\spinc+(n+1)\alpha_1^*)
\\
&=& \Delta'(\spinc).
\end{eqnarray*}
(note that $\Delta'$ has finite support). 

It is easy to see directly from the construction that $\Delta'$ and
the term $\Delta_{1,j}$ from Equation~\eqref{eq:TuraevsEquation} can
differ at most by a sign and a translation with $C_1\alpha_1^*+C_2\beta_j^*$, where
$C_1$ and $C_2$ are universal constants. Since $\tau(\spinc)$ and
$\chi(\HFp(Y,\spinc))$ are three-manifold invariants, by varying
$\beta_j^*$, it follows that $C_2=0$. A simple calculation in $S^1\times S^2$ shows that 
$C_1=0$, too. 
It follows that $\tau(\spinc)$ must agree with $\pm\chi(\HFp(Y,\spinc))$.
\end{proof}

\vskip.2cm
\noindent{\bf{Proof of Theorem~\ref{thm:EulerOne}.}}
This is now a direct consequence of Propositions~\ref{prop:DeltaSurj},
\ref{prop:IdentifyKerDelta} and
\ref{prop:TuraevTorsion}.
\qed
\vskip.2cm

\subsection{The Euler characteristic of $\HFp(Y,\spinc)$ when $b_1(Y)>1$, $\spinc$ is non-torsion}

\begin{theorem}
\label{thm:Euler}
If $\spinc$ is a non-torsion $\SpinC$ structure, over an oriented three-manifold
$Y$ with $b_1(Y)>1$,
then $\HFp(Y,\spinc)$ is finitely generated, and indeed,
$$\chi(\HFp(Y,\spinc))=\pm\tau(Y,\spinc),$$ where $\tau$ is Turaev's
torsion function.
\end{theorem}

The proof in subsection~\ref{subsec:BOneHPlus} applies, with the
following modifications. 

First of all, we use a Heegaard decomposition of $Y$ for which there
is a periodic domain $\PerDom$ containing $\alpha_1$ with multiplicity
one in its boundary, and with the property that the induced real
cohomology class $c_1(\spinc)$ is a non-zero multiple of
$\PD[\alpha_1^*]$. (This can be arranged after handleslides amongst
the $\alpha_i$.) The subgroup $c_1(\spinc)^\perp$ of $H_2(Y;\Z)$ which
pairs trivially with $c_1(\spinc)$ corresponds to the set of periodic
domains $\PerDom$ whose boundary contains $\alpha_1$ with multiplicity
zero. Let $\PerDom_2,...,\PerDom_{b}$ be a
basis for these domains. By winding normal to the
$\alpha_2,...,\alpha_g$, we can arrange for all of these periodic
domains to have both positive and negative coefficients with respect
to any possible choice of base-point on $\gamma$. It follows that the
Heegaard diagrams constructed above remain weakly admissible for any
$\SpinC$ structure. In the present case, the proof of
Lemma~\ref{lemma:TwoHoClasses} gives the following:

\begin{lemma}
\label{lemma:TwoHoClassesBig}
Fix $\spinc$ and an $n$ sufficiently large (in comparison with
$\langle c_1(\spinc),\PerDom\rangle$). Then, for each $\gamma$-induced
pair $\x^+$ and $\y^-$ inducing $\spinc$, there are at most two
homotopy classes modulo the action of $c_1(\spinc)^\perp$, 
$[\PhiIn], [\PhiOut]\in \pi_2(\x^+,\y^-)/c_1(\spinc)^\perp$ with Maslov
index one and with only non-negative multiplicities. Moreover, there
are no such classes in $\pi_2(\y^-,\x^+)$. 
\end{lemma}

Thus, Proposition~\ref{prop:LeftIsSubcomplex} holds in the present
context. In fact, the above lemma suffices to construct the
ordering. Note that there is no longer a unique map connecting $\x$ to
$\y$ with $\alpha_1$-boundary near $\gamma$, with specified
multiplicity at $z$ (the map $\phi_{[\x,i][\y,j]}$ from before), but
rather, any two such maps $\phi$ and $\phi'$ differ by the addition of
periodic domains in $c_1(\spinc)^\perp$. Thus, in view of
Theorem~\ref{HolDiskOne:thm:Grading} of~\cite{HolDisk}, the Maslov
indices of $\phi$ and $\phi'$ agree. If we choose the volume form on
$\Sigma$ so that all of $\PerDom_2,...,\PerDom_g$ have total signed
area zero (c.f. Lemma~\ref{HolDiskOne:lemma:EnergyZero}
of~\cite{HolDisk}), then the ordering defined by analogy with the
previous subsection is independent of the choice of $\phi$ or $\phi'$.

With these remarks in place, the proof of Theorem~\ref{thm:EulerOne}
applies, now proving that
$\chi(\spinc)=\pm \tau(\spinc)$, proving Theorem~\ref{thm:Euler}.

\section{Connected sums}
\label{sec:ConnectedSums}

In the second part of this section, we study the behaviour under
connected sums, as stated in Theorem~\ref{intro:ConnectedSum}. 
We begin with the simpler case of $\HFa$, and
then turn to $\HFm$.

\subsection{Connected sums and $\HFa$}
\label{subsec:HFaConnSum}

\begin{prop}
\label{prop:ConnSum}
Let $Y_1$ and $Y_2$ be a pair of oriented three-manifolds, and fix
$\spinc_1\in\SpinC(Y_1)$ and $\spinc_2\in\SpinC(Y_2)$. Let
$\CFa(Y_1,\spinc_1)$ and $\CFa(Y_2,\spinc_2)$ denote the corresponding
chain complexes for calculating $\HFa$. Then,
$$\CFa(Y_1\#Y_2,\spinc_1\#\spinc_2)\cong \CFa(Y_1,\spinc_1)\otimes_{\Z}
\CFa(Y_2,\spinc_2).$$
\end{prop}

In light of the universal coefficients theorem from algebraic
topology, the above result gives 
isomorphisms for all integers $k$:
\begin{eqnarray*}
\lefteqn{\HFa_{k}(Y_1\#Y_2,\spinc_1\#\spinc_2)} \\
&\cong&
\left(\bigoplus_{i+j=k} \HFa_i(Y_1,\spinc_1)\otimes
\HFa_j(Y_2,\spinc_2) \right)
\oplus\left(
\bigoplus_{i+j=k-1}\Tor(\HFa_i(Y_1,\spinc_1),\HFa_j(Y_2,\spinc_2))
\right)
\end{eqnarray*}
for some choice of absolute gradings on the complexes. (Of course,
this is slightly simpler with field coefficients, because in that case
all the $\Tor$ summands vanish.)

Note that Theorem~\ref{intro:ConnectedSumNonTriv} is an easy
consequence of this result, together with
Proposition~\ref{prop:NonVanishHFa}.

\vskip.2cm
\noindent{\bf{Proof of Proposition~\ref{prop:ConnSum}.}}
Fix weakly $\spinc_1$ and $\spinc_2$-admissible  pointed Heegaard diagrams 
$(\Sigma_1,\alphas,\betas,z)$ and $(\Sigma_2,\xis,\etas,z_2)$
for $Y_1$ and $Y_2$ respectively. Then, we form the pointed Heegaard
diagram $(\Sigma,\gammas,\deltas,z)$, where
$\Sigma$ is the connected sum of $\Sigma_1$ and $\Sigma_2$
at their distinguished points $z_1$ and $z_2$, $\gammas$ is the 
tuple of circles obtained by thinking of 
$\alphas\cup\xis$ as circles in $\Sigma$, and $\deltas$ are obtained
in the same way from $\betas\cup\etas$.
We place the basepoint $z$ in the connected sum region.
It is easy to see that $(\Sigma,\gammas,\deltas,z)$ is
represents $Y_1\#Y_2$. Moreover, there is an obvious identification
$$\Tc\cap\Td=(\Ta\cap\Tb)\times (\Tx\cap\Ty),$$
which is compatible with the relative gradings, in the sense that:
$$\gr(\x_1\times \x_2, \y_1\times \y_2)=\gr(\x_1,\y_1)+\gr(\x_2,\y_2).$$
Moreover, if $\phi\in \pi_2(\x_1\times \x_2, \y_1\times \y_2)$ has $n_z(\phi)=0$, then
$$\ModFlow_{J^{(1)}_s\# J^{(2)}_s}(\phi)\cong \ModFlow_{J^{(1)}_s}(\phi_1)
\times \ModFlow_{J^{(2)}_s}(\phi_2),$$ where
$\phi_i\in\pi_2(\x_i,\y_i)$ is the class with $n_{z_i}(\phi_i)=0$
(where $z_i\in\Sigma_i$ is the connected sum point), and $J^{(1)}_s$
and $J^{(2)}_s$ are families which are identified with
$\Sym^{(g)}(\sj_1)$ and $\Sym^{(g)}(\sj_2)$ near the connected sum
points, so we can form their connected sum
$J^{(1)}_s\#J^{(2)}_s$. Now, $\Mas(\phi)=1$ and $\ModFlow(\phi)$ is
non-empty, then the dimension count forces one of $\Mod(\phi_i)$ to be
constant. The proposition follows.\qed
\vskip.2cm

\subsection{Connected sums and $\HFm$}

We have seen how $\HFa$ behaves  under connected sum
(Proposition~\ref{prop:ConnSum}), and this suffices to give a
non-vanishing result for $\HFp$ under connected sums
(Theorem~\ref{intro:ConnectedSum}). The purpose of the present
subsection is to give a more precise description of the behaviour of
$\HFm$ and $\HFinf$ under connected sum. (Note that $\HFp$ can be
readily determined from $\HFm$ and $\HFinf$, using the long exact
sequence connecting these three $\Z[U]$-modules.)

Note that $\CFm(Y,{\mathfrak s})$, viewed as a $\Zmod{2}$-graded chain complex,
is finitely generated as a module over the ring $\Z[U]$.

\begin{theorem}
\label{thm:ConnSumHFm}
Let $Y_1$ and $Y_2$ be a pair of oriented three-manifolds, equipped
with $\SpinC$ structures ${\mathfrak s}_1$ and ${\mathfrak s}_2$ respectively. Then we have identifications:
\begin{eqnarray*}
\HFm(Y_1\# Y_2,{\mathfrak s}_1\#{\mathfrak s}_2)
&\cong& H_*\left(\CFm(Y_1,{\mathfrak s}_1)\otimes_{\Z[U]}
\CFm(Y_2,{\mathfrak s}_2)\right), \\
\HFinf(Y_1\# Y_2,{\mathfrak s}_1\#{\mathfrak s}_2)
&\cong& H_*\left(\CFinf(Y_1,{\mathfrak s}_1)\otimes_{\Z[U,U^{-1}]}
\CFinf(Y_2,{\mathfrak s}_2)\right). 
\end{eqnarray*}
\end{theorem}

Before proceeding with the proof of the above result, we give a
consequence for rational homology three-spheres $Y_1$ and $Y_2$, using
a field $\Field$ instead of the base ring $\Z$. In this case, since
$\HFm(Y,{\mathfrak s};\Field)$ is a finitely generated module over
$\Field[U]$, it splits as a direct sum of cyclic modules. Indeed, each
cyclic summand is either isomorphic to $\Field[U]$ or it has the form
$\Field[U]/U^{n}$ for some non-negative integer $n$, since if some
polynomial in $U$, $f(U)$, acts trivially on any element
$\xi\in\HFm(Y,{\mathfrak s})$, then clearly $U$ must divide $f$. We call
this exponent $n$ the {\em order} of the corresponding generator, i.e.
given a generator $\xi\in\HFm(Y,{\mathfrak s})$ as a $\Field[U]$-module, we
define its order $$\ord(\xi)=
\max\{i\in\Z^{\geq 0}\big| U^i\cm \xi\neq 0\}.$$ Note that by the
structure of $\HFinf(Y,{\mathfrak s})$, in any set of generators for
$\HFm(Y,{\mathfrak s})$ there is exactly one with infinite order.

\begin{cor}
Let $\Field$ be a field, and 
fix rational homology spheres $Y_1$ and $Y_2$. Let $\xi_i$ for
$i=0,...,M$ resp. $\eta_j$ for $j=0,...,N$ be generators of
$\HFm(Y_1,{\mathfrak s}_1;\Field)$ resp. $\HFm(Y_2,{\mathfrak s}_2;\Field)$ as a
$\Field[U]$-module. We order these so that
$\ord(\xi_0)=\ord(\eta_0)=+\infty$.  Then,
$\HFm(Y_1\#Y_2,{\mathfrak s}_1\#{\mathfrak s}_2;\Field)$ is generated as a
$\Field[U]$-module by generators $\xi_i\otimes \eta_j$ with $(i,j)\in
\{0,...,M\}\times \{0,...,N\}$ and also by generators
$\xi_i * \eta_j$ for $(i,j)\in \{1,...,M\}\times \{1,...N\}$. 
Moreover,  for all
$(i,j)\in \{0,...M\}\times \{0,...,N\}$, 
\begin{eqnarray*}
\ord (\xi_i\otimes
\eta_j)=\min (\ord(\xi_i),\ord(\eta_j))&{\text{and}}&
\gr (\xi_i\otimes\eta_j)=\gr(\xi_i)+\gr(\eta_j);
\end{eqnarray*}
while for all
$(i,j)\in \{1,...,M\}\times \{1,...,N\}$, we have that 
\begin{eqnarray*}
\ord (\xi_i * 
\eta_j)=\min (\ord(\xi_i),\ord(\eta_j)) &{\text{and}} &
\gr
(\xi_i * \eta_j)=\gr(\xi_i)+\gr(\eta_j)-1.
\end{eqnarray*}
In particular, we have that
$$\chi\left(\HFmRed(Y_1\# Y_2,{\mathfrak s}_1\#{\mathfrak s}_2)\right)=\chi 
\left(\HFmRed(Y_1,{\mathfrak s}_1)\right) + 
\chi\left(\HFmRed(Y_2,{\mathfrak s}_2)\right).$$
\end{cor}

\begin{proof}
This is an immediate application of Theorem~\ref{thm:ConnSumHFm} and
the K\"unneth formula for chain complexes over the principal ideal domain
$\Field[U]$. Specifically, we have that 
$$\HFm(Y_1\# Y_2,{\mathfrak s}_1\#{\mathfrak s}_2)
\cong \left(\HFm(Y_1,{\mathfrak s}_1)\otimes_{\Field[U]} \HFm(Y_1,{\mathfrak s}_2)\right)
\oplus \left(\HFm(Y_1,{\mathfrak s}_1) * \HFm(Y_1,{\mathfrak s}_2)\right),$$
where $A*B$ denotes the $\Tor$-complex, i.e.
$$(A*B)_{k}\cong \bigoplus_{i+j=k-1}\Tor_{\Field[U]}(A_i,B_j).$$
It is easy to see then that for any pair of non-negative integers $m$ and $n$,
$$(\Field[U]/U^m)\otimes_{\Field[U]} (\Field[U]/U^n)\cong\Field[U]/U^{\min(m,n)}\cong
\Tor_{\Field[U]}(\Field[U]/U^m, \Field[U]/U^n);$$
while for any $\Field[U]$-module $M$,
$\Field[U]\otimes_{\Field[U]} M \cong M$
and 
$\Tor_{\Field[U]}(\Field[U],M)=0$.

To see the Euler characteristic statement, we proceed as follows.
First, observe that to calculate the Euler characteristic of the 
graded $\Z$-module $\HFm(Y,{\mathfrak s})$ is the same as the Euler characteristic
of the $\Q$-vector space $\HFm(Y,{\mathfrak s};\Q)$.
From  above, we have that $\HFmRed(Y_1\#
Y_2,{\mathfrak s}_1\#{\mathfrak s}_2;\Q)$ is freely generated over $\Q$ by
$i,j\in\{0,...,M\}\times\{0,...,N\}-\{0,0\}$ with $U^m \xi_i\otimes
\eta_j$ where $m\in 0,...,\ord(\xi_i\otimes \eta_j)$ (observe that 
all generators of the form $U^m (\xi_0\otimes \eta_0)$ inject into
$\HFinf(Y_1\# Y_2,{\mathfrak s}_1\#{\mathfrak s}_2;\Field)$) and also 
generators $U^m (\xi_i * \eta_j)$ for $(i,j)\in\{1,...M\}\times\{1,...,N\}$ and $m\in
\{0,...,\ord(\xi_i*\eta_j)\}$. Observe in particular that when $i,j$
are both non-zero, $U^m(\xi_i \otimes \eta_j)$ has a corresponding element
$U^m (\xi_i*\eta_j)$ whose degree differs by one, so these cancel in the Euler characteristic. The only remaining elements are those of the form 
$U^m(\xi_i\otimes\eta_0)$ with $i>0$ and $m\in 0,..,\ord(\xi_i)$,
and also $U^n (\xi_0\otimes\eta_j)$ with $j>0$ and $n\in 0,...,\ord(\eta_j)$.
These contribute $\chi(\HFmRed(Y_1,\spinc_1))$ and $\chi(\HFmRed(Y_2,\spinc_2))$
to the Euler characteristic $\chi(\HFmRed(Y_1\# Y_2,\spinc_1\#\spinc_2))$
respectively.
\end{proof}

Before proving
Theorem~\ref{thm:ConnSumHFm}, we give the following special case.

\begin{prop}
\label{prop:ConnSumS1S2}
Let $\spinc_0$ be the $\SpinC$ structure on $S^2\times S^1$ with
$c_1(\spinc_0)=0$, and let $Y$ be  an oriented three-manifold, equipped with a $\SpinC$
structure $\spinc$.  There are 
isomorphisms:
\begin{eqnarray*}
\HFm(Y\#(S^2\times S^1),\spinc\#\spinc_0)&\cong &\HFm(Y,\spinc)\otimes
\wedge^* H^1(S^2\times S^1), \\
\HFinf(Y\#(S^2\times S^1),\spinc\#\spinc_0)&\cong &\HFinf(Y,\spinc)\otimes
\wedge^* H^1(S^2\times S^1), \\
\HFp(Y\#(S^2\times S^1),\spinc\#\spinc_0)&\cong &\HFp(Y,\spinc)\otimes
\wedge^* H^1(S^2\times S^1).
\end{eqnarray*}
For all other $\SpinC$ structures on $Y\#(S^2\times S^1)$, $\HFp$ vanishes.
\end{prop}

\begin{proof}
We consider first $\SpinC$ structures on $Y\#(S^2\times S^1)$ of the
form $\spinc\#\spinc_0$.  Let $(\Sigma,\alphas,\betas,z_1)$ be a
strongly $\spinc$-admissible pointed Heegaard diagram for
$Y$. Consider the Heegaard diagram for $S^2\times S^1$ discussed in
Section~\ref{subsec:GenusOne}, given by
$(E,\{\alpha_{g+1}\},\{\beta_{g+1}\},z_2)$, where $E$ is a genus one
surface and $\alpha_{g+1}$ and $\beta_{g+1}$ are a pair of exact
Hamiltonian isotopic curves meeting in a pair $x^+$ and $x^-$ of
intersection points. Choose the reference point $z_2$ so that the
exact Hamiltonian isotopy connecting the two attaching circles does
not cross $z_2$. Recall that there is a pair of
homotopy classes $\phi_1,\phi_2\in\pi_2(x^+,x^-)$ which contain
holomorphic representatives, indeed both containing a unique smooth,
holomorphic representative (for any constant complex structure on
$E$). We can form the connected sum diagram $(\Sigma\#
\alphas\cup\{\alpha_{g+1}\},\betas\cup\{\beta_{g+1}\},z)$, where we
form the connected sum along the two distingushed points, and let the
new reference point $z$ lie in the connected sum region.  This is
easily seen to be strongly $\spinc\#\spinc_0$-admissible.  Of course
$\Ta'\cap\Tb'=(\Ta\cap\Tb)\times\{x^+,x^-\}$; thus
$\CFp(Y_0,\spinc\#\spinc_0)$ is generated by $[\x,i]\otimes\{x^\pm\}$,
where $\x\in \Ta\cap\Tb$, and $\gr([\x,i]\otimes
\{x^+\},[\x,i]\otimes\{x^-\})=1$,
i.e. $\CFp(Y\#(S^2\times S^1),\spinc\#\spinc_0)\cong \CFp(Y,\spinc)\oplus
\CFp(Y,\spinc)$ (where the second factor is shifted in grading by one). 
We claim that when the neck is sufficiently
long, the differential respects this splitting.

Fix $\x,\y\in\Ta\cap\Tb$. 
First, we claim that for sufficiently long neck lengths, the only
homotopy classes $\phi'\in\pi_2(\x\times
\{x^+\},\y\times \{x^+\})$ 
with non-trivial holomorphic representatives are the ones which are
constant on $x^+$.  This follows from the following weak limit
argument. Suppose there is a homotopy class
$\phi'\in\pi_2(\{\x,x^+\},\{\y,x^+\})$ with $\Mas(\phi)\neq 0$ for
which the moduli space is non-empty for arbitrarily large connected
sum neck-length. Then, there is a limiting holomorphic disk in
$\Sym^g(\Sigma)\times E$. On the $E$ factor, the disk must be
constant, since $\pi_2(x^+,x^+)\cong \Z$ (here we are in the first
symmetric product of the genus one surface), and all non-constant
homotopy classes have domains with positive and negative coefficients.
Thus,  the limiting flow has the form $\phi\times \{x^+\}$ for some
$\phi\in\pi_2(\x,\y)$ (in $\Sym^g(\Sigma)$). 
Theorem~\ref{HolDiskOne:thm:Gluing} of~\cite{HolDisk} applies then to
give an identification $\ModFlow(\phi\times\{x^+\})\cong
\ModFlow(\phi')$.  Indeed, we have the same statement with $x^-$
replacing $x^+$.

Next, we claim
that (for generic choices) if $\phi'\in \pi_2(\x\times \{x^+\}, \y\times
\{x^-\})$ is any homotopy class with $\Mas(\phi')=1$, which contains a
holomorphic representative for arbitrarily long neck-lengths, then it
must be the case that $\x=\y$, and $\phi'=\{\x\}\times \phi_1$ or
$\phi'=\{\x\}\times \phi_2$. Again, this follows from weak limits.  If
it were not the case, we would be able to extract a sequence which
converges to a holomorphic disk in $\Sym^{g}(\Sigma)\times E$, which
has the form $\phi\times\phi_1$ or $\phi\times\phi_2$. Now, it is easy
to see that $\phi\times\{x^+\}*(\{\y\}\times \phi_i)=\phi'$ for $i=1$
or $2$ (by, say, looking at domains); hence,
$\Mas(\phi\times\{x^+\})=0$. It follows that as a flow in
$\Sym^g(\Sigma)$, $\Mas(\phi)=0$. Thus, there are generically no
non-trivial holomorphic representatives, unless $\phi$ is constant.
Observe, of course, that
$\#\UnparModFlow(\{\x\}\times\phi_1)=\#\UnparModFlow(\{\x\}\times\phi_2)=1$,
and also $n_z(\{\x\}\times\phi_1)=n_z(\{\x\}\times\phi_2)$.
With the appropriate orientation system, these flows cancel 
in the differential.

Putting these facts together, we have established that
$$\partial'([\x,i]\times\{x^\pm\})=(\partial [\x,i])\times \{x^\pm\}$$
(where $\partial'$ is the differential on
$\CFp(Y\#(S^2\times S^1),\spinc\#\spinc_0)$, and $\partial$ is the
differential on $\CFp(Y,\spinc)$. Indeed, it is easy to see the action
of the one-dimensional homology generator coming from $S^2\times S^1$
annihilates $[\x,i]\times\{x^-\}$, and sends $[\x,i]\times\{x^+\}$ to
$[\x,i]\times\{x^-\}$. 

When the first Chern class of the $\SpinC$ structure evaluates
non-trivially on the $S^2\times S^1$ factor, we can make $\alpha_{g+1}$
and $\beta_{g+1}$ disjoint, and have a Heegaard
diagram which is still weakly admissible for this $\SpinC$ structure. 
Since there are no intersection points, it follows that $\HFp$ in this case
is trivial.
\end{proof}

The proof of Theorem~\ref{thm:ConnSumHFm} is very similar to the proof
of Proposition~\ref{HolDiskOne:prop:Isomorphism} from~\cite{HolDisk}.
Like that proof, we find it convenient to subdivide the argument into
two cases depending on the first Betti number.

\vskip.2cm
\noindent{\bf{Proof of Theorem~\ref{thm:ConnSumHFm} when $b_1(Y_1\# Y_2)=0$}.}
First, we construct a chain map
$$\SumMap\colon \CFleq(Y_1,{\mathfrak s}_1)\otimes_{\Z[U]} \CFleq(Y_2,{\mathfrak s}_2)
\longrightarrow \CFleq(Y_1\# Y_2,{\mathfrak s}_1\#{\mathfrak s}_2).$$
To this end, consider pointed Heegaard diagrams
$(\Sigma_1,\alphas,\betas,z_1)$ and
$(\Sigma_2,\xis,\etas,z_2)$ for $Y_1$ and $Y_2$
respectively. Then there is connected sum Heegaard triple
$(\Sigma_1\#\Sigma_2,\alphas\cup \xis,
\betas\cup\xis,
\betas\cup\etas,z)$.
This triple describes a cobordism from
$Y_1\#(\#^{g_2}(S^2\times S^1))\coprod (\#^{g_1}(S^2\times S^1))\# Y_2$ to
$Y_1\# Y_2$ where where $g_1$ and $g_2$ are the genera of $\Sigma_1$
and $\Sigma_2$
respectively.
In fact, we let $\betas'$ and $\xis'$ be exact 
Hamiltonian translates of the $\betas$ and $\xis$ respectively,
so that the new triple
$$(\Sigma_1\#\Sigma_2,\alphas\cup \xis',
\betas\cup\xis,
\betas'\cup\etas,z),$$ is admissible. We let 
$\Theta_1\in\Tb\cap\TbPr$ and $\Theta_2\in\Tx\cap\TxPr$ denote the
``top'' intersection points in $\Sym^{g_1}(\Sigma_1)$
resp. $\Sym^{g_2}(\Sigma_2)$ between the tori corresponding to
$\betas$ and $\betas'$ resp. $\xis$ and $\xis'$.  In view of
Proposition~\ref{prop:ConnSumS1S2}, the maps $[\x,i]\mapsto [\x\times
\Theta_2,i]$ and $[\y,j]\mapsto [\Theta_1\times \y,j]$ give chain maps
$$\Phi_1\colon \CFleq(Y_1,{\mathfrak s}_1)\longrightarrow
\CFleq(Y_1\#^{g_2}(S^2\times S^1),{\mathfrak s}_1\#{\mathfrak s}_0)$$
and $$\Phi_2\colon \CFleq(Y_2,{\mathfrak s}_2)\longrightarrow
\CFleq(\#^{g_1}(S^2\times S^1) Y_2,{\mathfrak s}_0\#{\mathfrak s}_2)$$
are the chain maps considered in Proposition~\ref{prop:ConnSumS1S2}.
Now, we define $\SumMap$ to be the composite of $\Phi_1\otimes
\Phi_2$ with the map 
\begin{eqnarray*}
\lefteqn{F\colon \CFleq(Y_1\#\left(\#^{g_2}(S^2\times
S^1)\right),{\mathfrak s}_1\#\spinc_0) \otimes
\CFleq (\left(\#^{g_1}(S^2\times S^1)\right)\# Y_2,\spinc_0\#{\mathfrak s}_2) }
{\hskip3in}&& \\
&\longrightarrow&
\CFleq(Y_1\# Y_2,{\mathfrak s}_1\#{\mathfrak s}_2)
\end{eqnarray*}
defined by counting holomorphic triangles in the Heegaard triple
considered above.
Observe that
$F([\x,i-1]\otimes [\y,j])=F([\x,i]\otimes [\y,j-1])$, so that $F\circ(\Phi_1\otimes \Phi_2)$ is $\Z[U]$-bilinear, inducing the 
$\Z[U]$-equivariant chain map $\SumMap$.

Suppose that $\betas'$ is sufficiently close to the $\betas$. Then,
for each intersection point $\x\in\Ta\cap\Tb$, there is a unique
closest intersection point $\x'\in\Ta\cap\TbPr$; similarly, when $\xis'$ is
sufficiently close to $\xis$, each
intersection point $\y\in\Tx\cap\Ty$ corresponds to a unique 
closest intersection point $\y'\in\TxPr\cap\Ty$. In this case, 
there is an obvious map
$$\SumMap_0\colon \CFleq(Y_1,{\mathfrak s}_1)\otimes_{\Z[U]}\CFleq(Y_2,{\mathfrak s}_2)
\longrightarrow \CFleq(Y_1\# Y_2,{\mathfrak s}_1\#{\mathfrak s}_2)$$
defined by $$\SumMap_0([\x,i]\otimes[\y,j])=[\x'\times\y',i+j].$$ The map
$\psi_0$ is not necessarily a chain map, but it is clearly an
isomorphism of relatively $\Z$-graded groups.  Indeed, we claim that
when the total unsigned area $\epsilon$ in the regions between the
$\xi_i$ and the corresponding $\xi_i'$ (resp. $\beta_i$ and
corresponding $\beta_i'$) is sufficiently small, then, for the induced
energy filtration on (c.f. Section~\ref{HolDiskOne:sec:HandleSlides} 
of~\cite{HolDisk} and also Section~\ref{sec:Surgeries} below)
$\CFleq(Y_1\#Y_2,{\mathfrak s}_1\#{\mathfrak s}_2)$, we have
that $$\SumMap=\SumMap_0+{\text{lower order}}.$$ This is true because there
is an obvious small holomorphic triangle $\psi$ with $n_z(\psi)=0$,
$\Mas(\psi)=0$, and $\#\ModFlow(\psi)=1$ connecting
$\x\times\Theta_2$, $\Theta_1\times\y$, and $\x'\times\y'$. The total
area of this triangle is bounded by the total area $\epsilon$ (which
we can arrange to be smaller than any other triangle
$\psi'\in\pi_2(\x\times\Theta_2,\Theta_1\times\y, \w)$). Since the
energy filtration is bounded below in each degree (where now we view
the complexes as relatively $\Z$-graded modules over $\Z$), it follows
that $\Phi$ also induces an isomorphism in each degree. It follows that
$\SumMap$ induces an isomorphism of $\Z$-modules
$$\gamma\colon H_*\left(\CFleq(Y_1,{\mathfrak s}_1)
\otimes_{\Z[U]}\CFleq(Y_2,{\mathfrak s}_2)\right)
\longrightarrow \HFleq(Y_1\#Y_2,{\mathfrak s}_1\#{\mathfrak s}_2).$$

We have chosen to work with $\CFm$, but there is of course an 
identification $\CFleq\cong \CFm$ of complexes. Note also that the
above discussion also applies to prove the claim for $\CFinf$.
\qed
\vskip.3cm

For non-torsion $\SpinC$ structures $\spinc$, we must use the refined
filtration (again, as in Section~\ref{HolDiskOne:sec:HandleSlides}
of~\cite{HolDisk}). Specifically, given a strongly
$\spinc$-admissible Heegaard diagram, choose a volume form the
surface for which all $\spinc$-renormalized periodic domains have
total area zero.  Now, given $[\x,i]$ and $[\y,j]$ with the same
grading, we can find some disk $\phi\in \pi_2(\x,\y)$ with
$n_z(\phi)=i-j$ and $\Mas(\phi)=0$. We then define the filtration
difference to be the area of the domain associated to $\phi$:
$$\Filt([\x,i],[\y,j])=-\Area(\cald(\phi)).$$ Since any possible
choices of such disk $\phi$, $\phi'$ differ by a renormalized periodic
domain, it follows that the filtration defined above is independent of
the choice of disk.

Letting $\delta=\divis({\mathfrak s})$ be the grading indeterminacy of 
$\CFm(Y,\spinct)$, the filtration of $[\x,i]$ and $[\x,i+\delta]$
agree, since they can be connected by a Whitney disk $\phi$ whose
underlying domain is a renormalized periodic domain. Thus, the
filtration $\Filt$ is bounded below.

\vskip.2cm
\noindent{\bf{Proof of Theorem~\ref{thm:ConnSumHFm} when $b_1(Y_1\# Y_2)>0$}.}
When ${\mathfrak s}_1\#{\mathfrak s}_2$ is a torsion $\SpinC$
structure, the proof given under the assumption that $b_1(Y_1\# Y_2)=0$
adapts immediately in the present context. 

When ${\mathfrak s}_1\#{\mathfrak s}_2$ is non-torsion, we argue first
that the connected sum $Y_1\# Y_2$ can be endowed with a Heegaard
diagram which is both special in the above sense (each ${\mathfrak
s}_1\#{\mathfrak s}_2$-renormalied periodic domain has total area
zero), and it also splits as a sum of Heegaard diagrams
$(\Sigma_1\#\Sigma_2,\alphas\cup\xis,\betas\cup\etas,z)$. This is done
by winding the $\alphas$ within $\Sigma_1$, and the $\betas$ within
$\Sigma_2$. As in the proof of the theorem when $b_1(Y_1\# Y_2)=0$, we
consider the Heegaard triple $$(\Sigma_1\#\Sigma_2,\alphas\cup \xis',
\betas\cup\xis,
\betas'\cup\etas,z),$$ where $\xis'$ and $\betas'$ are obtained 
as sufficiently small Hamiltonian translates of the original $\xis$ and
$\betas$, letting $\epsilon$ denote the total (unsigned) areas in the regions
between the original curves and their Hamiltonian translates.

We claim that even when
${\mathfrak s}_1\#{\mathfrak s}_2$ is non-torsion, we can write
\begin{equation}
\label{eq:IsomorphismOfGamma}
\Gamma=\Gamma_0+{\text{lower order}},
\end{equation}
where now the lower order terms have lower order with respect to the
filtration $\Filt$ defined right before this proof.
To see this, suppose that $\psi$ is a holomorphic
triangle which contributes to $\Gamma$,
i.e. $\psi\in\pi_2(\x\times\y,\Theta_1\times\Theta_2,\p\times\q)$
satisfies $\Mas(\psi)=0$ and $\cald(\psi)>0$, while
$\psi_0\in\pi_2(\x\times\y,\Theta_1\times\Theta_2,\x'\times\y')$ is
the canonical small triangle. Assuming that $\x'\times\y'\neq
\p\times\q$, we argue that
$$\Filt([\x'\times\y',i],[\p\times\q,i-n_z(\psi)])<0.$$ To see this,
find some $\phi\in\pi_2(\x'\times\y',\p\times \q)$ with
$\Mas(\phi)=0$, so that both
$\psi,\psi_0*\phi\in\pi_2(\x\times\y,\Theta_1\times\Theta_2,\p\times\q)$
have $\Mas(\psi)=\Mas(\psi_0+\phi)=0$. Now, we claim that
$$\Area(\psi)=\Area(\psi_0+\phi),$$ since the difference is a
triply-periodic domain, while the $\xis'$ and $\etas'$ are obtained
from $\xis$ and $\etas$ by exact Hamiltonian translation. Since
$\Area(\psi)>\epsilon$, while $\Area(\psi_0)<\epsilon$, it follows
that $\Area(\phi)$ is positive.

Since the refined energy filtration is bounded below, the theorem now
follows as before.
\qed

\section{Adjunction Inequalities}
\label{sec:Adjunction}

\begin{theorem}
    \label{thm:Adjunction} 
    Let $\EmbSurf \subset Y$ be a connected
    embedded two-manifold of genus 
    $g(\EmbSurf)>0$ in an oriented three-manifold with $b_{1}(Y)>0$.
    If $\spinc$ is a $\SpinC$ structure for which
    $\HFp(Y,\spinc)\neq 0$, then 
    $$\big|\langle c_{1}(\spinc),[\EmbSurf]\rangle\big| \leq 2g(\EmbSurf)-2.$$
\end{theorem}

We can reformulate this result using Thurston's semi-norm, 
see~\cite{Thurston}. If 
$\EmbSurf=\bigcup_{i=1}^{k} \EmbSurf_{i}$ is 
a closed surface with $k$ connected
components, let
$$\chi_{-}(\EmbSurf)=\sum_{i=1}^{k} \max(0,-\chi(\EmbSurf_{i})).$$
The Thurston semi-norm of a homology class $\xi\in H_{2}(Y;\Z)$ is then 
defined by
$$\Theta(\xi)=\inf\{\chi_{-}(\EmbSurf)\big| \EmbSurf\subset Y, [\EmbSurf]=\xi\}.$$

In this language,  Theorem~\ref{thm:Adjunction} says the following:

\begin{cor}
If $\HFp(Y,\spinc)\neq 0$, then
$\big|\langle c_{1}(\spinc),\xi \rangle\big| \leq \Theta(\xi)$
for all $\xi\in H_{2}(Y;\Z)$.
\end{cor}

\begin{proof}
First observe that if $Z$ is an embedded sphere in $Y$, then for each
$\spinc$ for which $\HFp(Y,\spinc)\neq 0$, we have that $\langle
c_1(\spinc),[\EmbSurf]\rangle = 0$. This is a direct consequence of
Theorem~\ref{thm:Adjunction}:  attach a handle to $Z$ to get a
homologous torus $Z'$ and apply the theorem.

Now, let $\cup_{i=1}^k Z_i$ be a representative of $\xi$ whose
$\chi_-$ is minimal, labeled so that $Z_i$ for $i=1,...,\ell$ are the
components with genus zero. Then, $$ | \langle c_1(\spinc), \xi
\rangle | \leq
\sum_{i=1}^k |\langle c_1(\spinc), Z_i \rangle | 
\leq 
\sum_{i=\ell+1}^k (2g(Z_i)-2) \\
= \Theta(\xi).
$$
\end{proof}

Theorem~\ref{thm:Adjunction} is proved by constructing a special 
Heegaard diagram for $Y$, containing a periodic domain 
representative for $\EmbSurf$ with a particular form. The theorem then 
follows from a formula which calculates the evaluation of 
$c_{1}(\spinc)$ on $\EmbSurf$.

The following lemma, which is proved at the end of this subsection,
provides the required Heegaard diagram for $Y$.

\begin{lemma}
    \label{lemma:SpecialHeegaard} 
    Suppose $\EmbSurf\subset Y$ is a homologically non-trivial, 
    embedded two-manifold of genus $h=g(\EmbSurf)$, then $Y$ admits a
    genus $g$ Heegaard diagram $(\Sigma, \alphas, \betas)$, with
    $g>2h$, containing a periodic domain $\PerDom \subset \Sigma$
    representing $[\EmbSurf]$, all of whose multiplicities are one or
    zero.  Moreover, $\PerDom$ is a connected surface whose Euler
    characteristic is equal to $-2h$, and $\PerDom$ is bounded by $\beta_{1}$ and
    $\alpha_{2h+1}$.
\end{lemma}

Moreover, we have the following result, which follows from a more
general formula derived in Subsection~\ref{subsec:MaslovFormula}:

\begin{prop}
\label{prop:EasyMeasureCalc}
If $\x=\{x_1,...,x_g\}$ is an intersection point, and $z$ is chosen in
the complement of the periodic domain $\PerDom$ of Lemma~\ref{lemma:SpecialHeegaard},
then $$\langle c_1(\spinc_z(\x)), H(\PerDom)\rangle = 2-2h + 2
\#(\text{$x_i$ in the interior of $\PerDom$}).$$
\end{prop}

\vskip0.3cm
\noindent{\bf{Proof of Theorem~\ref{thm:Adjunction}.}}
If $\langle c_1(\spinc), [\EmbSurf]\rangle = 0$, then the inequality
is obviously true.

We assume that $\langle c_1(\spinc),[\EmbSurf]\rangle$ is non-zero.
If $\EmbSurf\subset Y$ is an embedded surface of genus
$g(\EmbSurf)=h$, then we consider a special Heegaard decomposition
constructed in Lemma~\ref{lemma:SpecialHeegaard}. Suppose that
$b_1(Y)=1$. Then this Heegaard decomposition is weakly admissible for
any non-torsion $\SpinC$ structure $\spinc$: there are no non-trivial
periodic domains $\cald$ with $\langle c_1(\spinc), H(\cald)
\rangle=0$.  Fix an intersection point $\x\in\Ta\cap\Tb$ which
represents $\spinc$. Clearly, of all $x_{i}\in\x$, exactly two must
lie on the boundary.  According to
Proposition~\ref{prop:EasyMeasureCalc}, then, $$\langle
c_{1}(\spinc),\PerDom \rangle = 2-2h +2\#(x_{i}\in \intPerDom);$$ i.e.
$$2-2h
\leq \langle c_{1}(\spinc),[\EmbSurf]\rangle.$$
If we consider the same inequality for $-\EmbSurf$ (or using the $J$ invariance), we get the stated bounds.

In the case where $b_1(Y)>1$, we must wind transverse to the
$\alpha_1,...,{\widehat {\alpha_{2h+1}}},...,
\alpha_g$ to achieve weak admissibility. Of course, 
we choose our transverse curves to be disjoint from one another (and
$\alpha_{2h+1})$. In winding along these curves, we leave the periodic
domain ${\mathcal P}$ representing $S$ unchanged.  Moreover, each
periodic domain ${\mathcal Q}$ which evaluates trivially on
$c_1(\spinc)$ must contain some $\alpha_j$ with $j\neq 2h+1$ on its
boundary; thus, by twisting sufficiently along the $\gamma$-curves, we
can arrange that the Heegaard decomposition is weakly admissible. The
previous argument when $b_1(Y)=1$ then applies.
\qed

We now return to the proof of Lemma~\ref{lemma:SpecialHeegaard}.

\vskip.3cm 
\noindent{\bf{Proof of Lemma~\ref{lemma:SpecialHeegaard}.}}
    The tubular neighborhood of $\EmbSurf$, identified with 
    $\EmbSurf\times [-1,1]$, has a handle decomposition with one
    zero-handle, $2h$ 
    one-handles, and one two-handle; i.e. the tubular neighborhood 
    admits a Morse 
    function $f$ with one index zero critical point $p$, $2h$ index one 
    critical points $\{a_{1},\ldots,a_{2h}\}$, and one index two 
    critical point $b_{1}$. Hence, we have a genus 
    $2h$ handlebody  $V_{2h}$, with an embedded circle on its boundary
    $\beta_{1}\subset \partial V_{2h}=\Sigma_{2h}$ (the descending 
    manifold of $b_{1}$). The circle 
    $\beta_{1}$ separates $\Sigma_{2h}$, and attaching a 
    two-handle to $V_{2h}$ along $\beta_{1}$ gives us the tubular neighborhood
    of $\EmbSurf$.
    Choose a component of the complement 
    of $\beta_{1}$, and denote its closure by
    $F_{2h}\subset \Sigma_{2h}$. Attaching the descending manifold of $b_{1}$ 
    along $\partial F_{2h}=\beta_{1}$,  we obtain 
    a representative of $[\EmbSurf]$ in this neighborhood.

    We claim that the Morse function $f$ can be extended 
    to all of $Y$, so that the 
    extension has one index three critical point and no additional 
    index zero critical points. To see this, extend $f$ to a Morse 
    function ${\widetilde f}$, and first cancel off all new index zero 
    critical points. This is a familiar argument from Morse 
    theory (see for instance~\cite{Milnor}): given another index zero 
    critical point $p'$, there is some index one critical point 
    $a$ which admits a unique flow to $p'$ (if there no 
    such index one critical points, then $p'$ would generate a 
    $\Z$ in the Morse complex for $Y$,
    which persists in $H_{0}(Y)$; but also, the sum of the other 
    index zero critical points would not lie in the image of 
    $\partial$, so it, too, would persist in homology, violating the 
    connectedness hypothesis of $Y$). Thus, we can cancel $p'$ and 
    the critical point $a$. 
    
    Next, we argue that the extension ${\widetilde f}$ need contain
    only one index three critical point, as well.  If there were two,
    call them $q$ and $q'$, we show that one of them can necessarily
    be canceled with an index two critical point other than
    $b_{1}$. If this could not be done, then both $q$ and $q'$ would
    have a unique flow-line to $b_{1}$. Thus, both $q$ and $q'$ would
    represent non-zero elements in $H_{3}(Y,\EmbSurf)\cong
    H^{0}(Y-\EmbSurf)$. But this is impossible since the complement
    $Y-\EmbSurf$ is connected, thanks to our homological assumption on
    $\EmbSurf$ (which ensures that $\EmbSurf$ admits a dual circle
    which hits it algebraically a non-zero number of times). In fact,
    the extension generically contains no flows between index $i$ and
    index $j$ critical points with $j\geq i$, hence giving us a
    Heegaard decomposition of $Y$.

    Thus, $Y$ has a handlebody decomposition
    $Y=U_{0}\cup_{\Sigma_{g}} U_{1}$, where $U_{0}$ is obtained from
    $V_{2h}$ by attaching a sequence of one-handles. The attaching
    regions for each of these one-handles consists of two disjoint
    disks in $\Sigma_{2h}$, which are disjoint from $\beta_1$. At
    least one of them has one component inside $F_{2h}$ and one
    outside.  This follows from the fact that $\beta_1$ is
    homologically trivial in $\Sigma_{2h}$, but homologically
    non-trivial in the final Heegaard surface $\Sigma$.  Let $\alpha_{2h+1}$ be the attaching
    circle for this one-handle.  After handleslides across
    $\alpha_{2h+1}$, we can arrange that all the other additional
    one-handles were attached in the complement of $F_{2h}$. The
    domain in $F_{2h}$ between and $\alpha_{2h+1}$ and $\beta_1$
    represents $\EmbSurf$.

\qed
\vskip.3cm

\subsection{The first Chern class formula}
\label{subsec:MaslovFormula}

Next, we give a proof Proposition~\ref{prop:EasyMeasureCalc}. Indeed,
we prove a more general result. But first, we introduce some data
associated to periodic domains.

A periodic domain $\PerDom$ is represented by an oriented two-manifold with 
boundary 
$\Phi \colon F\longrightarrow \Sigma$,
whose boundary maps under $\Phi$ into $\alphas \cup \betas$. 
We consider the pull-back bundle $\Phi^{*}(T\Sigma)$ over $F$. This bundle 
is canonically trivialized over the boundary: the velocity vectors of 
the attaching circles give rise to natural trivializations. We define the 
{\em Euler measure} of the periodic domain $\PerDom$ by the formula:
$$\chi(\PerDom)=\langle c_{1}(\Phi^{*}T\Sigma;\partial), F \rangle,$$
where $c_{1}(\Phi^{*}T\Sigma;\partial)$ is first Chern class of 
$\Phi^{*}T\Sigma$ relative to this boundary trivialization.
(It is easy to verify that $\chi(\PerDom)$ is independent of the 
representative $\Phi\colon F \longrightarrow \Sigma$.)

For example, if $\PerDom\subset \Sigma$ is a periodic domain all of whose coefficients 
are one or zero, with
$\partial \PerDom = \cup_{i=1}^{m}\gamma_{i}$ where the $\gamma_{i}$ are 
chosen among the $\alphas$ and the $\betas$, then
$\chi(\PerDom)$ agrees with the usual Euler characteristic of $\PerDom$, 
thought of as a subset of $\Sigma$.

Given a reference point $x\in\Sigma$, there is another quantity
associated to periodic domains, obtained from a natural generalization
of the local multiplicity $n_{x}(\PerDom)$ defined in
Section~\ref{HolDiskOne:sec:TopPrelim} of~\cite{HolDisk}. This
quantity, which we denote ${\overline n}_{x}(\PerDom)$, is defined by:
$$ {\overline n}_{x}(\sum_{i} a_{i} \cald_{i}) = \sum_{i} a_i
\left(\begin{array}{ll} 1 & {\text{if $x$ lies in the interior of
$\cald_i$}} \\
\OneHalf & {\text{if $x$ lies in the interior of some edge of
$\cald_i$}} \\
& {\text{or two vertices of $\cald_i$ are identified with
$x$}}  \\
\OneQuarter & {\text{if one vertex of $\cald_i$ is
identified with $x$}} \\
0 & {\text{if $x\not\in\cald_i$}}
\end{array}\right).
$$
Of course, if $x$ lies in $\CurveComp$, then ${\overline 
n}_{x}(\PerDom)=n_{x}(\PerDom)$. If $\PerDom$ has all multiplicities one or 
zero, and $x$ is contained in its boundary, then ${\overline 
n}_{x}(\PerDom)=\OneHalf$. 

\begin{prop}
    \label{prop:MeasureCalc}
    Fix a class $\xi\in H_{2}(Y;\Z)$, a base point $z\in\CurveComp$, 
    and a point $\x\in\Ta\cap \Tb$. Let $\PerDom$ be the periodic 
    domain associated to $z$ and $\xi$, and let $\spinc$ be the 
    $\SpinC$ structure $s_{z}(\x)$.
    Then the 
    evaluation of the first Chern class of $\spinc$ on $\xi$ 
    is calculated by
    $$\langle c_{1}(\spinc), \xi\rangle 
    = \chi(\PerDom) + 
    2\sum_{x_{i}\in \x} {\overline n}_{x_{i}}(\PerDom). $$
\end{prop}

Of course, Proposition~\ref{prop:EasyMeasureCalc} is a special case of
this result, since in that case, two of the $x_i$ are in the boundary
of $\PerDom$, so they have ${\overline n}_{x_i}=\OneHalf$. 

To prove the proposition, we need an explicit understanding of the
vector field belonging to the $\SpinC$ structure $s_{z}(\x)$.
Specifically, consider the normalized gradient vector field
$\frac{\grad f}{|\grad f|}$, restricted to the mid-level $\Sigma$ of
the Morse function $f$ (compatible with the given Heegaard
decomposition of $Y$). Clearly, the orthogonal complement of the
vector field is canonically identified with the tangent bundle of
$\Sigma$. Suppose, then, that $\gamma$ is a connecting trajectory
between an index one and an index two critical point (which passes
through $\Sigma$). We can replace the gradient vector field by another
vector field $v$ which agrees with $\frac{\grad f}{|\grad f|}$ outside
of a small three-ball neighborhood $B$, which meets $\Sigma$ in a disk
$\CDisk$.  Let $\tau$ be a trivialization of the two-plane field
$v^{\perp}|\partial \CDisk$ which extends as a trivialization of
$T\Sigma|\CDisk$. There is a well-defined relative first Chern class
$c_{1}(v,\tau)\in H^{2}(\CDisk,\partial\CDisk)$, which we can
calculate as follows:

\begin{lemma}
    \label{lemma:RelativeCOne}
    For $\CDisk$, $v$, and $\tau$ as above,
    the relative first Chern number is given by
    $$\langle c_{1}(v,\tau), [\CDisk,\partial \CDisk]\rangle=2$$
    (where we orient $\CDisk$ in the same manner as $\Sigma=\partial 
    U_{0}$).
\end{lemma}
    
\begin{proof}
    Using an 
    appropriate trivialization of the tangent bundle $TY|B$, we can view the 
    normalized gradient vector field $\frac{\grad f}{|\grad f|}$ as 
    constant over $\CDisk$. Let $S=\partial B$ be the boundary, which is 
    divided into two hemispheres $S=D_{1}\cup D_{2}$, so that the sphere $D_{1}\cup 
    \CDisk$ contains the index one critical point and $\CDisk\cup D_{2}$ 
    contains the index two critical point. We can replace $\frac{\grad 
    f}{|\grad f|}$ by another vector field $v$ which agrees with the 
    normalized gradient over $S$, and vanishes nowhere in $B$ (and hence 
    can be viewed as a unit vector field). With respect to the 
    trivialization of $TY|B$, we can think of the vector field as a map to 
    the two-sphere; indeed the 
    restriction
    $v \colon \CDisk \longrightarrow S^{2}$,
    is constant along the boundary circle, so it has a well-defined degree, which in the 
    present case is one, since
	$$
	-1 = \deg_{D_{1}}\left(\frac{\grad f}{|\grad f|}\right) + 
	\deg_{\CDisk}\left(\frac{\grad f}{|\grad f|}\right) 
	= \deg_{D_{1}}(v) 
	$$
    and
    \begin{eqnarray*}
	0 &=& \deg_{D_{1}}(v)+\deg_{\CDisk}(v).
    \end{eqnarray*}
    The line bundle we are considering, $v^{\perp}$, then, is the 
    pull-back of the tangent bundle to $S^{2}$, whose first Chern 
    number is the Euler characteristic for the sphere.
\end{proof}

\vskip.2cm
\noindent{\bf{Proof of Proposition~\ref{prop:MeasureCalc}}.}
    We find it convenient to consider domains with only non-negative 
    multiplicities; thus, we prove the following formula (for 
    sufficiently large $m$):
    \begin{equation}
	\label{eq:GenMasForm}
	\langle c_{1}(\spinc), \xi\rangle 
	= \chi(\PerDom+m[\Sigma]) + 
	2 \left(\sum_{x_{i}\in \x} {\overline 
	n}_{x_{i}}(\PerDom+m[\Sigma])\right)
	- 2 n_{z}(\PerDom+m[\Sigma]). 
    \end{equation}
    In fact, since 
    \begin{eqnarray*}
	\chi(\PerDom+m[\Sigma])&=&\chi(\PerDom)+m(2-2g), \\
	\sum_{x_{i}\in \x} {\overline n}_{x_{i}}(\PerDom+m[\Sigma])
	&=&
	mg+\sum_{x_{i}\in \x} {\overline n}_{x_{i}}(\PerDom) \\
	n_{z}(\PerDom+m[\Sigma])&=&m,
    \end{eqnarray*}
    Equation~\eqref{eq:GenMasForm} for any specific value of
    $m$ implies the formula 
    stated in the proposition. 
    
    The reformulation has the advantage that for $m$ sufficiently 
    large,
    $\PerDom+m[\Sigma]$ is represented by a map
    $\Phi \colon F \longrightarrow \Sigma $
    which is nowhere orientation-reversing, and whose restriction to 
    each boundary component is a diffeomorphism onto its 
    image (see Lemma~\ref{HolDiskOne:lemma:RepPerDom} of~\cite{HolDisk}).
    
\begin{figure}
	\mbox{\vbox{\epsfbox{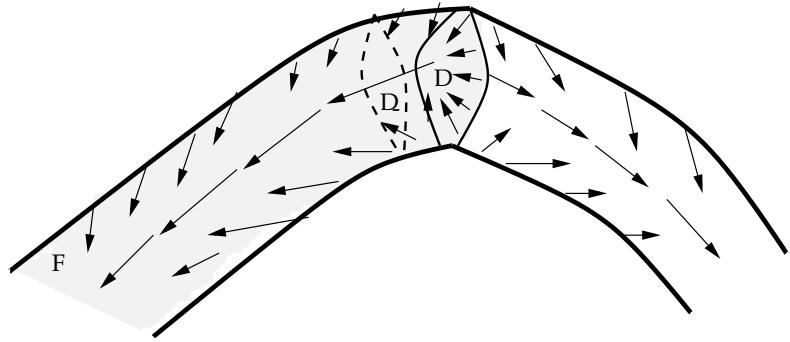}}}
	\caption{\label{fig:GradFlow}
	The gradient flow inside a one-handle. 
	The shaded region on the boundary of the one-handle is a piece of $F$; 
	the disk $D$ (with solid boundary, in the center) 
	goes through the index one critical point. Its translate $D_{-}$ 
	(with dotted boundary) does not, and the subregion of $F$ 
	terminating in the dotted circle, when capped off by $D_{-}$, is 
	transverse to the gradient flow.}
    \end{figure}

    Near each boundary component of $F$, we can identify a 
    neighborhood in $F$ with the half-open cylinder
    $[0,1)\times S^{1}$. Suppose that the image of the boundary 
    component is an 
    $\beta$ curve. The $\beta$ curve canonically bounds a disk in 
    $U_{1}$: this disk $D$ consists of points which
    flow (under $\grad f$) into the associated index two 
    critical point. Of course, we can glue this disk to $F$ along 
    the boundary, and correspondingly 
    extend $\Phi$ across 
    the disk as a map into $Y$, but then the gradient $\grad f$  vanishes at 
    some point of the extended map. To avoid this, we can back off 
    from the boundary of $F$: we delete a small neighborhood 
    $[0,\epsilon)\times S^{1}$ from $F$, to obtain a new 
    manifold-with-boundary $F^{-}$. In these local coordinates, now, 
    the boundary of $F^{-}$ is a translate of the $\beta$ curve 
    $\{\epsilon\}\times S^{1}$. Now, we can attach a translate of the 
    disk, $D_{-}$. Now, it is easy to see that (a smoothing of) the cap 
    $\left([\epsilon,1)\times S^{1}\right)\cup D_{-}$ is transverse to the gradient 
    flow $\grad f$.  (See the illustration 
    in Figure~\ref{fig:GradFlow}.)

    We can perform the analogous construction at the 
    $\alpha$-components of the boundary of $F$, only now, the $\alpha$ 
    curve bounds a disk $D$ in $U_{0}$, which  consists of points 
    flowing out of the corresponding index two critical point. By 
    cutting out a neighborhood of the boundary, and attaching a 
    translate of the $D$, we once again obtain a cap which is 
    transverse to the gradient flow $\grad f$.
    
    Observe that if $x_{i}\in\intPerDom$, then (if we chose the 
    above $\epsilon$ sufficiently small),
    \begin{equation}
	\label{eq:MultFormOne}
	{\overline n}_{x_{i}}(\PerDom)= \#\{x\in F^{-}\big| 
	\Phi(x)=x_{i}\}
    \end{equation}
    (with the same formula holding for $z$ in place of $x_{i}$).
    Moreover, if $x_{i}\in\partial\PerDom$, then
    \begin{equation}
	\label{eq:MultFormTwo}
	{\overline n}_{x_{i}}(\PerDom)
	= \OneHalf \# \{x\in \partial F| \Phi(x)=x_{i}\}
	+ \# \{x\in F^{-}| \Phi(x)=x_{i}\}.
    \end{equation}
    
    By adding the caps as above to $F^-$, we construct a closed,
    oriented two-manifold ${\widehat F}$ and a map
    $${\widehat \Phi}\colon{\widehat F}\longrightarrow Y,$$ which
    crosses the connecting trajectories between the index one and two
    critical points at each point $x\in F^{-}$ which maps under $\Phi$
    to $x_{i}$, and similarly, ${\widehat \Phi}$ crosses the
    connecting trajectory belonging to $z$ at those $x\in F^{-}$ which
    map under $\Phi$ to $z$.
    
    Away from these points, we have a canonical identification 
    $${\widehat \Phi}^{*}((\grad f)^{\perp})\cong 
    \Phi^{*}(v^{\perp}).$$
    By the local calculation from Lemma~\ref{lemma:RelativeCOne}, it 
    follows that
    \begin{equation}
	\label{eq:CompareGradV}
	\langle e\left({\widehat \Phi}^{*}(v^{\perp})\right), 
        {\widehat F} \rangle 
	= \langle e\left({\widehat \Phi}^{*}(\grad f^{\perp})\right), 
        {\widehat F} 
	\rangle + 2 \#\{x\in F^{-}| \Phi(x)=x_{i} \}
	- 2 \#\{x\in F^{-}| \Phi(x)=z\}.
    \end{equation}
    (Note that the term involving $z$ follows just as in the proof of 
    Lemma~\ref{lemma:RelativeCOne}, with the difference that now the 
    index of the vector field $v$
    around the corresponding critical point in $U_{0}$ is $+1$ rather than $-1$, 
    since the critical point has index zero rather than one.)
    
    Moreover, the Euler number of ${\widehat \Phi}^{*}(\grad 
    f^{\perp})$ is $\chi(\PerDom)$ plus the number of disks 
    which are attached to $F^{-}$ to obtain the closed manifold 
    ${\widehat F}$ (since each boundary disk is transverse 
    to the gradient flow, so 
    ${\grad f}^{\perp}$ is naturally 
    identified with the tangent bundle of the 
    disk, which has relative Euler number one relative to the 
    trivialization it gets from the bounding circle).
    But the number of such disks is simply $\#\{x\in\partial F| 
    \Phi(x)=x_{i}\}$. Combining this with 
    Equations~\eqref{eq:MultFormOne}, \eqref{eq:MultFormTwo}, and 
    \eqref{eq:CompareGradV}, we obtain Equation~\eqref{eq:GenMasForm}, 
    and hence proposition follows.
\qed

\section{Twisted Coefficients}
\label{sec:TwistedCoeffs}

We define here variants of the Floer homology groups constructed
in~\cite{HolDisk}: these are Floer homology groups with a ``twisted
coefficient system.'' The input here is is a three-manifold $Y$
equipped with a $\SpinC$ structure $\spinc$, and a module $M$ over
the group-ring $\Z[H^1(Y;\Z)]$. We begin with the definition in 
Subsection~\ref{subsec:DefTwisted}, discussing how the holomorphic triangle
construction needs to be modified in Subsection~\ref{subsec:Triangles}

\subsection{Twisted coefficients}
\label{subsec:DefTwisted}

We give first the ``universal construction'', using the free module
$M=\Z[H^1(Y;\Z)]$.  We need a surjective, additive assignment (in the
sense of Definition~\ref{HolDiskOne:def:Additive} of~\cite{HolDisk}):
$$A\colon \pi_2(\x,\y)\longrightarrow H^1(Y;\Z),$$ which is invariant
under the action of $\pi_2(\Sym^g(\Sigma))$.

We can construct such a map as follows.
A complete set of paths for $\spinc$ in the sense of 
Definition~\ref{HolDiskOne:def:CompSetPath} of~\cite{HolDisk}
gives rise to identifications for any $i,j$:
$$\pi_{2}(\x_{i},\x_{j})\cong \pi_{2}(\x_{0},\x_{0}),$$ by
$$\phi_{i}\csum\pi_{2}(\x_{i},\x_{j})\cong \pi_{2}(\x_{0},\x_{0})\csum
\phi_{j}.$$
These isomorphisms fit together in an additive manner, thanks to the 
associativity of $\csum$. We then use the splitting 
$\pi_{2}(\x_{0},\x_{0})\cong \Z\times H^{1}(Y;\Z)$ given by the
basepoint, followed by the natural projection to the second factor.

We can then define 
$$\uDel^\infty[\x,i] = \sum_{\y\in\Ta\cap\Tb}
\left(\sum_{\phi\in\pi_2(\x,\y)}
\#\ModFlow(\phi)e^{A(\phi)}{[\y,i-n_z(\phi)]}\right),
$$
which is a finite sum under the strong admissibility hypotheses.

Analogous constructions work for $\CFp$, $\CFm$, and $\CFa$, as well
(with, once again, weak admissibility sufficing for $\CFp$ and
$\CFa$). 

\begin{remark}
Note that there is a ``universal'' coefficient system for Lagrangian
Floer homology, with coefficients in a group-ring over $\pi_1(\Omega(L_0,L_1))$. In
fact, the construction we have here is a specialization of this: in
our case, the fundamental group of the configuration space is
$\Z\oplus H^1(Y,\Z)$, but the $\Z$ summand is already implicit in our
consideration of pairs $[\x,i]\in(\Ta\cap\Tb)\times \Z$.
\end{remark}

It is worth noting that, although the definition of the boundary map
still depends on a coherent system of orientations $\orient$, the
isomorphism class of the chain complex as a $\Z$-module 
does not: given a homomorphism
$\mu \colon H^1(Y;\Z) \longrightarrow \Zmod{2}$, the map
\begin{equation}
\label{eq:ChangeOrientations}
f(e^{h}[\x,i])=(-1)^{\mu(h)} e^{h} [\x,i]
\end{equation}
gives an isomorphism
from the chain complex using $\orient$ to the chain complex using
$\orient'$ with $\delta(\orient,\orient')=\mu$.

Note that as $\Z$-modules, all of these chain complexes have a natural
relative $\Z$ grading, which lifts the obvious relative
$\Zmod{{\mathfrak d}(\spinc)}$-grading. Specifically, given $g\otimes
[\x,i]$ and $h\otimes [\y,j]$ with $g,h\in H^1(Y;\Z)$, if we let
$\phi$ be the class with $A(\phi)=g-h$ and $n_z(\phi)=i-j$ (this now
uniquely specifies $\phi$), we let the relative grading between
$g\otimes [\x,i]$ and $h\otimes [\y,j]$ be given by the Maslov index
of $\phi$. In view of this, we can think of the corresponding
homologies as analogues of a construction of Fintushel and Stern, for
$\Z$ graded instanton homology (see~\cite{FintSternZGrading}).

For any $\Z[H^1(Y;\Z)]$-module $M$, we have homology groups defined by
$$\uHF(Y,\spinc;M)=H_*\left(\uCF(Y,\spinc)\otimes_{Z[H^1(Y;\Z)]}M\right)$$
(where $\uHF$ can be any of $\uHFinf$, $\uHFp$, $\uHFm$, or $\uHFa$).
The homology groups from~\cite{HolDisk} (with ``untwisted
coefficients'') are special cases of this construction, using the
module $M=\Z$, thought of as the trivial $\Z[H^1(Y;\Z)]$-module. (In
fact, writing $b=b_1(Y)$, 
the $2^b$ different choices of orientation systems over $\Z$
correspond to the $2^b$ different module structures on $\Z$,
induced from the $2^b$ ring homomorphisms $\Z[H^1(Y;\Z)]
\longrightarrow \Z$.)

Note also that the action of $H_1(Y;\Z)/\Tors$ on $\CFinf(Y,\spinc)$
has an interpretation in this world: the action
of $\zeta\in H_1(Y;\Z)$ on $[\x,i]\in\CFinf(Y,\spinc)$ as defined in
Subsection~\ref{HolDiskOne:subsec:DefAct} of~\cite{HolDisk} can be
represented by $\langle \uDel [\x,i], \zeta \rangle$,
where the angle brackets represent the
natural pairing
$\Z[H^1(Y;\Z)]\otimes (H_1(Y;\Z)/\Tors) \longrightarrow \Z$, 

A modification of the techniques from~\cite{HolDisk} gives the following:

\begin{theorem}
Let $Y$ be a three-manifold equipped with a $\SpinC$ structure
$\spinc$ and a $\Z[H^1(Y;\Z)]$-module $M$. Let
$(\Sigma,\alphas,\betas,z)$ be a strongly $\spinc$-admissible Heegaard
diagram for $Y$. Then the groups $\uHFinf(\alphas,\betas,\spinc,M)$,
$\uHFp(\alphas,\betas,\spinc,M)$, $\uHFm(\alphas,\betas,\spinc,M)$,
and $\uHFa(\alphas,\betas,\spinc,M)$ are invariant under changes of
almost complex structures and isotopies. These groups are all modules
over the group-ring $\Z[H^1(Y;\Z)]$.
\end{theorem}

Independence of complex structure follows exactly as
in~\cite{HolDisk}.  For isotopy invariance, observe that an isotopy
$\Psi_t$ as in Subsection~\ref{HolDiskOne:sec:Isotopies}
of~\cite{HolDisk} allows one to transfer an additive map $A$ from
$\pi_2(\x,\y)$ for $\x,\y\in\Ta\cap\Tb$ to an additive map on
$\pi_2(\x',\y')$ for $\x',\y'\in\Psi_1(\Ta)\cap\Tb$. Stabilization
follows as in~\cite{HolDisk}, while to understand handleslide
invariance, describe how to modify the holomorphic triangle
construction to take into account the twisted coefficient system.

\subsection{Triangles and twisted coefficients}
\label{subsec:Triangles}

To understand the triangle construction with twisted coefficients, we
set up some topological preliminaries concerning relative $\SpinC$ structures

\subsubsection{Relative $\SpinC$ structures}

Continuing notation from Subsection~\ref{HolDiskOne:sec:HolTriangles}
of~\cite{HolDisk}, let $(\Sigma,\alphas,\betas,\gammas,z)$ be a
pointed Heegaard triple, and let $X_{\alpha,\beta,\gamma}$ be the
induced cobordism between $Y_{\alpha,\beta}$, $Y_{\beta,\gamma}$, and
$Y_{\alpha,\gamma}$. Fix $\SpinC$ structures $\spinct_{\alpha,\beta}$,
$\spinct_{\beta,\gamma}$, $\spinct_{\alpha,\gamma}$ over the three
boundary components, with
$\epsilon(\spinct_{\alpha,\beta},\spinct_{\beta,\gamma},\spinct_{\alpha,\gamma})=
0$.  Fix complete sets of paths for each of these three $\SpinC$
structures (in the sense of
Definition~\ref{HolDiskOne:def:CompSetPath} of~\cite{HolDisk}). This
gives us identifications $$\pi_2(\x_0,\y_0,\w_0)=\pi_2(\x,\y,\w),$$
where $\x_0$ and $\x$ (resp. $\y_0$ and $\y$, resp. $\w_0$ and $\w$)
both represent $\spinct_{\alpha,\beta}$
(resp. $\spinct_{\beta,\gamma}$ resp. $\spinct_{\alpha,\gamma}$).

In effect, this allows us to think of $\pi_2(\x_0,\y_0,\w_0)$ as an
affine space for $H^2(X,Y;\Z)$
(c.f. Proposition~\ref{HolDiskOne:prop:CalcPiTwo} of~\cite{HolDisk}),
which maps onto the space of $\SpinC$ structures extending
$\spinct_{\alpha,\beta}$, $\spinct_{\beta,\gamma}$,
$\spinct_{\alpha,\gamma}$
(c.f. Proposition~\ref{HolDiskOne:prop:SpinCForTriangles} of~\cite{HolDisk}). 
When thinking of
$\pi_2(\x_0,\y_0,\w_0)$ in this manner, we refer to it as a space of
relative $\SpinC$ structures, and denote it by
$\RelSpinC(X_{\alpha,\beta,\gamma})$.

The subset of $\RelSpinC(X_{\alpha,\beta,\gamma})$
representing a fixed (absolute) $\SpinC$ structure structure
$\spinc_{\alpha,\beta,\gamma}$
will be denoted
$\SpinC(X_{\alpha,\beta,\gamma};\spinc_{\alpha,\beta,\gamma})$.

We will use this terminology for higher polygons, as well.

\subsubsection{The maps with twisted coefficients}

The
space of relative $\SpinC$ structures 
$\RelSpinC(X_{\alpha,\beta,\gamma};\spinc_{\alpha,\beta,\gamma})$ (which induce a given $\SpinC$ structure
$\spinc_{\alpha,\beta,\gamma}$ over $X_{\alpha,\beta,\gamma}$)
is a space with a natural action of $H^1(Y_{\alpha,\beta};\Z)\times
H^1(Y_{\beta,\gamma};\Z)\times H^1(Y_{\alpha,\gamma};\Z)$. As such, it
can be used to induce an $H^1(Y_{\alpha,\gamma};\Z)$-module from a
pair $M_{\alpha,\beta}$ and $M_{\beta,\gamma}$ of
$H^1(Y_{\alpha,\beta};\Z)$ and $H^1(Y_{\beta,\gamma};\Z)$-modules: 
$$\left\{M_{\alpha,\beta}\otimes M_{\beta,\gamma}\right\}^{\spinc_{\alpha,\beta,\gamma}}
= \frac{(m_{\alpha,\beta},m_{\beta,\gamma},\relspinc)\in M_{\alpha,\beta}\times
M_{\beta,\gamma}
\times \RelSpinC(X_{\alpha,\beta,\gamma},\spinc_{\alpha,\beta,\gamma})}
{(m_{\alpha,\beta},m_{\beta,\gamma},\relspinc)\sim
(h_{\alpha,\beta}\cdot m_{\alpha,\beta},
 h_{\beta,\gamma}\cdot m_{\beta,\gamma},
 (h_{\alpha,\beta}\times h_{\beta,\gamma}\times 0)\cdot \relspinc)
},
$$
where $h_{\alpha,\beta}$ and $h_{\beta,\gamma}$ are arbitrary 
elements of $H^1(Y_{\alpha,\beta};\Z)$ and $H^1(Y_{\beta,\gamma};\Z)$ respectively.

Fix a $\SpinC$ structure $\spinc$ over $X_{\alpha,\beta,\gamma}$, whose restriction to $Y_{\alpha,\beta}$ and $Y_{\beta,\gamma}$ is $\spinct_{\alpha,\beta}$ and $\spinct_{\beta,\gamma}$ respectively. 
We can now define a map
\begin{eqnarray*}
\lefteqn{\ufinf(~\cdot~,\spinc)\colon 
\uCFinf(Y_{\alpha,\beta},\spinct_{\alpha,\beta};M_{\alpha,\beta})
\otimes
\uCFinf(Y_{\beta,\gamma},\spinct_{\beta,\gamma};M_{\beta,\gamma})} \\
&&\hskip2in
\longrightarrow
\uCFinf(Y_{\alpha,\gamma},\spinct_{\alpha,\gamma};\left\{
M_{\alpha,\beta}\otimes M_{\beta,\gamma}\right\}^{\spinc}),
\end{eqnarray*}
by the formula:
\begin{eqnarray}
\lefteqn{\ufinf(m_{\alpha,\beta}[\x,i]\otimes m_{\beta,\gamma}[\y,j];\spinc)
=}  \nonumber 
\\
&& \sum_{\w\in\Ta\cap\Tb} 
\sum_{\{\psi\in\pi_2(\x,\y,\w)\big| \spinc_z(\psi)=\spinc\}}
\left(\#\ModFlow(\psi)\right)
\{m_{\alpha,\beta}\otimes
m_{\beta,\gamma}\otimes \sRelSpinC_z(\psi)\}\cm [\w,i+j-n_z(\psi)].
\label{eq:DefTriangleTwisted}
\end{eqnarray}
The braces above indicate the natural map
$$\left\{\cdot\otimes \cdot\otimes \cdot\right\} \colon
M_{\alpha,\beta}\otimes M_{\beta,\gamma}\otimes
\RelSpinC(X_{\alpha,\beta,\gamma},\spinc)
\longrightarrow \{M_{\alpha,\beta}\otimes M_{\beta,\gamma}\}^{\spinc}.$$

The following analogue of Theorem~\ref{HolDiskOne:thm:HolTriangles}
of~\cite{HolDisk} holds in the present context:

\begin{theorem}
\label{thm:HolTrianglesTwisted}
Let $(\Sigma,\alphas,\betas,\gammas,z)$ be a pointed Heegaard
triple-diagram, which is strongly $\spinc$-admissible for some $\SpinC$
structure $\spinc$ over the underlying  four-manifold $X$, and fix modules
$M_{\alpha,\beta}$ and $M_{\beta,\gamma}$ for $H^1(Y_{\alpha,\beta};\Z)$ and
$H^1(Y_{\beta,\gamma};\Z)$ respectively.
Then the sum on
the right-hand-side of Equation~\eqref{eq:DefTriangleTwisted} is finite,
giving rise to a chain map which also induces maps on homology:
\begin{eqnarray*}
\lefteqn{\uFinf{\alpha,\beta,\gamma}(\cdot,\spinc_{\alpha,\beta,\gamma})
\colon \uHFinf(Y_{\alpha,\beta},\spinct_{\alpha,\beta};M_{\alpha,\beta})
\otimes \uHFinf(Y_{\beta,\gamma},\spinct_{\beta,\gamma};M_{\beta,\gamma})} \\
&&\hskip2in
\longrightarrow 
\uHFinf(Y_{\alpha,\gamma},\spinct_{\alpha,\gamma};\{M_{\alpha,\beta}\otimes
M_{\beta,\gamma}\}^{\spinc_{\alpha,\beta,\gamma}}) \\
\lefteqn{\uFleq{\alpha,\beta,\gamma}(\cdot,\spinc_{\alpha,\beta,\gamma})
\colon \uHFleq(Y_{\alpha,\beta},\spinct_{\alpha,\beta};M_{\alpha,\beta})
\otimes \uHFleq(Y_{\beta,\gamma},\spinct_{\beta,\gamma};M_{\beta,\gamma})} \\
&&\hskip2in  \longrightarrow
\uHFleq(Y_{\alpha,\gamma},\spinct_{\alpha,\gamma};\{M_{\alpha,\beta}\otimes M_{\beta,\gamma}\}^{\spinc_{\alpha,\beta,\gamma}}).
\end{eqnarray*}
The induced chain map
\begin{eqnarray*}
\lefteqn{
\ufp{\alpha,\beta,\gamma}(\cdot,\spinc_{\alpha,\beta,\gamma})\colon
\uCFp(Y_{\alpha,\beta},\spinct_{\alpha,\beta};M_{\alpha,\beta})
\otimes \uCFleq(Y_{\beta,\gamma},\spinct_{\beta,\gamma};M_{\beta,\gamma})} \\
&&\hskip2in
\longrightarrow \uCFp(Y_{\alpha,\gamma},\spinct_{\alpha,\gamma};
\{M_{\alpha,\beta}\otimes M_{\beta,\gamma}\}^{\spinc_{\alpha,\beta,\gamma}}) 
\end{eqnarray*}
gives a well-defined  chain map when the triple diagram is only weakly
admissible, and the Heegaard diagram $(\Sigma,\betas,\gammas,z)$ is
strongly admissible for $\spinct_{\beta,\gamma}$. 
In fact, the induced map
\begin{eqnarray*}
\lefteqn{
\ufa(\cdot,\spinc_{\alpha,\beta,\gamma})\colon
\uCFa(Y_{\alpha,\beta},\spinct_{\alpha,\beta};M_{\alpha,\beta})
\otimes \uCFa(Y_{\beta,\gamma},\spinct_{\beta,\gamma};M_{\beta,\gamma})
} \\
&& \hskip2in \longrightarrow
\uCFa(Y_{\alpha,\gamma},\spinct_{\alpha,\gamma};\{M_{\alpha,\beta}\otimes M_{\beta,\gamma}\}^{\spinc_{\alpha,\beta,\gamma}})
\end{eqnarray*}
gives a well-defined chain map when the diagram is weakly admissible.
There are induced maps on homology:
\begin{eqnarray*}
\lefteqn{\uFa{\alpha,\beta,\gamma}(\cdot,\spinc_{\alpha,\beta,\gamma})
\colon \uHFa(Y_{\alpha,\beta},\spinct_{\alpha,\beta};M_{\alpha,\beta})
\otimes \HFa(Y_{\beta,\gamma},\spinct_{\beta,\gamma};M_{\beta,\gamma}) } \\
&&\hskip2in \longrightarrow \HFa(Y_{\alpha,\gamma},\spinct_{\alpha,\gamma};\{M_{\alpha,\beta}\otimes M_{\beta,\gamma}\}^{\spinc_{\alpha,\beta,\gamma}}) \\
\lefteqn{F^+_{\alpha,\beta,\gamma}(\cdot,\spinc_{\alpha,\beta,\gamma})
\colon \uHFp(Y_{\alpha,\beta},\spinct_{\alpha,\beta})
\otimes \uHF^{\leq 0}(Y_{\beta,\gamma},\spinct_{\beta,\gamma})} \\
&&\hskip2in \longrightarrow \uHFp(Y_{\alpha,\gamma},\spinct_{\alpha,\gamma};\{M_{\alpha,\beta}\otimes M_{\beta,\gamma}\}^{\spinc_{\alpha,\beta,\gamma}})).
\end{eqnarray*}
\end{theorem}

Independence of complex structure and isotopy invariance of this map
proceeds exactly as in~\cite{HolDisk}
(c.f. Propositions~\ref{HolDiskOne:prop:TrianglesJIndep}
and~\ref{HolDiskOne:prop:TriangleIsotopyInvariance} of~\cite{HolDisk}
respectively). Associativity, on the other hand, can be given a
the following sharper statement.

Observe first that there is a canonical gluing
$$\RelSpinC(X_{\alpha,\beta,\gamma},\spinc_{\alpha,\beta,\gamma})\times
\RelSpinC(X_{\alpha,\gamma,\delta},\spinc_{\alpha,\gamma,\delta})
\longrightarrow \RelSpinC(X_{\alpha,\beta,\gamma,\delta}) $$
which maps onto the set of all relative $\SpinC$
structures over
$X_{\alpha,\beta,\gamma,\delta}$ whose restrictions to
$X_{\alpha,\beta,\gamma}$ and $X_{\alpha,\gamma,\delta}$ represent
$\SpinC$ structures $\spinc_{\alpha,\beta,\gamma}$ and
$\spinc_{\alpha,\gamma,\delta}$ respectively. Thus, the set of $\SpinC$ 
induced structures in $X_{\alpha,\beta,\gamma,\delta}$ under this map consists of a 
$\delta H^1(Y;\Z)$-orbit.
Using this
gluing, we obtain an identification 
\begin{eqnarray*}
\lefteqn{\{\{M_{\alpha,\beta}\otimes
M_{\beta,\gamma}\}^{\spinc_{\alpha,\beta,\gamma}}\otimes
M_{\gamma,\delta}\}^{\spinc_{\alpha,\beta,\delta}}} \\
&\cong & 
\coprod_{\{\spinc\in\SpinC(X_{\alpha,\beta,\gamma,\delta})\big|
\spinc|X_{\alpha,\beta,\gamma}=\spinc_{\alpha,\beta,\gamma},
\spinc|X_{\alpha,\gamma,\delta}=\spinc_{\alpha,\gamma,\delta}\}}
\{M_{\alpha,\beta}\otimes M_{\beta,\gamma}\otimes
M_{\gamma,\delta}\}^{\spinc},
\end{eqnarray*} where 
$\left\{M_{\alpha,\beta}\otimes M_{\beta,\gamma}
\otimes M_{\gamma,\delta}\right\}^{\spinc}$
denotes the $H^1(Y_{\alpha,\delta};\Z)$-module induced from $M_{\alpha,\beta}$, $M_{\beta,\gamma}$, 
$M_{\gamma,\delta}$ and 
the set of relative $\SpinC$ structures inducing the given $\SpinC$ structure $\spinc$
over the four-manifold $X_{\alpha,\beta,\gamma,\delta}$.

\begin{theorem}
\label{thm:AssociativityTriangles}
Let $(\Sigma,\alphas,\betas,\gammas,\deltas,z)$ be a
pointed Heegaard quadruple which is strongly
${\mathfrak S}$-admissible, where ${\mathfrak S}$ is
a $\delta H^1(Y_{\beta,\delta})+ \delta
H^1(Y_{\alpha,\gamma})$-orbit in
$\SpinC(X_{\alpha,\beta,\gamma,\delta})$.
Fix also modules
$M_{\alpha,\beta}$, $M_{\beta,\gamma}$, and $M_{\gamma,\delta}$
for $H^1(Y_{\alpha,\beta};\Z)$,
$H^1(Y_{\beta,\gamma};\Z)$, $H^1(Y_{\beta,\gamma};\Z)$, and $H^1(Y_{\gamma,\delta};\Z)$ 
respectively. 

Then,
\begin{eqnarray*}
\lefteqn{\sum_{\spinc\in{\mathfrak S}}
\uFstar{\alpha,\gamma,\delta}
(\uFstar{\alpha,\beta,\gamma}
(\xi_{\alpha,\beta}\otimes
\theta_{\beta,\gamma}; 
\spinc_{\alpha,\beta,\gamma})\otimes 
\theta_{\gamma,\delta};\spinc_{\alpha,\gamma,\delta})} \\
&=&\sum_{\spinc\in{\mathfrak S}}
\uFstar{\alpha,\beta,\delta}(\xi_{\alpha,\beta}\otimes 
\uFleq{\beta,\gamma,\delta}
(\theta_{\beta,\gamma}\otimes\theta_{\gamma,\delta};\spinc_{\beta,\gamma,\delta});
\spinc_{\alpha,\beta,\delta}),
\end{eqnarray*}
where $\uFstar{}=\uFinf{}$, $\uFp{}$ or $\uFm{}$; also,
\begin{eqnarray*}
\lefteqn{\sum_{\spinc\in{\mathfrak S}}
\uFa{\alpha,\gamma,\delta}(\uFa{\alpha,\beta,\gamma}
(\xi_{\alpha,\beta}\otimes
\theta_{\beta,\gamma}; 
\spinc_{\alpha,\beta,\gamma})\otimes 
\theta_{\gamma,\delta};\spinc_{\alpha,\gamma,\delta})} \\
&=&\sum_{\spinc\in{\mathfrak S}}
\uFa{\alpha,\beta,\delta}(\xi_{\alpha,\beta}\otimes 
\uFa{\beta,\gamma,\delta}
(\theta_{\beta,\gamma}\otimes\theta_{\gamma,\delta};\spinc_{\beta,\gamma,\delta});
\spinc_{\alpha,\beta,\delta}),
\end{eqnarray*}
where we are taking coefficients in 
coefficients in 
$\coprod_{\spinc\in{\mathfrak S}}
\{M_{\alpha,\beta}\otimes M_{\beta,\gamma}\otimes 
M_{\gamma,\delta}\}^{\spinc}$ over $Y_{\alpha,\delta}$.
\end{theorem}

\begin{proof}
The proof is the same as the proof of 
Theorem~\ref{HolDiskOne:thm:Associativity} of~\cite{HolDisk}, only keeping track now of the homotopy classes
of the corresponding triangles.
\end{proof}

\subsubsection{Handleslide invariance}
\label{subsec:HandleSlidesTwisted}

With the holomorphic triangles in place, the proof of handleslide
invariance proceeds as it did in~\cite{HolDisk}, with the following remarks.

Recall that the map given by a handleslide (as in
Theorem~\ref{HolDiskOne:thm:HandleslideInvariance} of~\cite{HolDisk})
is induced from a Heegaard triple $(\Sigma,\alphas,\betas,\gammas,z)$,
which represents the cobordism $X_{\alpha,\beta,\gamma}$ obtained from
$[0,1]\times Y$ by deleting a bouquet of circles. Here,
$Y_{\alpha,\beta}\cong Y$, $Y_{\beta,\gamma}\cong \#^g(S^1\times
S^2)$, and $Y_{\alpha,\gamma}\cong Y$. Now, our input includes an
arbitrary $\Z[H^1(Y;\Z)]$ module $M$. For the handleslide map, we
consider the trivial $H^1(Y_{\beta,\gamma};\Z)$-module
$M_{\beta,\gamma}\cong\Z$ (so that $\uHFleq(Y_{\beta,\gamma},M)\cong
\HFleq(\#^g(S^1\times S^2))$ is equipped with its top-dimensional
generator $\Theta_{\beta,\gamma}$).  
It is easy to see that for this choice of $M_{\beta,\gamma}$,
there is also a canonical identification of $\Z[H^1(Y;\Z)]$-modules
$$M\cong \{M\otimes M_{\beta,\gamma}\},$$ where the pairing here uses
the cobordism $X_{\alpha,\beta,\gamma}$. 

\section{Surgery exact sequences}
\label{sec:Surgeries}

We investigate how surgeries on a three-manifold affect its
invariants.  We consider first the effect on $\HFp$ of $+1$ surgeries
on integral homology three-spheres, then a generalization which holds
for arbitrary (closed, oriented) three-manifolds, and then the case of
fractional $1/q$-surgeries on an integral homology three-sphere.  This
latter case uses the homology theories with twisted coefficients.  We
then give analogous results for $\HFa$.  After this, we present a
surgery formula for integer surgeries.  In the final subsection, we
consider a $+1$ surgery formula with twisted coefficients.

\subsection{$+1$ surgeries on an integral homology three-sphere}

We start with the case of a homology three-sphere $Y$. Let $K\subset Y$
be a knot. Let $Y_0$ be the manifold obtained by $0$-surgery on $K$,
and $Y_{1}$ be obtained by $(+1)$-surgery.
Let 
$$\HFp(Y_0)\cong \bigoplus_{\spinc\in\SpinC(Y_0)}
\HFp(Y_0,\spinc),$$
viewed as a $\Zmod{2}$-relatively graded group.
In fact, we will view the homology groups $\HFp(Y)$ and $\HFp(Y_1)$ as
$\Zmod{2}$-graded, as well.

\begin{theorem}
\label{thm:ExactOne}
There is a $U$-equivariant exact sequence of relatively
$\Zmod{2}$-graded complexes:
$$\begin{CD} ... @>>> \HFp(Y)@>{F_1}>>
\HFp(Y_0)@>{F_2}>> \HFp(Y_1) @>{F_3}>> ...
\end{CD}
$$ 
In fact, if we give $\HFp(Y)$ and $\HFp(Y_1)$ absolute $\Zmod{2}$-gradings so that
$\chi(\HFa(Y))=\chi(\HFa(Y_1))=+1$, then $F_3$ preserves degree.
\end{theorem}

The maps in Theorem~\ref{thm:ExactOne} are constructed with the help of 
holomorphic triangles. Thus, we must construct
compatible Heegaard decompositions for all three manifolds $Y$, $Y_0$,
and $Y_1$. Exactness is then proved using a filtration on the homology
groups above, together with the homological-algebraic constructions
used in establishing the surgery sequences for instanton Floer
homology (see~\cite{FloerTriangles},
\cite{BraamDonaldson}). The proof occupies the rest of the present subsection.

\begin{lemma}
\label{lemma:HeegaardDiagrams}
There is a pointed Heegaard multi-diagram
$$(\Sigma,\alphas,\betas,\gammas,\deltas,z)$$
with the property that
\begin{list}
        {(\arabic{bean})}{\usecounter{bean}\setlength{\rightmargin}{\leftmargin}}
\item 
\label{item:WhatTheyDescribe}
the Heegaard diagrams $(\Sigma,\alphas,\betas)$,
$(\Sigma,\alphas,\gammas)$, and $(\Sigma,\alphas,\deltas)$ describe
$Y$, $Y_0$, and $Y_1$ respectively,
\item 
\label{item:BetasGammasDeltas}
for each $i=1,...,g-1$, the curves $\beta_i$, $\gamma_i$, and
$\delta_i$ are small isotopic translates of one another, each pairwise
intersecting in a pair of canceling transverse intersection points
(where the isotopies are supported in the complement of $z$),
\item 
\label{item:SurgeryCurves}
the curve $\gamma_g$ is isotopic to the juxtaposition
of $\delta_g$ and $\beta_g$ (with appropriate orientations),
\item 
\label{item:PosAndNeg}
every non-trivial multi-periodic domain has both positive and
negative coefficients.
\end{list}
\end{lemma}

\begin{proof}
Consider a Morse function on $Y-\Nbd{K}$ with one index zero critical
point, $g$ index one critical points and $g-1$ index two critical
points. Let $\Sigma$ be the $3/2$-level of this function, $\alphas$ be
the curves where $\Sigma$ meets the ascending manifolds of the index
one critical points in $\Sigma$, and let $\beta_1,...,\beta_{g-1}$ be
the curves where $\Sigma$ meets the
descending manifolds of the index two critical points. By gluing
in the solid torus in three possible ways, we get the manifolds $Y$,
$Y_0$, $Y_1$. Extending the given Morse function to the glued in solid
tori, (by introducing an additional index two and index three critical
point), we obtain Heegaard decompositions for the manifolds $Y$,
$Y_0$, and $Y_1$. We let $\gamma_i$ and $\delta_i$ be small
perturbations of $\beta_i$ for $i=1,...,g-1$. In this manner, we have
satisfied
Properties~(\ref{item:WhatTheyDescribe})-(\ref{item:SurgeryCurves}).

To satisfy Property~(\ref{item:PosAndNeg}), we wind to achieve weak
admissibility for all $\SpinC$ structures for the Heegaard subdiagram
$(\Sigma,\alphas,\gammas,z)$: in fact, we can use a volume form over
$\Sigma$ for which all such doubly-periodic domains have zero signed area
(c.f. Lemma~\ref{HolDiskOne:lemma:EnergyZero} of~\cite{HolDisk}). Then, for the
$\{\beta_1,...\beta_{g-1}\}$ and $\{\delta_1,...,\delta_{g-1}\}$, we
use small Hamiltonian translates of the
$\{\gamma_1,...,\gamma_{g-1}\}$ (ensuring that the corresponding new
periodic domains each have zero energy). There is a triply-periodic
domain which forms the homology between $\beta_g$, $\gamma_g$, and
$\delta_g$ in a torus summand of $\Sigma$ containing no other
$\beta_i$ or $\gamma_i$ (for $i\neq g$). By adjusting the areas of the
two triangles with non-zero area, we can arrange for the signed area
of the triply-periodic domain to vanish.

% export at 65% 
% maybe 50% better
\begin{figure}
\mbox{\vbox{\epsfbox{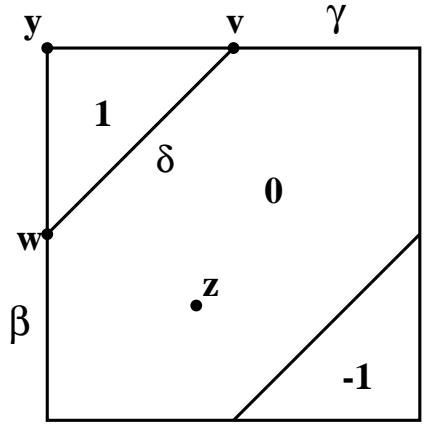}}}
\caption{\label{fig:onesurg}
This picture takes place in the torus, with the usual edge
identifications. The integers denote multiplicities for a
triply-periodic domain.}
\end{figure}

\end{proof}

\begin{figure}
\mbox{\vbox{\epsfbox{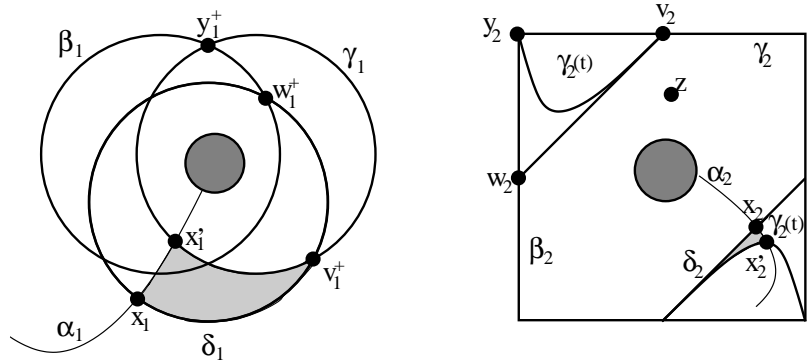}}}
\caption{\label{fig:LabelledSurgeries}
$+1$-surgery, $g=2$. The left side takes place in an annulus, the
right side in a torus minus a disk; both are pieces of our genus two
surface $\Sigma$ (the central disk missing from the annulus and the
disk removed from the torus are both indicated by large grey circles).
We have curves $\{\beta_1,\beta_2\}$, $\{\gamma_1,\gamma_2\}$ and
$\{\delta_1,\delta_2\}$ as in Lemma~\ref{lemma:HeegaardDiagrams}, with
intersection points $\Theta_{\beta,\gamma}=\{y_1^+,y_2\}$,
$\Theta_{\gamma,\delta}=\{v_1^+,v_2\}$, and
$\Theta_{\beta,\delta}=\{w_1^+,w_2\}$. The curve $\gamma_2(t)$ is
isotopic to $\gamma_2$, but it approximates the juxtaposition of
$\beta_2$ and $\delta_2$.  We have also pictured arcs in $\alpha_1$
and $\alpha_2$. There is an intersection point $\x=\{x_1,x_2\}$ for
$\Ta\cap\Td$, and its nearest point $\Ta\cap{\mathbb T}_{\gamma(t)}$,
$\{x_1',x_2'\}=\rho(\x)$. Observe the two lightly shaded triangles:
they correspond to the canonical triangle in
$\pi_2(\rho(\x),\Theta_{\gamma,\delta},\x)$.}
\end{figure}

For $i=1,...,g-1$, label $$y_i^\pm = \beta_i\cap\gamma_i,~~~~~~~~
v_i^\pm=\gamma_i\cap\delta_i, ~~~~~~~~ w_i^\pm=\beta_i\cap\delta_i,$$
where the sign indicates the sign of the intersection point.  Also,
let $$y_g = \beta_g\cap\gamma_g,~~~~~~~~ v_g=\gamma_g\cap\delta_g,
~~~~~~~~ w_g=\beta_g\cap\delta_g.$$ Then, let
$\Theta_{\beta,\gamma}=\{y_1^+,...,y^+_{g-1},y_g\}$,
$\Theta_{\gamma,\delta}=\{v_1^+,...,v^+_{g-1},v_g\}$,
$\Theta_{\beta,\delta}=\{w_1^+,...,w^+_{g-1},w_g\}$ denote the
corresponding intersection points between $\Tb\cap\Tc$, $\Tc\cap\Td$
and $\Tb\cap\Td$. (See Figure~\ref{fig:LabelledSurgeries} for an
illustration.)

\begin{prop}
\label{prop:EltsAreCycles}
The elements $\theta_{\beta,\gamma}=[\Theta_{\beta,\gamma},0]$, 
$\theta_{\gamma,\delta}=[\Theta_{\gamma,\delta},0]$, 
$\theta_{\beta,\delta}=[\Theta_{\beta,\delta},0]$ are cycles in
$\CFinf(\Tb,\Tc)$, $\CFinf(\Tc,\Td)$ and $\CFinf(\Tb,\Td)$ respectively.
\end{prop}

\begin{proof}
Note that the three-manifolds described here are $(g-1)$-fold connected sums of
$S^1\times S^2$, so the result follows from Proposition~\ref{prop:ConnSumS1S2}
(or, alternatively, see Section~\ref{HolDiskOne:sec:HandleSlides} 
of~\cite{HolDisk}).
\end{proof}

We can reduce the study of holomorphic triangles belonging to
$X_{\beta,\gamma,\delta}$ to holomorphic triangles in the first
symmetric product of the two-torus, with the help of the following
analogue of the gluing theory used to establish stabilization
invariance of the Floer homology groups. 

\begin{theorem}
\label{thm:GlueTriangles}
Fix a pair of Heegaard diagrams 
\begin{eqnarray*}
(\Sigma,\betas,\gammas,\deltas,z) &{\text{and}}&
(E,\beta_0,\gamma_0,\delta_0,z_0), 
\end{eqnarray*}
where $E$ is a Riemann surface of genus one. We will form the
connected sum $\Sigma\# E$, where the connected sum points are near
the distinguished points $z$ and $z_0$ respectively. Fix intersection
points $\x,\y,\w$ for the first diagram and a class
$\psi\in\pi_2(\x,\y,\w)$, and intersection points $x_0$, $y_0$, and
$w_0$ for the second, with a triangle $\psi_0\in\pi_2(x_0,y_0,w_0)$
with $\Mas(\psi)=\Mas(\psi_0)=0$.  Suppose moreover that
$n_{z_0}(\psi_0)=0$.  Then, for a suitable choice of complex structures
and perturbations, we
have a diffeomorphism of moduli spaces: $$\ModFlow(\psi')\cong
\ModFlow(\psi)\times
\ModFlow(\psi_0),$$
where $\psi'\in \pi_2(\x\times x_0, \y\times y_0,\w\times w_0)$ is 
the triangle for $\Sigma\# E$ whose domain on the $\Sigma$-side 
agrees with $\cald(\psi)$, and whose domain on the $E$-side agrees with
$\cald(\psi_0)+n_{z}(\psi)[E]$. 
\end{theorem}

\begin{proof}
The proof is obtained by suitably modifying 
Theorem~\ref{HolDiskOne:thm:Gluing} of~\cite{HolDisk}. 

Suppose that $u$ and $u_0$ are holomorphic representatives
of $\psi$ and $\psi_0$ respectively. We obtain a nodal pseudo-holomorphic
disk $u\vee u_0$ in the singular space 
$\Sym^{g+1}(\Sigma\vee E)$
specified as follows:
\begin{itemize}
\item At the stratum $\Sym^{g}(\Sigma)\times
\Sym^{1}(E)$, $u\vee u_0$
is the product map $u\times u_0$.
\item At the stratum $\Sym^{g-1}(\Sigma)\times \Sym^2(E)$,
$u\vee u_0$ is given by
$n_{z}(\psi)$ pseudo-holomorphic spheres which are constant on the first
factor. More precisely, for each $p\in\Delta$ for which
$u(p)=\{z,x_2,...,x_{g}\}$ (where the $x_i\in\Sigma-\{z\}$ are
arbitrary), there is a component of
$u\vee u_0$ mapping into $\Sym^{g-1}(\Sigma)\times \Sym^2(E)$,
consisting of the product of the constant map
$\{x_2,...,x_{g}\}$ with the sphere in $\Sym^{2}(E)$
which passes through $\{z\}\times u_0(p)$.
\item The map $u\vee u_0$ misses all other strata of
$\Sym^{g+1}(\Sigma\vee E)$.
\end{itemize}
As in Theorem~\ref{HolDiskOne:thm:Gluing} of~\cite{HolDisk}, we can
splice to obtain an approximately holomorphic disk $u\#u_0$ (a
triangle) in $\Sym^{g+1}(\Sigma\# E)$. When the connected sum tube is
sufficiently long, the the inverse function theorem can be used to
find the nearby pseudo-holomorphic triangle. The domain belonging to
$u\#u_0$ is clearly given by $\psi\#\psi_0$ described
above. Conversely, by Gromov's compactness (see also
Proposition~\ref{HolDiskOne:prop:GromovCompactness} of~\cite{HolDisk}),
any sequence of pseudo-holomorphic representatives
$u_i\in\pi_2(\x\times x_0,\y\times y_0, \w\times w_0)$ for arbitrarily
long connected sum neck must limit to a pseudo-holomorphic
representative for $\psi'\#\psi_0'$, where
$\cald(\psi_0')-\cald(\psi_0)=k[E]$ for some $0\leq k \leq n_z(\psi)$.
However, since $\pi_2(E)=0$, it follows that $k=0$. Thus, the gluing
map covers the moduli space.
\end{proof}

\begin{prop}
\label{prop:HoClassesCancel}
There are homotopy classes of triangles $\{\psi_k^\pm\}_{k=1}^\infty$
in
$\pi_2(\Theta_{\beta,\gamma},\Theta_{\gamma,\delta},\Theta_{\beta,\delta})$
for the triple-diagram $(\Sigma,\betas,\gammas,\deltas,z)$ satisfying
the following properties:
\begin{eqnarray*}
\Mas(\psi^\pm_k)&=&0, \\
n_z(\psi^\pm_k)&=&\frac{k(k-1)}{2}.
\\
\end{eqnarray*}
Moreover, each triangle in $\pi_2(\Theta_{\beta,\gamma},\Theta_{\gamma,\delta},\Theta_{\beta,\delta})$
is $\SpinC$
equivalent to some $\psi^\pm_k$.
Furthermore, there is a choice of perturbations and
complex structure on $\Sigma$ with the property that
for each $\Psi\in\pi_2(\Theta_{\beta,\gamma},\Theta_{\gamma,\delta},
\x)$ (where $\x\in\Tb\cap\Td$) with $\Mas(\Psi)=0$, we have that
$$
\#\Mod(\Psi)=
\left\{\begin{array}{ll}
\pm 1 & {\text{if $\Psi\in\{\psi^\pm_k\}_{k=1}^\infty$}} \\
0 & {\text{otherwise}}
\end{array} \right..
$$
\end{prop}

\begin{proof}
First observe that the space of $\SpinC$ structures over
$X_{\beta,\gamma,\delta}$ extending a given one on the boundary is
identified with $\Z$. In particular, modulo doubly-periodic domains
for the three boundary three-manifolds, every triangle
$\psi\in\pi_2(\Theta_{\beta,\gamma},\Theta_{\gamma,\delta},\Theta_{\beta,\delta})$
can uniquely be written as $\psi_1+ a[S] + b[\PerDom]$ for some pair
of integers $a$ and $b$, where $\PerDom$ is the generator of the space
of triply-periodic domains: in fact, the integer $a$ is determined by
the intersection number $n_z$, and $b$ can be determined by the signed
number of times the arc in $\beta_g$ obtained by restricting $\psi$ to
its boundary crosses some fixed $\tau\in\beta_g$. For the triangles
$\{\psi_k^\pm\}$ this signed count can be any arbitrary integer, so
these triangles represent all possible $\SpinC$-equivalence classes of
triangles.

The other claims are straightforward in the case where $g=1$. In this
case, the curves $\beta$, $\gamma$, $\delta$ lie in a surface of genus
one, so the holomorphic triangle can be lifted to the complex
plane. Hence, by standard complex analysis, it is smoothly cut out,
and unique.

The fact that $\#\Mod(\psi_k^\pm)=\pm 1$ for higher genus follows from
induction, and the gluing result,
Theorem~\ref{thm:GlueTriangles}. Specifically, if the result is known
for genus $g$, then we can add a new torus $E$ to $\Sigma$ which
contains three curves $\beta_0$, $\gamma_0$, $\delta_0$ which are
small Hamiltonian translates of one another (and the basepoint is
chosen outside the support of the isotopy). The torus $E$ contains a
standard small triangle $\psi_0\in \pi_2(y_0^+,v_0^+,w_0^+)$, for
which it is clear that $\#\ModFlow(\psi_0)=1$.  Gluing this triangle
to the $\{\psi_k^\pm\}$ in $\Sigma$, we obtain corresponding triangles
in $\Sigma\#E$ satisfying all the above hypotheses.

The fact that $\#\Mod(\Psi)=0$ for
$\Psi\not\in\{\psi_k^\pm\}_{k=1}^\infty$ follows similarly, with the
observation that the other moduli spaces of triangles on the
torus are empty.
\end{proof}

We can define the map $$F_1\colon \HFp(Y)\longrightarrow \HFp(Y_0)$$
by summing: $$F_1(\xi)=\sum_{\spinc\in\SpinC(X_{\alpha,\beta,\gamma})}
\pm \Fp{\alpha,\beta,\gamma}(\xi\otimes \theta_{\beta,\gamma},\spinc).$$ 
On the chain level, $F_1$ is induced from a map:
$$f_1([\x,i])=\sum_{\w\in\Ta\cap\Tc}\sum_{\{\psi\in\pi_2(\x,\Theta_{\beta,\gamma},\w)\big|\Mas(\psi)=0\}}
\left(\#\Mod(\psi)\right)\cm [\w,i-n_z(\psi)],$$
where $\#\Mod(\psi)$ is calculated with respect to a particular
choice of coherent orientation system 
(see Proposition~\ref{prop:CompSquareZero} below).  
It is important to note here that the sum on the right
hand side will have only finitely many non-zero elements for each
fixed $\xi\in
\CFp(Y)$. The reason for this is that all the multi-periodic
domains have both positive and negative coefficients.
Similarly, we define 
$$f_2([\x,i])=\sum_{\{\psi\in\pi_2(\x,\Theta_{\gamma,\delta},\w)\big|\Mas(\psi)=0\}}
\left(\#\Mod(\psi)\right)\cm [\w,i-n_z(\psi)],$$
letting $F_2$ be the induced map on homology.

Observe that the maps $f_1$ and $f_2$ preserve the relative
$\Zmod{2}$-grading.  The reason for this is that the parity of the
Maslov index of a triangle $\psi\in\pi_2(\x,\y,\w)$ depends only on
the sign of the local intersection numbers of the $\Ta\cap\Tb$,
$\Tb\cap\Tc$, and $\Ta\cap \Tc$ at $\x$, $\y$, and $\w$.  (Although
each local intersection number is calculated using some choice of
orientations on the three tori, their product is independent of these
choices.)

\begin{prop}
\label{prop:CompSquareZero}
For any coherent system of orientations for 
$Y_0$, we can find coherent systems of orientations for 
the triangles defining $f_1$ and $f_2$ so that the composition $F_2\circ F_1=0$. 
\end{prop}

\begin{proof}
For any system of coherent orientations, associativity, 
together with Proposition~\ref{prop:HoClassesCancel},
can be interpreted as saying that
$$\sum_{s_{\beta,\gamma,\delta}\in S_{\beta,\gamma,\delta}}
f^{\leq 0}_{\beta,\gamma,\delta}(\theta_{\beta,\gamma}\otimes \theta_{\gamma,\delta})
= \sum_{k=1}^\infty
\left[\Theta_{\beta,\delta},-\frac{k(k-1)}{2}\right]\pm
\left[\Theta_{\beta,\delta},-\frac{k(k-1)}{2}\right]$$
(up to an overall sign), as a formal sum.

Of course, if we are using only $\Zmod{2}$ coefficients, the proof is
complete.

More generally, the orientation system for $Y_{\beta,\delta}$ is
chosen so that $\Theta_{\beta,\delta}$ is a cycle. But this leaves the
orientation system over $Y_{\alpha,\gamma}$ unconstrained, and any
choice of such orientation system determines the choice over
$X_{\alpha,\beta,\gamma}$ (up to an overall sign depending on the
$\SpinC$ structure used over $Y_{\alpha,\gamma}$).  Now, the relative
sign appearing above corresponds to the orientation of the triangles
$\psi_k^+$ vs. the triangles $\psi_k^-$ over
$X_{\beta,\delta,\gamma}$, and each such pair of triangles belongs to
different $\delta H^1(Y_{\alpha,\delta})+\delta
H^1(Y_{\beta,\delta})$-orbits for the square
$X_{\alpha,\beta,\gamma,\delta}$. Thus, we can modify the relative
sign at will. We choose it so that the terms pairwise cancel.
\end{proof}

We can choose a one-parameter family of $\gamma$-curves $\gamma_i(t)$
with the property that $\lim_{t\goesto 0}\gamma_i(t)=\beta_i$ for
$i=1,...,g-1$, and $\lim_{t\goesto 0}\gamma_g(t)=\delta_g * \beta_g$
(juxtaposition of curves), and we choose our basepoint $z$ to lie
outside the support of the homotopies $\gamma_i(t)$. We choose another one-parameter family of $\delta$-curves $\delta_i(t)$ for $i=1,...,g-1$ with 
$\lim_{t\goesto 0}\delta_i(t)=\beta_i$. We assume that all $\alpha_i$
are disjoint from the $\beta_g\cap\delta_g$.
Then, if $t$ is
sufficiently small, then there is a natural partitioning of
$\Ta\cap{{\mathbb T}_{\gamma(t)}}$ into two subsets, those which are nearest to points in
$\Ta\cap\Tb$, and those which are nearest to points in
$\Ta\cap\Td(t)$. (See Figure~\ref{fig:LabelledSurgeries} for an
illustration.) Correspondingly, we have a splitting $$\CFp(Y_0)\cong
\CFp(Y)\oplus \CFp(Y_1);$$ or, a short exact sequence of graded groups
$$\begin{CD} 0@>>>\CFp(Y) @>{\iota}>> \CFp(Y_0)
@>{\pi}>>\CFp(Y_1)@>>>0
\end{CD}$$ with splitting 
$$ R\colon \CFp(Y_1)\longrightarrow \CFp(Y_0), $$
where the maps $\iota$, $\pi$, and $R$ are not necessarily chain maps.
Our goal is to construct a short exact sequence as above, which is
compatible with the boundary maps.

% export at 50% 
% Must modify .eps file: 
% move the big polylines to before the ellipses.
% & make sure every ellipse is: 30.000 slw

\begin{prop}
\label{prop:ShortExactSequence}
The map $f_1$ is chain homotopic to a $U$-equivariant
chain map $g_1$ with the
property that
$$
\begin{CD}
0@>>> \CFp(Y)@>{g_1}>>\CFp(Y_0)@>{f_2}>>\CFp(Y_1)@>>> 0.
\end{CD}
$$
is a short exact sequence of chain complexes.
\end{prop}

Theorem~\ref{thm:ExactOne} is a consequence of this proposition: the
associated long exact sequence is the exact sequence of
Theorem~\ref{thm:ExactOne}. 

For the construction of $g_1$, we need the following ingredients:
\begin{itemize}
\item lower-bounded filtrations on the $\CFp(Y)$, $\CFp(Y_0)$, and
$\CFp(Y_1)$,
which are strictly decreasing for the  boundary maps; i.e. each chain
complex is generated by elements with $\partial \xi < \xi$. 
\item an injection $\iota$ and splitting map $R$ as above, both of which
respect the filtrations
\item decompositions of $f_1=\iota + \text{lower order}$ and
$f_2=\pi+\text{lower order}$, where, here, lower order is with respect
to the filtrations. More precisely $\CFp(Y)$ is generated by elements
$\xi$ with the property that $f_1(\xi)-\iota(\xi)<\iota(\xi)$, and
$\CFp(Y_1)$ is generated by elements $\eta$ with
$\eta-f_2\circ R(\eta)<\eta$. 
\item $f_2\circ f_1$ is chain homotopic to zero by a $U$-equivariant 
homotopy
$$ H\colon \CFp(Y)\longrightarrow \CFp(Y_1) $$  which decreases
filtrations, in the sense that $R\circ H < \iota $. 
\end{itemize}

Following Lemma~9 of~\cite{BraamDonaldson}, we define 
a right inverse $R'$ for $f_2$ by
$$R'=R\circ \sum_{k=0}^{\infty} (\Id - f_2\circ R)^{\circ k},$$
and let
$$g_1 = f_1 -
(\partial (R'\circ H) + (R'\circ H) \partial);$$ so that our hypotheses
ensure that $g_1 = \iota + {\text{lower order}}$. It follows that if
$L$ is the left inverse of $\iota$ induced from $R$, then $L\circ g_1 $ is
invertible, as $L\circ g_1 (\xi) = \xi - N(\xi)$, where $N$
decreases filtration (so we can define $$(L\circ g_1)^{-1}(\xi) =
\sum_{k=0}^\infty N^{\circ k}(\xi),$$ as the sum on the right contains
only finitely many non-zero terms for each $\xi\in\CFp(Y)$); thus,
$(L\circ g_1)^{-1}\circ L$ is a left inverse for $g_1$.

A similar argument shows surjectivity of $f_2$, and exactness at the
middle stage (see~\cite{BraamDonaldson}). 

We will use a compatible energy filtration on $\CFp(Y_0)$ defined
presently.  First, fix an $\x_0\in\Ta\cap\Tb$.  If $[\y,j]\in
\CFp(Y_0)$, let $\psi\in\pi_2(\x_0,\Theta_{\beta,\gamma},\y)$ be a
(homotopy class of) triangle, with $n_z(\psi)=-j$. We then define
$$\Filt_{Y_0}([\y,j])=-\Area(\psi).$$ (Note that
$\pi_2(\x_0,\Theta_{\beta,\gamma},\y)$ is non-empty.)  As in
Lemma~\ref{HolDiskOne:lemma:EnergyZero} of~\cite{HolDisk}, 
the topological hypothesis from
Lemma~\ref{lemma:HeegaardDiagrams} allows us to use a volume form on
$\Sigma$ for which every periodic domain for $Y_0$ has zero area:
every periodic domain for $(\Tb,\Td)$, $(\Tb,{\mathbb T}_{\gamma(t)})$
and also the triply-periodic domain for
$(\beta_g,\gamma_g(t),\delta_g)$ has area zero. (For example, we can
start with the area form constructed in the proof of
Lemma~\ref{lemma:HeegaardDiagrams} for the initial $t=0$
$\gamma$-curves, and then move those curves by an exact Hamiltonian isotopy.)
Now, the real-valued function $\Filt_{Y_0}$ on the generators of
$\CFp(Y_0)$ gives the latter group an obvious partial ordering.

We will assume now that the $\gamma_g(t)$ is sufficiently close to the
juxtaposition of $\beta_g$ and $\delta_g$, in the following sense. Let
$\PerDom$ be a triply-periodic domain between $\gamma_g(t)$,
$\beta_g$, and $\delta_g$ which generates the group of such periodic
domains (this is the domain pictured in Figure~\ref{fig:onesurg},
before $\gamma_g$ was isotoped); and for $i=1,...,g-1$, let
$\PerDom_i$ be the doubly-periodic domains with
$\partial\PerDom_i=\beta_i-\gamma_i(t)$.  We let $\epsilon(t)$ be the
sum of the absolute areas of all these periodic domains:
$$\epsilon(t)=\Area\Big(|\cald(\PerDom)|\Big)+\sum_{i=1}^{g-1}\Area\Big(|\cald(\PerDom_i)|\Big),$$
where here the absolute signs denote the unsigned area.  Note that
$\lim_{t\goesto 0}\epsilon(t)=0$.  Also, let $M$ be the minimum of the
area of any domain in
$\Sigma-\alpha_1-...-\alpha_g-\beta_1-...-\beta_g-\delta_g$.  We
choose $t$ small enough that $\epsilon(t)<M/2$.  We assume that the
absolute (unsigned) area of the periodic domain ${\mathcal Q}_i$ with
$\partial ({\mathcal Q}_i)=\beta_i-\delta_i(t)$ agrees with the
absolute area of ${\mathcal P}_i$.

\begin{lemma}
For sufficiently small $t$, the function $\Filt_{Y_0}$ induces a filtration on $\CFp(Y_0)$. In particular, 
$$\partial [\y,j]<[\y,j].$$
\end{lemma}

\begin{proof}
It is important to observe that the area filtration defined above is
indeed well-defined. The reason for this is that if $\psi, \psi'$ are
a pair of homotopy classes in $\pi_2(\x_0,\Theta_{\beta,\gamma},\y)$
with $n_z(\psi)=n_z(\psi')$, then $\cald(\psi)-\cald(\psi')$ is a
triply-periodic domain. It follows from above that it must have total
area zero.

Suppose that we have a pair of generators $[\y,j]$ and $[\y',j']$ which
are connected by a flow $\phi$. If 
$\psi\in\pi_2(\x_0,\Theta_{\beta,\gamma},\y)$ is a class with
$n_z(\psi)=-j$, then, of course,
$\psi+\phi\in\pi_2(\x_0,\Theta_{\beta,\gamma},\y')$ is a class with
$n_z(\psi+\phi)=-j'$; thus,
$\Filt_{Y_0}([\y',j']) - \Filt_{Y_0}([\y,j])=-\Area(\phi)$; but 
$\Area(\phi)>0$, as all of its coefficients are non-negative (and at
least one is positive). 
\end{proof}

The filtration on $\CFp$, together with the data $\iota$, $\pi$, and
$R$, endow $\CFp(Y)$ and $\CFp(Y_1)$ with a filtration as well.

\begin{lemma}
\label{lemma:CompatibleFiltrations}
For $t$ sufficiently small, the orderings induced on $\CFp(Y)$ and
$\CFp(Y_1)$ give filtrations. 
\end{lemma}

\begin{proof}
There is a natural filtration on $Y$, defined by
$\Filt_{Y}([\x,i])=-\Area(\phi)$, where $\phi\in\pi_2(\x_0,\x)$ is the
class with $n_z(\phi)=-i$. This is a filtration, in view of the usual
positivity of holomorphic disks (see
Lemma~\ref{HolDiskOne:lemma:NonNegativity}); indeed, the filtration
decreases by at least $M$ along flows.

The filtration induced by $\Filt_{Y_0}$ and the map $\iota$, defined
by $\Filt_Y'([\x,i])=\Filt_{Y_0}(\iota[\x,i])$ very nearly agrees with
this natural filtration, for sufficiently small $t$.  To see this,
note that there is a unique 
``small'' triangle 
$\psi_0\in
\pi_2(\x,\Theta_{\beta,\gamma},\iota(\x))$ which has non-negative coefficients and is supported inside the
support of ${\mathcal P}+\mathcal P_1+...{\mathcal P}_{g-1}$.
Clearly, $\Area(\psi_0)<\epsilon(t)$, and $n_z(\psi_0)=0$. Now, if
$\phi\in\pi_2(\x_0,\x)$ is the class with $n_z(\phi)=-i$ the
juxtaposition of $\psi_0+\phi\in
\pi_2(\x_0,\Theta_{\beta,\gamma},\iota(\x))$ can be used to calculate
the $Y_0$ filtration of $\iota(\x)$; thus
$|\Filt_{Y}([\x,i])-\Filt_Y'([\x,i])|<\epsilon(t)$. In particular,
since $\Filt_Y$ decreases by at least $M$ along flowlines,
$\Filt_{Y_0}\circ \iota$, too, must decrease along flows.

For $Y_1$, there is another filtration, this one induced by squares.
Given $[\y,i]\in(\Ta\cap\Td)\times \Z^{\geq 0}$, consider
$\varphi\in\pi_2(\x_0,\Theta_{\beta,\gamma},\Theta_{\gamma,\delta},\y)$
with $n_z(\varphi)=-i$,
and let $$\Filt_{Y_1}''([\y,i])=-\Area(\cald(\varphi)).$$ Indeed, if
$M'$ is the minimum area of any
domain in $\Sigma-\alpha_1-...-\alpha_g-\delta_1(t)-...-\delta_{g-1}(t)-\delta_g$, then
$\Filt_{Y_1}''$ decreases by at least $M'$ along each flowline. 
Note that $M'>M-\epsilon(t)$.

Now, we claim that $\Filt_{Y_1}''$ nearly agrees with the filtration
$\Filt_{Y_1}'$ induced by $\Filt_{Y_0}$ and the right inverse $R$:
$\Filt_{Y_1}'([\y,j])=\Filt_{Y_0}(R[\y,j])$. Again, if we let
$\rho(\y)$ denote the point in $\Ta\cap\T\gamma(t)$ closest to
$\y\in\Ta\cap\Td$, there is a unique small triangle $\psi_0\in
\pi_2(\rho(\y),\Theta_{\gamma,\delta},\y)$. If
$\psi\in\pi_2(\x_0,\Theta_{\beta,\gamma},\rho(\y))$ is a triangle with
$n_z(\psi)=-j$ (i.e. used to calculate $\Filt_{Y_0}\circ R$), then,
the juxtaposition $\psi+\psi_0$ is a square which can be used to
calculate $\Filt_{Y_1}''([\y,j])$. But
$|\Area(\psi+\psi_0)-\Area(\psi)|\leq \epsilon(t)$, so since
$\Filt_{Y_1}''$ decreases by at least $M'$ for non-trivial flows, it
follows that $\Filt_{Y_0}\circ R$, too, must
decrease along flows.
\end{proof}

\begin{lemma}
The maps $f_1$ and $f_2$ have the form:
\begin{eqnarray*}
f_1=\iota + {\text{lower order}}, &&
f_2|_{\Image R}=\pi + {\text{lower order}}
\end{eqnarray*}
\end{lemma}

\begin{proof}
The map $f_1([\x,i])$ counts the number of
holomorphic triangles in homotopy classes with
$\psi\in\pi_2(\x,\Theta_{\beta,\gamma},\y)$, with 
$\y\in\Ta\cap {\mathbb T}_{\gamma(t)}$ and $\Mas(\psi)=0$. 
One of these triangles, of course is the canonical small triangle
$\psi_0\in\pi_2(\x,\Theta_{\beta,\gamma},\iota(\x))$. One can
calculate that 
$\#\Mod(\psi_0)=1$. This gives the $\iota$ component of $f_1$. Now,
no other  homotopy class $\psi\in\pi_2(\x,\Theta_{\beta,\gamma},\y)$
with ${\mathcal D}(\psi)\geq 0$
has its domain $\cald(\psi)$ contained inside the support of
${\mathcal P}+\PerDom_1+...\PerDom_{g-1}$; thus, if
$\Mod(\psi)$ is non-empty, then $\Area(\psi)>M-\epsilon(t)>M/2$. Moreover, in
the 
proof of Lemma~\ref{lemma:CompatibleFiltrations}, we saw
that if $\phi\in\pi_2(\x_0,\x)$ is the homotopy class with
$n_z(\phi)=-i$, then
$$|\Filt_{Y_0}(\iota([\x,i]))+\Area(\phi)|<\epsilon(t).$$ But now
$\psi+\phi$ can be used to calculate the filtration 
$\Filt_{Y_0}([\y,i-n_z(\psi)])$. Thus,
$$\Filt_{Y_0}([\y,i-n_z(\psi)])-\Filt_{Y_0}(\iota[\x,i])\leq 
-\Area(\psi)+\epsilon(t)<0.$$

Next, we consider $f_2$.
As before, if $\y\in\Ta\cap\Td$, we let $\rho(\y)\in\Ta\cap{\mathbb T}_{\gamma(t)}$
denote the intersection point closest to $\y$. Suppose that 
$f_2([\rho(\y),i])$ has a non-zero component in $[\w,j]$ with $[\y,i]\neq
[\w,j]$; thus, we have a 
$\psi\in\pi_2(\rho(\y),\Theta_{\gamma,\delta},\w)$ with $n_z(\psi)=i-j$,
which supports a 
holomorphic triangle. Again, $\psi$ cannot be supported inside the
support of ${\mathcal P}+\PerDom_1+...+\PerDom_{g-1}$, so $\Area(\psi)>M/2$.
Fix 
$\psi_\w\in\pi_2(\x_0,\Theta_{\beta,\gamma},\rho(\w))$ (for $\Ta,\Tb,\Tc$)
with $n_z(\psi_\w)=-j$, and
$\psi_{\y}\in
\pi_2(\x_0,\Theta_{\beta,\gamma},\rho(\y))$ with
$n_z(\psi_\y)=-i$. 
Clearly, the
juxtaposition
$\psi_\y +
\psi\in\pi_2(\x_0,\Theta_{\beta,\gamma},\Theta_{\gamma,\delta},\w)$ is
a square whose area must agree with the square
$\psi_\w+\psi_0$, where
$\psi_0\in\pi_2(\rho(\w),\Theta_{\gamma,\delta},\w)$ is the canonical small
triangle, so
$$\Area(\psi_\w)=\Area(\psi_\y)-\Area(\psi_0)+\Area(\psi),$$
and hence $\Filt([\rho(\y),i])>\Filt([\rho(\w),j])$.
\end{proof}

\begin{lemma}
For sufficiently small $t$, there is a null-homotopy $H$ of
$f_2\circ f_1$ satisfying $R\circ H<\iota$. 
\end{lemma}

\begin{proof}
Theorem~\ref{HolDiskOne:thm:Associativity} of~\cite{HolDisk} provides
the null-homotopy $H$: the $[\y,j]$ coefficient of $H[\x,i]$ counts
holomorphic squares $\varphi \in
\pi_2(\x,\Theta_{\beta,\gamma},\Theta_{\gamma,\delta},\y)$ with
$n_z(\varphi)=i-j$. 

Our aim here is to
prove that if the $[\y,j]$ component of $H[\x,i]$ is non-zero then
$\iota[\x,i]>R[\y,j]$.
Now, the
filtration difference between $\iota([\x,i])$ and $R[\y,j]$ is
calculated (to within $\epsilon(t)$) by $\Area(\psi)$, where
$\psi\in\pi_2(\x,\Theta_{\beta,\gamma},\rho(\y))$ has
$n_z(\psi)=i-j$. Adding   the smallest triangle in
$\pi_2(\rho(\y),\Theta_{\gamma,\delta},\y)$ (and hence changing the area
by no more than $\epsilon(t)$), we obtain another square 
$\varphi'\in\pi_2(\x,\Theta_{\beta,\gamma},\Theta_{\gamma,\delta},\y)$
with $n_z(\varphi')=i-j$,
whose area must agree with the area of $\varphi$. Now if $t$ is
sufficiently small ($\epsilon(t)<M/4$), it follows that the filtration
difference between $\iota[\x,i]$ and $R[\y,j]$ is positive.
\end{proof}

\vskip.2cm
\noindent{\bf{Proof of Theorem~\ref{thm:ExactOne}.}}
Theorem~\ref{thm:ExactOne} is now a consequence of the long exact
sequence associated to the short exact sequence from
Proposition~\ref{prop:ShortExactSequence}, with a few final
observations regarding the $\Zmod{2}$ grading. 

Orient the $\alpha_1,...,\alpha_{g}$, the $\beta_1,...,\beta_{g-1}$
arbitrarily (hence inducing orientations on the
$\gamma_1,...,\gamma_{g-1}$ and the $\delta_1,...,\delta_{g-1}$). The
orientation on $\beta_g$ is then forced on us by the requirement that
$$1= \chi(\HFa(Y))=\#(\Ta\cap\Tb),$$ where we orient the tori $\Ta$
and $\Tb$ in the obvious manner. Similarly, the orientation on
$\delta_g$ is forced; indeed, so that $$\delta_g=\beta_g\pm\gamma_g$$
We can orient $\gamma_g$ so that the above sign is positive. It is
then clear with these conventions (by looking at the small triangles)
that $F_1$ preserves the absolute $\Zmod{2}$ grading, while $F_2$
reverses it. It follows then that $F_3$ preserves degree as claimed.

\qed
\vskip.2cm

\subsection{A generalization}
\label{subsec:GenExact}

Let $Y$ be an oriented three-manifold, and let $K\subset Y$ be a knot.
Let $m$ be the meridian of $K$, and let $h\in
H_1(\partial(Y-\nbd{K}))$ be a homology class with $m\cm h=1$ (here,
the torus is oriented as the boundary of the neighborhood of $K$). We
let $Y_h$ denote the three-manifold obtained by attaching a solid
torus to $Y-\nbd{K}$, with framing specified by $h$.

Fix a $\SpinC$ structure $\spinc_0$ over $Y-K$. We let
$$\HFp(Y_h,[\spinc_0])
=\bigoplus_{\{\spinc\big|\spinc|_{Y-K}=\spinc_0\}}\HFp(Y_h,\spinc).$$
We define $\HFp(Y,[\spinc_0])$ similarly.

The following is a direct generalization of Theorem~\ref{thm:ExactOne}
(the case where $Y$ is an integer homology three-sphere, and $h$ is the 
``longitude'' of $K$):

\begin{theorem}
\label{thm:GeneralSurgery}
For each $\SpinC$ structure $\spinc_0$ on $Y-K$, we have the $U$-equivariant
exact sequence:
$$\begin{CD}
... @>>> \HFp(Y,[\spinc_0])@>>> \HFp(Y_h,[\spinc_0])@>>> \HFp(Y_{h+m},
[\spinc_0]) @>>> ...
\end{CD}
$$
\end{theorem}

\begin{cor}
Let $Y$ be an integer homology three-sphere with a knot $K\subset Y$,
and let $Y_n$ be the three-manifold obtained by $n$ surgery on $K$
where $n>0$, then there is a $U$-equivarant long exact sequence $$
\begin{CD}
... @>>> \HFp(Y) @>>> \HFp(Y_n) @>>> \HFp(Y_{n+1}) @>>> ...
\end{CD}
$$
\end{cor}

The proof given in the previous section adapts to this context, after
a few observations. 

Note first that the map from $Y$ to $Y_h$ defined by counting
triangles is naturally partitioned into equivalence classes. To see
the decomposition agrees with what we have stated, we observe
the following. Let $X$ be
the pair-of-pants cobordism connecting $Y$, $Y_h$, and 
$\#^{g-1}(S^2\times S^1)$. The four-manifold obtained by filling the
last component with $\#^{g-1}(D^3\times S^1)$ 
is the cobordism $W_h$ from $Y$ to $Y_h$ obtained by attaching a 
two-handle to $Y$ along $K$ with framing $h$.

Now, $\SpinC$-equivalence classes of triangles for $\Ta$, $\Tb$,
$\Tc$ agree with $\SpinC$ structures on the cobordism $W_h$,
since $\spinc_z(\Theta_{\beta,\gamma})$ is a torsion $\SpinC$
structure over $\#^{g-1}(S^2\times S^1)$ (which extends uniquely over
$\#^{g-1}(D^3\times S^1)$). But two $\SpinC$ structures on $Y$ and
$Y_h$ extend over $W_h$ if and only if they agree on the knot
complement $Y-K$ (thought of as a subset of both $Y$ and $Y_h$).

With this said, the maps $f_1$ and $f_2$ partition according to
$\SpinC$ structures on $Y-K$.

Next, we observe that there are in principle many periodic domains 
for the triple $(\Ta,\Tb,\Tc)$. By twisting normal to the $\alphas$,
however, we can arrange that the triple is admissible. 
%\smargin{talk about 4-periodic domain}
By choosing the
volume form on $\Sigma$ appropriately, we can arrange that they all
have zero signed area. 

We can define the filtrations as before. Fix any $\x_0\in\Ta\cap\Tb$ 
so that $\spinc_z(\x_0)$ restricts to $\spinc_0$ on $Y-K$. 
The triangle connecting
$\x_0$, $\Theta_{\beta,\gamma}$ and any intersection point $\y\in
\Ta\cap \Tc$ with $\spinc_z(\y)|Y-K=\spinc_0$ is guaranteed to exist,
since the corresponding $\SpinC$ structures extend over $W_h$. The
area the of the domain of any such triangle can be used to define
$\Filt_{Y_h}([\y,i])$. The proof given before, then, applies.

\subsection{Fractional Surgeries}
\label{subsec:FracSurg}

There are other directions to generalize
Theorem~\ref{thm:ExactOne}. We consider presently the case of fractional ($1/q$) 
surgeries on an integral homology three-sphere. 

Let $Y$ be an integer homology three-sphere, and $K\subset Y$ be a
knot. Let $Y_0$ be the manifold obtained by zero-surgery on $K$, and
let $Y_{1/q}$ be obtained by $1/q$ surgery on $K$, where $q$ is a
positive integer.

We fix a representation $$H^1(Y;\Z) \longrightarrow \Zmod{q}$$ taking
generators to generators, and let $$\uHFp(Y_0,\Zmod{q})\cong
\bigoplus_{\spinc\in\SpinC(Y_0)}\uHFp(Y_0,\spinc)$$ denote the
corresponding homology group with twisted coefficient ring (in the
sense of Section~\ref{sec:TwistedCoeffs}).

\begin{theorem}
\label{thm:ExactFrac}
Let $Y$ be an integral homology three-sphere and let $q$ be a positive integer.
Then, there is a $U$-equivariant exact sequence
$$\begin{CD}
... @>>>  \uHFp(Y_0;\Zmod{q})@>>> \HFp(Y_{1/q})@>>> \HFp(Y)
@>>> ...
\end{CD}$$
\end{theorem}

The proof of Lemma~\ref{lemma:HeegaardDiagrams} in the present context
gives us a generalized pointed Heegaard diagram
$(\Sigma,\alphas,\betas,\gammas,\deltas,z)$
with the property that:
\begin{itemize}
\item the Heegaard diagrams $(\Sigma,\alphas,\betas)$,
$(\Sigma,\alphas,\gammas)$, and $(\Sigma,\alphas,\deltas)$ describe
$Y$, $Y_0$, and $Y_{1/q}$ respectively,
\item for each $i=1,...,g-1$, the curves $\beta_i$, $\gamma_i$, and
$\delta_i$ are small isotopic translates of one another, each pairwise
intersecting in a pair of canceling transverse intersection points
\item the curve $\delta_g$ is isotopic to the juxtaposition
of $\beta_g$ with the 
$q$-fold juxtaposition of $\gamma_g$.
\end{itemize}

We can think concretely about $\uCFp(Y_0;\Zmod{q})$ as
follows. Let $\zeta=e^{\frac{2\pi i}{q}}$, and fix a reference point
$\tau\in\gamma_g$, which we choose to be disjoint from all the other
$\alphas$, $\betas$, and $\deltas$. This gives rise to a codimension-one submanifold
$$V=\gamma_1\times...\times\times \gamma_{g-1}\times \{\tau\}\subset \Tc.$$
Then, $\uCFp(Y_0;\Zmod{q})$ is generated over
$\Z$ by the basis $[\x,i]\otimes \zeta^j$ where of course, $\x$ is an 
intersection point $\x\in\Ta\cap\Tc$ in the appropriate equivalence
class, $i$ is a non-negative integer, and $j\in\Zmod{q}$. The boundary
map then is given by
\begin{equation}
\label{eq:DefBoundaryTwist}
\partial\left([\x,i]\otimes
\zeta^j\right)=\sum_{\y\in\Ta\cap\Tc}\sum_{\{\phi\in\pi_2(\x,\y)\big|\Mas(\phi)=1\}}
\left(\#\UnparModFlow(\phi)\right)\cm[\y,i-n_z(\phi)]\otimes \zeta^{j+\#(V\cap
\partial_{\gamma}(\phi))}.
\end{equation}
The quantity $V\cap\partial_{\gamma}(\phi)$ is the intersection number
between the codimension-one submanifold $V\subset \Tc$ with
the path in $\Tc$ obtained by restricting $\phi$ to the
appropriate edge. 

Again, we let $v_g$ be the intersection point between $\delta_g$ and
$\gamma_g$. We now have $q$ different intersection points between
$\delta_g$ and $\beta_g$, of which we choose one, labelled $w_g$, in
the following Proposition~\ref{prop:OneQHoClassesCancel}. We will have
no need for the $q-1$ other intersection points. Let
$\Theta_{\beta,\gamma}$, $\Theta_{\gamma,\delta}$, and
$\Theta_{\beta,\delta}$ be as before.

As in Proposition~\ref{prop:EltsAreCycles}, if we let
$\theta_{\beta,\delta}=[\Theta_{\beta,\delta},0]$, then
$\theta_{\beta,\delta}$ is a cycle in $\CFinfty(\Tb,\Td)$. Note that
the three-manifold described by the pair $(\Sigma,\betas,\deltas)$ is
now a sum $L(q,1)\#\left(\#_{i=1}^{g-1}(S^1\times S^2)\right)$ (where
$L(q,1)$ is a lens space).  

\begin{prop}
\label{prop:OneQHoClassesCancel}
For an appropriate choice $w_g\in \beta_g\cap\delta_g$ 
for $\beta_g$ with
$\delta_g$,
there are homotopy classes of triangles 
$\{\psi^\pm_k\}_{k=1}^\infty\in\pi_2(\Theta_{\beta,\gamma},
\Theta_{\gamma,\delta},\Theta_{\beta,\delta})$
satisfying the following properties (for each $k$):
\begin{eqnarray*}
\Mas(\psi^\pm_k)&=&0, \\
n_z(\psi^+_k)&=&n_z(\psi^-_k), \\
n_z(\psi^+_k)&<&n_z(\psi^+_{k+1}), \\
\end{eqnarray*}
Moreover, each triangle in
$\pi_2(\Theta_{\beta,\gamma},\Theta_{\gamma,\delta},\Theta_{\beta,\delta})$
is $\SpinC$ equivalent to some $\psi^\pm_k$. Also, the congruence
class modulo $q$ of the intersection number
$\#(V\cap\partial_\gamma(\psi))$ is independent of the choice of
$\psi\in\pi_2(\Theta_{\beta,\gamma},\Theta_{\gamma,\delta},
\Theta_{\beta,\delta})$.
Furthermore, there is a choice of perturbations and complex structure with the property that 
for each $\Psi\in\pi_2(\x,\Theta_{\gamma,\delta},
\Theta_{\beta,\delta})$ (where $\x\in\Tb\cap\Tc$) with $\Mas(\Psi)=0$, we have that
$$
\#\Mod(\Psi)=
\left\{\begin{array}{ll}
\pm 1 & {\text{if $\Psi\in\{\psi^\pm_k\}_{k=1}^\infty$}} \\
0 & {\text{otherwise}}
\end{array} \right..
$$
\end{prop}

\begin{proof}
The proof follows along the lines of
Proposition~\ref{prop:HoClassesCancel}. In this case, letting
$\PerDom$ be the generating periodic domain in the torus, we have that
$$\partial \PerDom = \beta_g + q \gamma_g - \delta_g.$$ We must choose
$w_g$ so that it is the $\beta_g$-$\delta_g$ corner point for the
domain containing the basepoint $z$.  Note that $\partial\PerDom$
meets the reference point $\tau\in\gamma$ with multiplicity $q$. This
proves the $q$ independence of the intersection number
$\#(V\cap\partial_\gamma(\psi))$ of the choice of
$\psi\in\pi_2(\Theta_{\beta,\gamma},\Theta_{\gamma,\delta},\Theta_{\beta,\delta})$.
(See Figure~\ref{fig:OneQSurgery}.)

\begin{figure}
\mbox{\vbox{\epsfbox{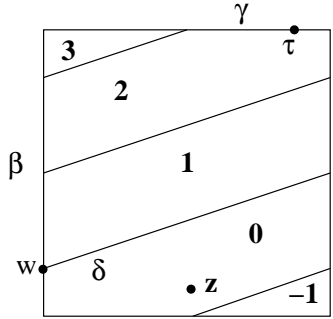}}}
\caption{\label{fig:OneQSurgery}
The triply-periodic domain in the torus 
relevant for $1/q$ surgery, with $q=3$.}
\end{figure}
\end{proof}

Our choice of basepoint $z$ and the intersection point
$\Theta_{\beta,\delta}$, from the above proposition give us a $\SpinC$ structure
$\spinct_{\beta,\delta}\in\SpinC(L(q,1)\#\left(\#_{i=1}^{g-1}(S^1\times
S^2)\right))$.

We consider the chain map
$$f_2\colon \uCFp(Y_0,\Zmod{q})\longrightarrow \CFp(Y_{1/q})$$
defined by
$$f_2(\xi)=\sum_{\{\spinc\in\SpinC(X_{\alpha,\gamma,\delta})}
\ufp{\alpha,\gamma,\delta}(\xi\otimes
\theta_{\gamma,\delta},\spinc).$$
In the present context,
$$
\ufp{\alpha,\gamma,\delta}\left([\x,i]\otimes \zeta^k \otimes [\y,j];\spinc\right)
= \nonumber \\
\sum_{\w\in\Ta\cap\Td} 
\sum_{\{\psi\in\pi_2(\x,\y,\w)\big| \#(V\cap
\partial_\gamma \psi)=-k, \spinc_z(\psi)=\spinc\}}
\left(\#\ModFlow(\psi)\right)\cdot[\w,i+j-n_z(\psi)].
$$

We define
$$f_3\colon \CFp(Y_{1/q})\longrightarrow \CFp(Y)$$ by
$$f_3(\xi)=\sum_{\{\spinc\in\SpinC(X_{\alpha,\delta,\beta})\big|
\spinc|Y_{\beta,\delta}=\spinct_{\beta,\delta}\}}
f^+_{\alpha,\delta,\beta}(\xi\otimes \theta_{\beta,\delta},\spinc).$$

This gives us maps:
$$\begin{CD}
\uCFp(Y_0,\Zmod{q}) @>{f_2}>> \CFp(Y_{1/q}) @>{f_3}>> \CFp(Y).
\end{CD}
$$

It follows, once again, from associativity, together with the
Proposition~\ref{prop:OneQHoClassesCancel} that the maps on homology
$F_3\circ F_2=0$. Note that the chain homotopy evaluated on
$\zeta^k\times[\x,i]$ is constructed by counting squares in
$\varphi\in\pi_2(\x,\Theta_{\gamma,\delta},\Theta_{\delta,\beta},\y)$
with $V\cap \partial_\gamma(\varphi)=-k$.

\begin{figure}
\mbox{\vbox{\epsfbox{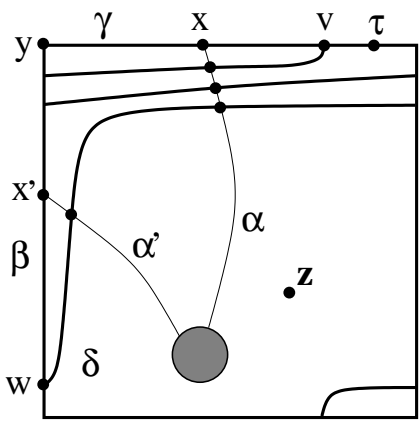}}}
\caption{\label{fig:OneQSurgPert}
The analogue of Figure~\ref{fig:LabelledSurgeries}, only for $1/q$
surgery with $q=3$. We have pictured here only the part of the
surface taking place in the final torus summand, and correspondingly
dropped the $g$ subscripts. There are two $\alpha$-curves crossing the
region here, labelled $\alpha$ and $\alpha'$: the first of 
these 
meets $\gamma$ at $x$, the second meets $\beta$ at $x'$. Observe the three
intersection points of $\alpha\cap\delta$ and the 
intersection
point of $\alpha'\cap\delta$ corresponding to $x$ and $x'$ respectively.}
\end{figure}

We homotope the $\delta$-curve to the juxtaposition of the $\beta_g$ with the 
$q$-fold juxtaposition of $\gamma_g$.
This gives a short exact sequence of
graded groups
$$
\begin{CD}
0@>>> \uCFp(Y_0,\Zmod{q}) @>{\iota}>> \CFp(Y_{1/q}) @>{\pi}>> 
\CFp(Y) @>>> 0.
\end{CD}
$$ 
To
see the inclusion, note that each intersection point $\x$ of
$\Ta\cap\Tc$ corresponds to $q$ distinct intersection points between
$\Ta\cap\Td$, labelled $(\x_1,...,\x_q)$. For each of these
intersection points, there is a unique smallest triangle $u_1,...,u_q$,
with $u_i\in\pi_2(\x,,\Theta_{\gamma,\delta},\x_j)$. We claim that the
$q$ integers $\# (V\cap u_i)$ each lie in different congruence classes
modulo $q$. This gives the inclusion.
To see surjection, note that each intersection point of $\x'\in
\Ta\cap\Tb$ gives rise to a unique intersection point $\rho(\x')$ between 
$\Ta\cap\Td$, which can be joined by a small triangle in $\pi_2(\rho(\x'),\Theta_{\beta,\delta},\x')$.
(See Figure~\ref{fig:OneQSurgPert} for an illustration.)

% export at 50% 

With this said, then, the energy filtration is defined as before,
calculating the energy of classes
$\psi\in\pi_2(\x_0,\Theta_{\gamma,\delta},\y)$. Thus we obtain the
required long exact sequence.

\subsection{$\HFa$}
\label{subsec:SurgeriesHFa}

Let $Y$ be an oriented three-manifold, $K\subset Y$ be a knot, and
$\spinc_0$ be a fixed $\SpinC$ structure over $Y-K$.

\begin{theorem}
\label{thm:GeneralSurgeryHFa}
For each $\SpinC$ structure $\spinc_0$ on $Y-K$, we have the exact sequence:
$$\begin{CD}
... @>>> \HFa(Y,[\spinc_0])@>>> \HFa(Y_h,[\spinc_0])@>>> \HFa(Y_{h+m},
[\spinc_0]) @>>> ...
\end{CD}
$$
\end{theorem}

Similarly, we have:

\begin{theorem}
\label{thm:ExactFraca}
Let $Y$ be an integral homology three-sphere and let $q$ be a positive integer.
Then, there is a $U$-equivariant exact sequence
$$\begin{CD}
... @>>>  \uHFa(Y_0;\Zmod{q})@>>> \HFa(Y_{1/q})@>>> \HFa(Y)
@>>> ...
\end{CD}$$
\end{theorem}

For the proofs of these results, 
Proposition~\ref{prop:OneQHoClassesCancel} (or
Proposition~\ref{prop:HoClassesCancel}, for the case of
$+1$-surgeries) is replaced by the comparatively simpler:

\begin{prop}
\label{prop:OneQHoClassesCancela}
There are two homotopy classes of triangles  
$\psi^+$ and $\psi^-$ in
$\pi_2(\Theta_{\beta,\gamma},\Theta_{\gamma,\delta},\Theta_{\beta,\delta})$
with
\begin{eqnarray*}
\Mas(\psi^\pm)&=&0, \\
n_z(\psi^\pm)&=&0, \\
\#(\partial_\gamma \psi^+)&=&\#(\partial_\gamma \psi^-) + q.
\end{eqnarray*}
Indeed, these are the only two triangles with $\cald(\psi)\geq 0$ and
$n_z(\psi)=0$. Also, each moduli space consists of a single, smooth isolated point. 
\end{prop}

\begin{proof}
This now follows directly from the picture in the torus. In
particular, in the present case, there is no need for
Theorem~\ref{thm:GlueTriangles}.
\end{proof}

\vskip.2cm
\noindent{\bf{Proof of Theorems~\ref{thm:GeneralSurgeryHFa} and
\ref{thm:ExactFraca}.}}
The proofs here are now obtained by copying the earlier proofs for 
$\HFp$, with the obvious notational changes.
\qed

\subsection{Integer surgeries}
\label{subsec:IntSurg}

Another generalization of
Theorem~\ref{thm:ExactOne} involves integer surgeries.

Let $Y$ be an integer homology three-sphere, and $K\subset Y$ be a
knot. Let $Y_0$ be the manifold obtained by zero-surgery on $K$, and
let $Y_p$ be obtained by $+p$ surgery on $K$, where $p$ is a positive
integer. 

\begin{theorem}
\label{thm:ExactP}
There is a surjective map $Q\colon \SpinC(Y_0)\longrightarrow \SpinC(Y_p)$ with the property
that for each $\SpinC$ structure $\spinct\in\SpinC(Y_p)$, we have a $U$-equivariant exact sequence
$$\begin{CD}
... @>{F_1}>>  \HFp(Y_0,[\spinct])@>{F_2}>> \HFp(Y_{p},\spinct)@>{F_3}>> \HFp(Y)
@>>> ...,
\end{CD}
$$
where 
$$\HFp(Y_0,[\spinct])=\bigoplus_{\{\spinct_0\big| Q(\spinct_0)=\spinct\}}
\HFp(Y_0,\spinct_0).$$
Moreover, $F_3$ preserves $\Zmod{2}$ degree, chosen so that
$$\chi(\HFa(Y_p,\spinct))=\chi(\HFa(Y))=1.$$
In particular, there is a $U$-equivariant exact sequence 
$$\begin{CD}
... @>>> \HFp(Y_0) @>>> \HFp(Y_{p}) @>>> \bigoplus_{i=1}^p\HFp(Y) @>>> ...,
\end{CD}
$$
\end{theorem}

\begin{remark}
Indeed, a modification of the following proof can also be given to construct an exact sequence
$$\begin{CD}
... @>{F_2}>>  \HFp(Y) @>{F_3}>> \HFp(Y_{-p},\spinct)@>{F_1}>> 
\HFp(Y_0,[\spinct])
@>>> ...,
\end{CD}
$$
where $F_3$ preserves the $\Zmod{2}$ degree.

In another direction, Theorem~\ref{thm:ExactP} readily generalizes to
the case where $Y$ is not an integral homology sphere. For example, if
$K$ is a null-homologous knot in $Y$, there is still a notion of
integral surgery, and we obtain sequences as above, only now there is
one for each fixed $\SpinC$ structure over $Y$.
\end{remark}

\begin{proof}
This time, the curve $\delta_g$ is isotopic to the juxtaposition of
the $p$-fold juxtaposition of $\beta_g$ with the $\gamma_g$.

Now, we have $p$ different intersection points between $\delta_g$ and
$\gamma_g$. We choose one (so that the analogue of
Proposition~\ref{prop:OneQHoClassesCancel} holds, for our given choice of
basepoint), and label it $v_g$. We will have no need for the remaining
$p-1$ intersection points. Let $w_g$ denote the intersection point
between $\beta_g$ and $\delta_g$, and let $\Theta_{\beta,\gamma}$,
$\Theta_{\gamma,\delta}$, and $\Theta_{\beta,\delta}$ be as
before. We have a corresponding $\SpinC$ structure
$\spinct_{\gamma,\delta}$ corresponding to $\Theta_{\gamma,\delta}$.

If $\spinct'\in\SpinC(Y_0)$, there is a unique $\SpinC$ structure
$\spinct\in\SpinC(Y_p)$ with the property that there is a $\SpinC$
structure $\spinc$ on $X_{\alpha,\gamma,\delta}$ with
$\spinc|Y_0=\spinct'$,
$\spinc|Y_{\gamma,\delta}=\spinct_{\gamma,\delta}$, and
$\spinc|Y_{\alpha,\delta}=\spinct$. We let $Q(\spinct')=\spinct$.

Fix a $\SpinC$ structure $\spinct$ over $Y_{p}$. 
We consider the chain map
$$f_2\colon \CFp(Y_0)\longrightarrow \CFp(Y_{p},\spinct)$$
defined by
$$f_2(\xi)=\sum_{\{\spinc\in\SpinC(X_{\alpha,\beta,\delta})\big|
\spinc|Y_{\alpha,\delta}=\spinct,~\spinc|_{\gamma,\delta}=\spinct_{\gamma,\delta}\}}
f^+_{\alpha,\gamma,\delta}(\xi\otimes
\theta_{\gamma,\delta},\spinc).$$

We define $f_3$ as follows. 
Consider
$$f_3(\xi)=\sum_{\{\spinc\in\SpinC(X_{\alpha,\delta,\beta})\big|
\spinc|Y_{p} = \spinct\} }
f^+_{\alpha,\delta,\beta}(\xi\otimes \theta_{\beta,\delta}).$$

This gives us maps:
$$\begin{CD}
\CFp(Y_0,[\spinct]) @>{f_2}>> \CFp(Y_{p},\spinct) @>{f_3}>> \CFp(Y).
\end{CD}$$

It follows once again from associativity, together with 
the analogue of Proposition~\ref{prop:OneQHoClassesCancel},
that $F_3\circ F_2=0$.

We homotope the $\delta$-curve to the juxtaposition of the $p$-fold
multiple of $\beta_g$ with $\gamma_g$. This gives a short exact sequence of
graded groups
$$
\begin{CD}
0@>>> \CFp(Y_0,[\spinct])@>{\iota}>> \CFp(Y_p,\spinct) @>{\pi}>> \CFp(Y)
@>>> 0.
\end{CD}
$$ 
The inclusion follows as before: each intersection point $\x$ of
$\Ta\cap\Tc$ corresponds a unique intersection point between
$\Ta\cap\Td$, which can be canonically connected by a small triangle. 
To see surjection, note that each intersection point of $\y\in
\Ta\cap\Tb$ gives rise to $p$ different intersection points between
$\Ta\cap\Td$, which we label $(\y_1,...,\y_p)$.  Note, however, that
$\epsilon(\y_i,\y_j)=(i-j)\PD[\beta_g^*]$. Now, $\PD[\beta_g^*]\in
H^2(Y_{p})$ is a generator, so there will always be a unique induced
intersection point representing the $\SpinC$ structure $\spinct$ over $Y_p$.  
The rest follows as before.
\end{proof}

\subsection{$+1$ surgeries for twisted coefficients}
\label{subsec:SurgeriesTwisted}

There is also surgery exact sequence for $+1$ surgeries which
uses twisted coefficients. 

For simplicity, we state it in the case where we begin with a
three-manifold $Y$ which is an integer homology sphere. In that case,
if we let $T$ be a generator for $H^1(Y_0;\Z)$, then we can think of
$\Z[H^1(Y_0;\Z)]$ as $\Z[T,T^{-1}]$.  Given any $\Z[U]$ module $M$,
let $M[T,T^{-1}]$ denote the induced module over $\Z[U,T,T^{-1}]$.

\begin{theorem}
\label{thm:ExactOneTwist}
There is a $\Z[U,T,T^{-1}]$-equivariant long exact sequence:
$$\begin{CD} 
... @>>> \HFp(Y)[T,T^{-1}] @>{\uFp{1}}>>
\uHFp(Y_0)@>{\uFp{2}}>> \HFp(Y_1)[T,T^{-1}] @>{\uFp{3}}>> ...
\end{CD}
$$
\end{theorem}

We will think of $\uHFp(Y_0)$ like we did in
Subsection~\ref{subsec:FracSurg}: we fix a reference point
$\tau\in\gamma_g$, and let the boundary map record, in the power of
$T$, the multiplicity with which $\phi$ meets $\tau$ along its
boundary, as in Equation~\eqref{eq:DefBoundaryTwist} (with the
difference that now we use a formal variable $T$ rather than a root of
unity $\zeta$).

We will similarly use a reference point $\tau'\in\delta_g$, again
defining the boundary map for $Y_1$ which records the intersection
with $\tau'$ in the power of $T$, to obtain a chain complex for $Y_1$,
which we write as: $\CFp(Y_1,\Z[T,T^{-1}])$. Note that (by contrast
with the case of $Y_0$) this has little effect on the
homology. Indeed, it is easy to construct an isomorphism of chain
complexes (over $\Z[U,T,T^{-1}]$): $$\CFp(Y_1)\otimes_\Z \Z[T,T^{-1}]
\cong \CFp(Y_1,\Z[T,T^{-1}]).$$ Moreover, it is clear that $$H_*(\CFp(Y_1)\otimes_\Z
\Z[T,T^{-1}]) \cong \HFp(Y_1)\otimes_{\Z}\Z[T,T^{-1}].$$ However, this
device will be convenient in constructing the chain maps.

We choose $\tau'$ to lie on the boundary of $\psi^-$ and $\tau$ to lie
on the boundary of $\psi^+$ (where $\psi^\pm=\psi^\pm_1$ from
Proposition~\ref{prop:HoClassesCancel}), and let $V$, $V'$ be the corresponding
codimension one subsets of $\Tc$ and $\Td$ respectively. We then let 
$$\ufp{1}([\x,i])=
\sum_{\w\in\Ta\cap\Tc}
\sum_{\{\psi\in\pi_2(\x,\Theta_{\beta,\gamma}.\w)\big|\Mas(\psi)=0\}}
c(\x,\w,\psi)\cm [\w,i-n_z(\psi)],$$
and
$$\ufp{2}([\x,i])=
\sum_{\w\in\Ta\cap\Td}
\sum_{\{\psi\in\pi_2(\x,\Theta_{\gamma,\delta},\w)\big|\Mas(\psi)=0\}}
c(\x,\w,\psi)\cm [\w,i-n_z(\psi)];$$
where in both cases $c(\x,\w,\psi)\in \Z[T,T^{-1}]$ is given by
$$c(\x,\w,\psi)=(\#\Mod(\psi))\cm 
\left(T^{\#(\partial_\gamma\psi\cap
V)+\#(\partial_\delta\psi\cap V')}\right).$$

We have the following analogue of
Proposition~\ref{prop:CompSquareZero}:

\begin{prop}
\label{prop:CompSquareZeroTwist}
The composition $\uFp{2}\circ \uFp{1}=0$. 
\end{prop}

\begin{proof}
Observe that for the homotopy classes $\{\psi^\pm_k\}_{k=1}^\infty$ from
Proposition~\ref{prop:HoClassesCancel}, we have that
$$\#(\partial_\beta\psi^+_k\cap V)+\#(\partial_\delta\psi^+_k\cap V')
= \#(\partial_\beta\psi^-_k\cap V)+\#(\partial_\delta\psi^-_k\cap V')
= 1$$
This implies that the formal sum
\begin{eqnarray*}
\sum_{\spinc_{\beta,\gamma,\delta}\in S_{\beta,\gamma,\delta}}
\ufleq{\beta,\gamma,\delta}(\theta_{\beta,\gamma}\otimes\theta_{\gamma,\delta},
\spinc_{\beta,\gamma,\delta})
&=&
\sum_{k=1}^\infty
T\otimes \left(\left[\Theta_{\beta,\delta},-\frac{k(k-1)}{2}\right]-
\left[\Theta_{\beta,\delta},-\frac{k(k-1)}{2}\right]\right) \\
&=&
0.
\end{eqnarray*}
Thus, the proof follows from associativity as before.
\end{proof}

\vskip.2cm
\noindent{\bf{Proof of Theorem~\ref{thm:ExactOneTwist}.}}
With Proposition~\ref{prop:CompSquareZeroTwist} replacing
Proposition~\ref{prop:CompSquareZero}, the proof proceeds as the proof
of Theorem~\ref{thm:ExactOne}.
\qed

We have also the generalization for integer surgeries:

\begin{theorem}
\label{thm:ExactPTwist}
Let $Y$ be an integral homology three-sphere, let $K\subset Y$ be a 
knot in $Y$, and fix a positive integer $p$.
For each $\SpinC$ structure $\spinct\in\SpinC(Y_p)$, 
we have a $\Z[U,T,T^{-1}]$-equivariant exact sequence
$$\begin{CD}
... @>{F_1}>>  \uHFp(Y_0,[\spinct])@>{F_2}>> 
\HFp(Y_{p},\spinct)[T,T^{-1}]@>{F_3}>> \HFp(Y)[T,T^{-1}]
@>>> ...,
\end{CD}$$
where
$$\uHFp(Y_0,[\spinct])=\bigoplus_{\{\spinct_0\big| Q(\spinct_0)=\spinct\}}
\uHFp(Y_0,\spinct_0),$$
using the map
$Q\colon \SpinC(Y_0)\longrightarrow \SpinC(Y_p)$ be the map
from Theorem~\ref{thm:ExactP}.
\end{theorem}

\begin{proof}
Combine the refinements from Theorem~\ref{thm:ExactP} with those of Theorem~\ref{thm:ExactOneTwist}.
\end{proof}

\section{Calculation of $\HFinf$}
\label{sec:HFinfty}

The main result of the present section is the complete calculation of
$\uHFinf(Y)$ purely in terms of the homological data
of $Y$. We also give the following 
similar calculation of
$\HFinf(Y)$ when $b_1(Y)\leq 2$. 
We start with the latter construction, establishing the following:

\begin{theorem}
\label{thm:HFinfGen}
Let $Y$ be a closed, oriented three-manifold with $b_1(Y)\leq 2$.  
Then, there is
an equivalence class of orientation system over $Y$ with the following
property.  If $\spinc_0$ is torsion, then $$\HFinfty(Y;\spinc_0)\cong
\Z[U,U^{-1}]\otimes_\Z \Wedge^* H^1(Y;\Z)$$ as a $\Z[U]\otimes_\Z
\Wedge^*(H_1(Y;\Z)/\Tors)$-module.  Furthermore, if $\spinc$ is not
torsion, $$\HFinfty(Y;\spinc)\cong
\left(\Z[U]/{U^n-1}\right) \otimes_\Z \Wedge^* c_1(\spinc)^\perp,$$ where
$c_1(\spinc)^\perp\subset H^1(Y;\Z)$ is the subgroup pairing trivially
with $c_1(\spinc)$, and $n={\mathfrak d}(\spinc)/2$.
\end{theorem}

\begin{remark}
Of course, in the above statement, we think of the usual cohomology
$H^1(Y;\Z)$ (with constant coefficients); but it will be apparent from
the proof that for each choice of locally constant $\Z$ coefficient
system, we obtain an orientation system for $\HFinf$ for which the
analogous isomorphism holds: this gives an identification between
locally constant $\Z$ coefficient systems over $Y$ and equivalence
classes of orientation system over $Y$. 
\end{remark}

The proof in some important special cases is given in
Subsection~\ref{subsec:HFinfSmallY}, and the general case is proved in
Subsection~\ref{subsec:HFinfGen}. We give also the ``twisted''
analogue in Subsection~\ref{subsec:HFinfTwist} which holds for
arbitrary $b_1(Y)$.

The theorem describes the module structure of $\HFp(Y,\spinc_0)$ in
sufficiently large degree, when $\spinc_0$ is a torsion $\SpinC$
structure and $b_1(Y)\leq 2$. It also allows us to pay off several
other debts: first, it allows us to define an absolute $\Zmod{2}$
grading on the homology groups; then, combined with the discussion of
Section~\ref{sec:EulerCharacteristic}, it allows us to relate
$\chi(\HFm(Y,\spinc))$ with Turaev's torsion in
Subsection~\ref{subsec:EulerHFm} (though an alternative calculation
could also be given by modifying directly the discussion in
Section~\ref{sec:EulerCharacteristic}). It also allows us to extend
the Euler characteristic calculations for $\HFp$ to the case where the
$\SpinC$ structure is torsion,
c.f. Subsection~\ref{subsec:TruncEuler}.  Finally, the result allows
us to identify a ``standard'' orientation system for $Y$: the one for
which Theorem~\ref{thm:HFinfGen} holds, with the usual $H^1(Y;\Z)$ on
the right-hand-side. (This justifies our practice of dropping the
coefficient system from the notation for $\HFinf$, and the other
related groups.) Since the analogue of Theorem~\ref{thm:HFinfGen} in
the twisted case (Theorem~\ref{thm:HFinfTwist}) holds without
restriction on the Betti numbers of $Y$, it can be used to identify a
canonical coherent system of orientations for any oriented
three-manifold $Y$.

\subsection{$\HFinfty(Y)$ when $H_1(Y;\Z)=0$ or $\Z$}
\label{subsec:HFinfSmallY}

\begin{theorem}
\label{thm:HFinfZHS}
Theorem~\ref{thm:HFinfGen} 
holds when $Y$ is an integer homology three-sphere;
i.e. over $\Z$, $\HFinfty(Y)$ is freely generated by generators
$y_i$ for $i\in\Z$, with $U y_i = y_{i-1}$.
\end{theorem}

\begin{theorem}
\label{thm:HFinfBigger}
Theorem~\ref{thm:HFinfGen} holds when the three-manifold in question
$Y_0$ satisfies $H_1(Y_0)\cong \Z$. More concretely, let $H\in H^2(Y_0;\Z)$ be a
generator, and let $\spinc_0$ denote the $\SpinC$ structure with
trivial first Chern class. Then if $\spinc=\spinc_0\pm n\cm H$ with
$n>0$, then $\HFinf(Y_0,\spinc)$ is freely generated by generators
$x_i$ for $i\in 1,...,n$, with $U x_i=x_{i-1}$, $U x_1= x_n$.  Moreover,
$\HFinf(Y_0,\spinc_0)$ is freely generated by generators $x_i$, $y_i$
for $i\in\Z$, with $U y_i = y_{i-1}$, $U x_i = x_{i-1}$ and
$\gr(x_i,y_i)=1$; also, $\PD[H]\cm x_i=y_i$.
\end{theorem}

The main ingredient in the proof of the above results is the following:

\begin{prop}
\label{prop:HFinfProp}
Let $Y$ be an integer homology three-sphere, and $K\subset Y$ be a knot, then
there is an identification:
$$\HFinf(Y_0,\spinc)\cong \HFinf(Y,\spinc_0)/(U^n-1),$$
where $Y_0$ is the three-manifold obtained by zero-surgery on $K$, 
and where the divisibility of $c_1(\spinc)$ is $2n$.
\end{prop}

This is proved in several steps.

We start with a Heegaard diagram $(\Sigma,\alphas,\betas,z)$
describing $Y_0$, with the property that
$(\Sigma,\{\alpha_2,...,\alpha_g\}, \{\beta_1,...,\beta_g\})$
describes the knot complement. Let $\gamma$ be a curve which
intersects $\alpha_1$ once and is disjoint from
$\{\alpha_2,...,\alpha_g\}$, so that
$(\Sigma,\{\gamma,\alpha_2,...,\alpha_g\}, \{\beta_1,...,\beta_g\})$
represents $Y$. Indeed, we let $\gamma_2,...,\gamma_g$ be small isotopic translates of $\alpha_2,...\alpha_g$,
with $\gamma_i\cap\alpha_i$ for $i=2,...,g$ consisting of a canceling pair of points $w^\pm_i$.
Such a diagram can always be found (compare 
Lemma~\ref{lemma:HeegaardDiagrams}). 
We twist $\alpha_1$ along $\gamma$,
and let $\RightFinf(\spinc)$ resp. $\LeftFinf(\spinc)$ denote the subset of
$\CFinfty(Y_0,\spinc)$, generated by the $\gamma$-induced intersection
points to the right resp. the left of the curve $\gamma$. Recall
that if we twist sufficiently, then $\LeftFinf(\spinc)$ is a
subcomplex (c.f. Proposition~\ref{prop:LeftIsSubcomplex}).

We relate $\HFinfty$ for $Y$ with
$H_*(\RightFinf)$, as follows:

\begin{lemma}
\label{lemma:RightInfty}
There is an isomorphism $H_*(\RightFinf)\cong \HFinfty(Y)$.
\end{lemma}

\begin{proof}
Let $\Theta_{\alpha,\gamma}\in\Ta\cap\Tc$ be the intersection point
$\{\gamma\cap\alpha_1,w_2^+,...,w_g^+\}$.  It follows as in the proof
of Proposition~\ref{prop:LeftIsSubcomplex} that there are no triangles
$\psi\in
\pi_2(\Theta_{\alpha,\gamma},\x,\y)$ with $\x\in\LeftFinf$,
$\y\in\Tc\cap\Tb$ and $\cald(\psi)\geq 0$, and $\Mas(\psi)=0$.  Hence,
counting holomorphic triangles whose $\Ta\cap\Tc$-vertex is
$\Theta_{\alpha,\gamma}$, we obtain a map
$H_*(\RightFinf)\longrightarrow \HFinfty(Y)$. On the chain level, this
map has the form $\iota + {\text{lower order}}$, where
$\iota[\x,i]=[\x',i-n_z(\psi_\x)]$ where $\x'$ is the intersection
point on $\Tc\cap \Tb$ closest to $x\in\Ta\cap\Tb$, $\psi_\x$ is the
unique small triangle (supported in the neighborhood of $\gamma$ and
the support of the isotopies between $\gamma_i$ and $\alpha_i$ with
non-negative multiplicities) and lower order is taken with respect to
the energy filtration on $Y$. Moreover, there is a relative
$\Z$-grading on both complexes, given by the Maslov index (where we
take an ``in'' domain for $Y_0$). The map preserves this
grading. Moreover, there are only finitely many generators in each
degree. It follows then that the induced map is an isomorphism.
\end{proof}

We have seen that the map $H_*(\RightFinf)\longrightarrow
H_*(\LeftFinf)$ naturally splits into two pieces, $\delta_1$ and
$\delta_2$, where $\delta_1$ uses the domains $\PhiIn$
from Lemma~\ref{lemma:TwoHoClasses}.

\begin{lemma}
\label{lemma:DeltaOne}
The map $\delta_1$ is an isomorphism.
\end{lemma}

\begin{proof}
This follows from the fact that on the chain level, $\delta_1$  has the form
$$\delta_1[\x^+,i]=[\x^-,i-n_z(\phi_{\x^+,\x^-})]+\text{lower order}.$$
(Lemma~\ref{lemma:DecomposeDelta}), together with the fact that $\delta_1$ preserves the relative $\Z$ grading.
\end{proof}

\begin{lemma}
\label{lemma:DeltaTwo}
The map $\delta_2$ is an isomorphism.
\end{lemma}

\begin{proof}
Fix an equivalence class of intersection points
between $\Ta\cap \Tb$, all of which are $\gamma$-induced. According
to Section~\ref{sec:EulerCharacteristic}, if we wind sufficiently many times along
$\gamma$ and move the basepoint $z$ sufficiently close to $\gamma$,
then $\langle c_1(\spinc), H \rangle$ can be made arbitrarily
large. By moving the basepoint to change the $\SpinC$ structure, we
have that the complexes $\LeftFp(\spinc)$ and $\LeftFp(\spinc')$
(resp. $\RightFp(\spinc)$ and $\RightFp(\spinc')$) are
identical. Moreover, if $\spinc$ and $\spinc'$ are sufficiently
positive, then the map $\delta_2^+$ is independent of the $\SpinC$
structure. 

Choose a degree $i$ sufficiently large that $H_i(\RightFp)\cong
H_i(\RightFinf)$ and $H_i(\LeftFp)\cong H_i(\LeftFinf)$, and note
under these identifications, the map induced on homology
$$ \delta^+_2\colon H_i(\RightFp)\longrightarrow H_{i-1}(\LeftFp)$$
agrees with $\delta_2$. For fixed $i$ and sufficiently large $\spinc$,
$\delta_1^+$ on $H_i(\RightFp(\spinc))$ vanishes. Since $\HFp(Y,\spinc)$ is zero for all sufficiently large
$\spinc$, it follows from the long exact sequence induced from
$$\begin{CD}
0@>>> \LeftFp(\spinc) @>>> \CFp(Y,\spinc) @>>> \RightFp(\spinc) @>>> 0
\end{CD}
$$ that $\delta=\delta^+_1+\delta^+_2\colon
H_*(\RightFp(\spinc))\longrightarrow H_*(\LeftFp(\spinc))$ is an
isomorphism. 
It follows  that the kernel of $\delta^+_2$ in degree
$i$ is trivial. 
From this, it follows in turn
that the kernel of $\delta^+_2$ is
trivial in all larger degrees. Since $\delta_1^+$ decreases degree
more than $\delta_2^+$, it is easy to see that the cokerenel of
$\delta_2^+$ in dimension $i$ is trivial, as well.
The lemma then follows.
\end{proof}

\vskip.2cm
\noindent{\bf{Proof of Proposition~\ref{prop:HFinfProp}.}}
Note that $\delta_1$ and and $\delta_2$ are both isomorphisms, and
$$\gr(\delta_1([\x,i]),\delta_2([\x,i]))=\pm 2n$$ for each generator
$[\x,i]$ for $\CFp(Y,\spinc)$.  It
follows that: $$\HFinfty(Y_0,\spinc)\cong H_*(\RightFinf)/(U^n-1).$$ Thus,
the proposition follows from Lemma~\ref{lemma:RightInfty}.
\qed

\vskip.2cm
\noindent{\bf Proof of Theorem~\ref{thm:HFinfZHS}.} 
Since multiplication by $U$ is an isomorphism on $\HFinf(Y,\spinc_0)$,
Proposition~\ref{prop:HFinfProp} shows that $\HFinf(Y)\cong
\HFinf(Y_1)$, where $Y_1$ denotes the $+1$ surgery on any knot
$K\subset Y$. Since any two integer homology three-spheres
can be connected by sequences of $\pm 1$ surgeries, 
it follows that $\HFinf(Y)\cong
\HFinfty(S^3)$, which we know has the claimed form.
\qed

\vskip.3cm
\noindent{\bf Proof of Theorem~\ref{thm:HFinfBigger}.} This is a
direct consequence of Theorem~\ref{thm:HFinfZHS} and
Proposition~\ref{prop:HFinfProp} when $c_1(\spinc)$ is non-torsion.
In the torsion case, the induced maps on homology satisfy either
$\delta_1=\delta_2$, or $\delta_1=-\delta_2$,
according to the two possible orientation conventions for $Y$. The two
possibilities give two different homology groups (over $\Z$). We
define the standard orientation convention to be the one for which
$\delta_1=-\delta_2$.

Finally, note that the action of $h\in H_1(Y_0;\Z)$ is given by $\pm
\delta_1$, as can be easily seen from the geometric representative for
the circle action (see Remark~\ref{HolDiskOne:rmk:GeoRep}
of~\cite{HolDisk}).
\qed

\subsection{The general case of Theorem~\ref{thm:HFinfGen}}
\label{subsec:HFinfGen}

\begin{defn}
Let $Z$ be a compact three-manifold with $\partial Z = T^2$.
The kernel of the map 
$$H_1(\partial Z) \longrightarrow H_1(Z)$$
is cyclic, generated by $d\ell$, where $\ell\subset T^2$ is a simple, 
closed curve. We call such a curve $\ell$ a {\em longitude}, and $d$ the {\em divisibility} of $Z$.
\end{defn}

\begin{prop}
\label{prop:InductiveStep}
Suppose that $b_1(Z)=1$, and let $h_1$, $h_2$ be primitive homology
classes in $H_1(T^2;\Z)$ and with $h_1\cm \ell$ and $h_2\cm \ell$ positive
with $h_1\cm h_2 = 1$. Then, if $\HFinf$ of $Y_{h_1}$ and $Y_{h_2}$
satisfy the property Theorem~\ref{thm:HFinfGen}, then so does
$Y_{h_1+h_2}$.
\end{prop}

\begin{proof}
Recall that the Floer homologies of a rational homology three-sphere have an absolute
$\Zmod{2}$ grading, specified by $$\chi(\HFa(Y))=|H_1(Y;\Z)|.$$ From
the exact sequence of Theorem~\ref{thm:GeneralSurgery}, we have
that $$
\begin{CD}
..@>>>\HFp(Y_{h_1}) @>{F_1}>> \HFp(Y_{h_2}) @>{F_2}>> \HFp(Y_{h_1+h_2}) @>>>...
\end{CD}
$$ The hypothesis in the sign guarantees that the degree shift occurs
at $F_1$ (using the absolute $\Zmod{2}$ grading on each group). It
follows that $\HFinfty(Y_{h_1+h_2})$ vanishes in all odd
degrees. Indeed, since this is true when we take coefficients in
$\Zmod{p}$ for all $p$; hence, $\HFinfty(Y_{h_1+h_2})$ has no torsion
in even degrees. Since $\chi(\HFinfty(Y,\spinc)/(U-1))=1$ for all rational
homology three-spheres, the result follows. 
\end{proof}

\begin{prop}
\label{prop:IndStepTwo}
Suppose that $Z$ be an oriented three-manifold with torus boundary.
For each $h$ with the property that $h\cm \ell = 1$, we have an
identification $$\HFinfty(Y_{\ell},\spinc)\cong
\HFinfty(Y_{h},\spinc_0)/(U^n-1) $$ where $\spinc_0$ is a torsion $\SpinC$
structure, $\spinc_0|_Z=\spinc|_Z$, and ${\mathfrak d}(\spinc)=2n$.
\end{prop}

\begin{proof}
We adapt the proof of Proposition~\ref{prop:HFinfProp}.  We start with
$(\Sigma,\{\alpha_2,...,\alpha_g\},\{\beta_1,...,\beta_g\})$
representing the knot complement $Z$, and then choose $\alpha_1$ to
represent $\ell$ and $\gamma$ to represent $h$:
i.e. $(\Sigma,\alphas,\betas)$ represents $Y_\ell$ and
$(\Sigma,\{\gamma,\alpha_2,...,\alpha_g\},\{\beta_1,...,\beta_g\})$
represents $Y_{h}$. There is an added feature now, since the
divisibility $d$ of $Z$ could be greater than one.  It is still the
case that for sufficiently large winding, all the intersection points
are represented from $\RightFinf(\spinc)$ or $\LeftFinf(\spinc)$, and,
as in Lemma~\ref{lemma:TwoHoClasses}, all homotopy classes of maps $\phi$
with $\Mas(\phi)=1$ admitting holomorphic representatives (connecting
any two intersection points) satisfy that the property that
$\partial_\alpha \phi$ uses the central point $p=\alpha_1\cap\gamma$
either once or zero times. Recall $\delta_1$ is the map defined using
those homotopy classes which meet $p$ once.  Now, there is a
difference map $$\eta\colon (\Ta\cap\Tb)\times
(\Ta\cap\Tb)\longrightarrow
\Zmod{d},$$ which is defined by
$$\eta(\x,\y)=\#\left(\partial_{\alpha_1}\phi\cap p\right) 
\pmod{d}.$$
There are corresponding splittings
\begin{eqnarray*}
\LeftFinf(\spinc)=\LeftFinf_1,...,\LeftFinf_d &{\text{and}}&
\RightFinf(\spinc)=\RightFinf_1,...,\RightFinf_d.
\end{eqnarray*}
labeled so that $\eta(\x,\y)=1$ if $\x\in\RightFinf_i$ and
$\y\in\RightFinf_{i+1}$, and 
$\delta_1(\RightFinf_i)\subset \LeftFinf_{i+1}$.
and $\delta_2(\RightFinf_i)\subset \LeftFinf_{i}$. 

The proof of Lemma~\ref{lemma:RightInfty} gives us that
$H_*(\RightFinf_i)\cong \HFinfty(Y,\spinc_0)$ (for $i=1,...,d$). 
Also, analogues of Lemmas~\ref{lemma:DeltaOne} and \ref{lemma:DeltaTwo} still hold: both $\delta_1$ and $\delta_2$ are isomorphisms.
Now, the proposition easily follows as before. 
\end{proof}

\vskip.2cm
\noindent{\bf{Proof of Theorem~\ref{thm:HFinfGen}.}}
We begin with the case where $b_1(Y)=0$, and prove the claim by
induction on $|H_1(Y;\Z)|$. The base case is, of course,
Theorem~\ref{thm:HFinfZHS}.  For the inductive step, we choose a knot
$K\subset Y$ which represents a non-trivial homology class. With
appropriate orientation, we have that $m\cdot \ell >0$. If $m\cdot
\ell>1$, the inductive step follows from
Proposition~\ref{prop:InductiveStep}, since $m$ can be decomposed as
$m=h_1+h_2$ with $h_1\cm h_2 = 1$, $h_1\cm \ell, h_2\cm \ell >1$. Note
also that if $h\cm\ell>0$, then $|H_1(Y_h)|$ depends linearly on $h\cm
\ell$.

If $m\cm \ell=1$, then since $K$ is homologically non-trivial, we must
have that $d>1$. Also, $|\Tors H_1(Y_\ell)|=\frac{1}{d}|\Tors
H_1(Y)|$.  Applying Proposition~\ref{prop:IndStepTwo} along a
different knot in $Y_\ell$ which represents a generator for
$H_1(Y_\ell)/\Tors$, we see that $$\HFinf(Y_\ell,\spinc)\cong
\HFinfty(Y',\spinc')/(U^n-1),$$ where $|H_1(Y';\Z)|<|H_1(Y;\Z)|$. Applying the
proposition again, and the induction hypothesis, 
we obtain that $\HFinf(Y)\cong
\Z[U,U^{-1}]$.

The proof for general $b_1(Y)=1$ or $2$ follows from an induction on
$b_1(Y)$. Let $Y$ be an oriented three-manifold with $b_1(Y)=1$ or
$2$. Choose a knot $K\subset Y$ whose image in $H_1(Y;\Z)/\Tors$ is
primitive. (This implies that in $Y-K$, the divisibility $d=1$.)  If
$\spinc$ is a non-torsion $\SpinC$ structure on $Y_\ell$, then the
result follows from Proposition~\ref{prop:IndStepTwo}. The other case
follows from the fact that we have two maps $\delta_1$ and $\delta_2$
from $\RightFinf(\spinc)$ to $\LeftFinf(\spinc)$, and both of these
maps are isomorphisms of $\Z[U]\otimes_\Z \Wedge^*
\left(H_1(Y_h;\Z)/\Tors\right)$-modules (between two modules are, in turn,
isomorphic
to $\Z[U^{-1}]\otimes_{\Z}\Wedge^*H^1(Y_h;\Z)$). Now, observe that the
automorphism of $\Z[U]\otimes_\Z \Wedge^*
\left(H_1(Y_h;\Z)/\Tors\right)$-module $\Z[U]\otimes_\Z \Wedge^*
\left(H_1(Y_h;\Z)/\Tors\right)$ is determined by its action on the
determinant line $\Wedge^b (H_1(Y_h;\Z)/\Tors)\cong \Z$, where it is
either multiplication by $+1$ or $-1$. Thus, the maps $\delta_1$
and $\delta_2$ either cancel (for one orientation convention) or
they do not (for the other one). The convention where $\delta_1+\delta_2=0$ is the one
for which the theorem follows; it is, in this case, the standard orientation convention for $Y$.
\qed

\subsection{The twisted case}
\label{subsec:HFinfTwist}

We state a version of Theorem~\ref{thm:HFinfGen} which holds for arbitrary first Betti number.

Observe that the proof of Theorem~\ref{thm:HFinfGen}
breaks down when $b_1(Y)\geq 3$, since
now the module $\Z[U]\otimes_\Z \Wedge^*(\Z^{b-1})$ has non-trivial
automorphisms, so that $\delta_1$ and $\delta_2$ do not necessarily
cancel. Indeed, it is proved in~\cite{HolDiskThree} that
$$\HFinf(T^3,\spinc_0) \cong \Z[U,U^{-1}]\otimes_\Z \Big(H^1(T^3)\oplus
H^2(T^3)\Big)$$
where $\spinc_0$ is the $\SpinC$ structure with
$c_1(\spinc_0)=0$.

There is, however, a version which holds for twisted coefficient systems.

Observe first that the twisted homology group
$\uHFinf(Y,\spinc)$ is a module over the group-ring
$\Z[H^1(Y;\Z)]\otimes_\Z \Z[U,U^{-1}]$ (which can be thought of as a
ring of Laurent polynomials in $b_1(Y)+1$ variables). To make the
ring structure respect the relative grading, we give
$\uHFinf(Y,\spinc_0)$ a relative $\Zmod{2}$ grading.

\begin{theorem}
\label{thm:HFinfTwist}
Let $Y$ be a closed, oriented three-manifold. Then, there is a unique
equivalence class of orientation system for which we have a
$\Z[U,U^{-1}]\otimes_\Z\Z[H^1(Y;\Z)]$-module isomorphism for
each torsion $\SpinC$ structure $\spinc_0$ on $Y$:
$$\uHFinf(Y,\spinc_0)\cong
\Z[U,U^{-1}],$$ 
where here the latter group is endowed with a trivial action by $H^1(Y;\Z)$.
\end{theorem}

\begin{proof}
The proof is obtained by modifying the above proof of
Theorem~\ref{thm:HFinfGen}, with minor modifications, which we outline presently.

For the case where $H_1(Y_0;\Z)\cong \Z$, we adapt the proof of
Theorem~\ref{thm:HFinfBigger}, thinking of $\Z[H^1(Y;\Z)]$ as
$\Z[T,T^{-1}]$. In this case, Lemma~\ref{lemma:RightInfty} is replaced
by an isomorphism $H_*(\RightFinf)\cong \HFinfty(Y)[T,T^{-1}]$ (with
the same proof). Next, we observe that rather than having $\delta_1$
and $\delta_2$ cancel, as in the proof of
Theorem~\ref{thm:HFinfBigger}, we have that $\delta_1=\pm \delta_2\cm
T$. In fact, for some choice of orientation convention, we can arrange for
$\delta_1=-\delta_2$.
The result then follows easily from the long exact sequence
connecting $\uLeftFinf(Y,\spinc)$, $\uHFinf(Y,\spinc)$, and
$\uRightFinf(\spinc)$ observing that the map
$$\Z[T,T^{-1}]\stackrel{1-T}{\longrightarrow} \Z[T,T^{-1}]$$
injective, with cokernel $\Z$ (with trivial action by $T$).

The same modifications work to prove the general case
(arbitrary $b_1(Y)$) as well.

We now turn to the uniqueness assertion on the orientation system.
For the various equivalence classes of orientation systems, it is
always true that $\uHFinf(Y,\spinc_0)\cong \Z[U,U^{-1}]$ as a $\Z$
module. In fact, we saw (c.f. Equation~\eqref{eq:ChangeOrientations})
that as a $\Z$ module, the isomorphism class of the chain complex
$\uCFinf(Y,\spinc_0)$ is independent of the choice of
orientation system. Moreover, from Equation~\eqref{eq:ChangeOrientations},
it is clear that the
$2^{b_1(Y)}$ different equivalence classes of coherent orientation
system give rise to all $2^{b_1(Y)}$ different $\Z[H^1(Y;\Z)]$-module
structures on $\Z[U,U^{-1}]$ which correspond
naturally to $\Hom(H^1(Y;\Z),\Zmod{2})$, with a distinguished module
for which the action by $H^1(Y;\Z)$ is trivial.
\end{proof}

\begin{remark} In fact, the above argument shows in general that
for any $\SpinC$ structure over $Y$, there is an identification of
$\Z[U,U^{-1}]$ modules $\uHFinf(Y,\spinc_0)\cong\Z[U,U^{-1}]$. However,
the action of $\xi\in H^1(Y;\Z)$ will, in general, be given by
multiplication by $U^k$, where
$k$ is given by $2 k=\langle \xi\cup c_1(\spinc),[Y]\rangle $.
\end{remark}

\subsection{Absolute $\Zmod{2}$ gradings}
\label{subsec:AbsoluteGradings}

With the help of Theorem~\ref{thm:HFinfTwist}, we can define an absolute
$\Zmod{2}$ grading on $\CFinf(Y,\spinc)$ (and hence all the other
associated chain complexes), for all $\SpinC$ structures,
simultaneously.

We declare the non-zero generators of $\uHFinf(Y,\spinc)$ to
have even degree.  Note that for a rational homology three-sphere,
this orientation convention agrees with that used before, i.e.
$\chi(\HFa(Y))=|H_1(Y;\Z)|$. (In fact, if we orient $\Ta$ and $\Tb$ so that the intersection number
$\#(\Ta\cap\Tb)=|H_1(Y;\Z)|$, then the $\Zmod{2}$ grading at a generator $[\x,i]$ is $+1$ if and only if the
local intersection number of $\Ta$ and $\Tb$ at $\x$ is $+1$.)

%       Our convention is the following. Let $\omega\in
%       \Lambda^{b} (H_1(Y;\Z)/\Tors)\cong \Z$ be a generator, 
%       and choose a torsion $\SpinC$ structure $\spinc_0$.  We declare
%       non-zero elements in $\omega\cm \HFinf(Y,\spinc)$ to 
%       have even degree. Note that for a rational homology
%       three-sphere, this orientation convention agrees with that used
%       before, i.e.  $\chi(\HFa(Y))=|H_1(Y;\Z)|$.  With this convention, we
%       can orient $\Ta$ and $\Tb$ so that the local intersection number
%       between $\Ta$ and $\Tb$ at a point $\x$ is $+1$, when $[\x,i]$ is an
%       even-dimensional element of $\CFinfty(Y,\spinc_0)$.  Indeed, these
%       orientations for $\Ta$ and $\Tb$, induce an absolute $\Zmod{2}$
%       grading on all the $\SpinC$ structures simultaneously.

With this orientation convention, we have the following 
refinement of Corollary~\ref{cor:AlexPoly}:

\begin{prop}
\label{prop:PreciseChi}
Let $Y_0$ be an oriented three-manifold with $b_1(Y_0)=1$, and
$\spinc$ be a non-torsion $\SpinC$ structure, then
$$\chi(\HFp(Y_0,\spinc_0+nH))=-\tau_t(Y_0,\spinc),$$ where $t$ is the
component containing $c_1(\spinc)$, and the sign on
$\tau_t(Y_0,\spinc)$ is specified by
$$\tau_{-t}(\spinc)-\tau_t(\spinc)=n.$$ In particular, if $Y_0$ is
obtained by zero-surgery on a knot $K$ in a homology three-sphere,
whose symmetrized Alexander polynomial is
$$\Delta_K=a_0+\sum_{i=1}^d a_i(T^i+T^{-i}),$$ then 
$$ \chi(\HFp(Y_0,\spinc_0+nH)=-\sum_{j=1}^d j a_{|n|+j}.$$
\end{prop}

\begin{proof}
First observe that the sign comparing $\chi(\HFp(Y_0))$ and $\tau_t$
in Theorem~\ref{thm:EulerOne} is universal, depending on the relative
sign between $\Delta_{i,j}$ and $\Delta'_{i,j}$. Checking these signs
for $S^1\times S^2$, the Proposition follows.
\end{proof}

\subsection{The Euler characteristic of $\HFm$}
\label{subsec:EulerHFm}

The following is an immediate consequence of
Theorem~\ref{thm:EulerOne}, together with
Theorem~\ref{thm:HFinfBigger} (though a more direct proof can be given
by modifying the discussion in Section~\ref{sec:EulerCharacteristic}):

\begin{cor}
Let $Y$ be an oriented three-manifold with $b_1(Y)=1$, and $\spinc\in
\SpinC(Y)$ be a non-torsion $\SpinC$ structure. Then,
$\chi(\HFm(Y,\spinc))=\tau_{-t}(\spinc)$, where $t$ is the component of
$H^2(Y;\Z)-0$ containing $c_1(\spinc)$
\end{cor}

\begin{proof}
The short exact sequence 
$$
\begin{CD}
0 @>>>\CFm(Y,\spinc)@>>>\CFinf(Y,\spinc)@>>>\CFp(Y,\spinc)@>>>0.
\end{CD}
$$
induced a long exact sequence in homology
$$
\begin{CD}
@>>> \HFm(Y,\spinc)@>>>\HFinf(Y,\spinc)@>>>\HFp(Y,\spinc)@>>>...,
\end{CD}
$$ which shows that
$$\chi(\HFinf(Y,\spinc))=\chi(\HFp(Y,\spinc))+\chi(\HFm(Y,\spinc)).$$
Moreover, Theorem~\ref{thm:HFinfGen} implies that
$$\chi(\HFinf(Y,\spinc))=n,$$ where $2n$ is the divisibility of
$c_1(\spinc)$ in $H^2(Y,\spinc)/\Tors$. The result now follows from
the ``wall-crossing formula'':
$$\tau_{-t}(Y,\spinc)-\tau_{t}(Y,\spinc)=n$$ for Turaev's torsion
(see~\cite{Turaev}).
\end{proof}

\begin{cor}
If $Y$ is an oriented three-manifold with $b_1(Y)=1$ or $2$ and
$\spinc\in\SpinC(Y)$ is a non-torsion $\SpinC$ structure, then
$\chi\left(\HFm(Y,\spinc)\right)=\pm\tau(\spinc)$.
\end{cor}

\begin{proof}
This follows in the same manner as the previous corollary, except that
now $c_1(\spinc)^\perp$ is a non-trivial vector space, so its exterior
algebra has Euler characteristic zero: thus, $\chi(\HFinfty(Y,\spinc))=0$.
\end{proof}

\subsection{The truncated Euler characteristic}
\label{subsec:TruncEuler}

In Theorem~\ref{thm:EulerOne}, we worked with a non-torsion $\SpinC$
structure. The reason for this, of course, is given
Theorem~\ref{thm:HFinfGen}: if $\spinc_0$ is torsion and $Y_0$ is a
three manifold with $0<b_1(Y)=b\leq 2$, then in all sufficiently large degrees
$i$, $\HFp_i(Y_0,\spinc_0)\cong \HFinf_i(Y_0,\spinc_0)\cong
\Z^{2^{b_1(Y)-1}}$. This shows, however, that for all sufficiently large $n$,
the Euler characteristic of the graded Abelian group $\HFp_{\leq
n}(Y_0,\spinc_0)$ takes on two possible values, 
depending on the parity of $n$ (and the difference between the two
values is $2^{b_1(Y)-1}$). In fact, we have the following:

\begin{theorem}
\label{thm:TruncEuler}
Let $Y$ be a three-manifold with $b_1(Y)=1$ or $2$, equipped with a torsion $\SpinC$ structure $\spinc_0$.
Then, when $b_1(Y)=1$, then for all sufficiently large $n$
$$\chi(\HFp_{\leq n}(Y,\spinc_0))=
\left\{\begin{array}{ll}
-\tau(Y) & {\text{for odd $n$}} \\
-\tau(Y)+1 & {\text{for even $n$}}
\end{array}\right.$$
When $b_1(Y)=2$, then in all sufficiently large degrees,
$$\chi(\HFp_{\leq n}(Y,\spinc_0))=\pm \tau(Y) + 
(-1)^n.$$
\end{theorem}

\begin{proof}
As before, we have a short exact sequence
$$
\begin{CD}
0@>>>\LeftFp @>>>\CFp(Y_0,\spinc_0)@>>>\RightFp @>>> 0,
\end{CD}
$$
and hence a long exact sequence:
$$
\begin{CD}
... @>>> H_i(\LeftFp) @>>> \HFp_i(Y,\spinc_0) @>>> H_i(\RightFp) @>{\delta}>> ...
\end{CD}
$$ 
Note that we are using a relative $\Z$ grading here, which we can do since $\spinc_0$ is torsion.
When $i$ is sufficiently large, the coboundary map $\delta$ is
zero, since on $\HFinfty$, the map $H_*(\LeftFinf)\longrightarrow
\HFinf(Y)$ is an injection.

It follows that for all sufficiently large $n$,
\begin{equation}
\label{eq:LeftRightEuler}
\chi(\HFp_{\leq n}(Y))=\chi(H_{\leq n}(\LeftFp))+\chi(H_{\leq n}(\RightFp)).
\end{equation}

On the other hand, we still have a short exact sequence:
$$ 
\begin{CD}
0 @>>> {\ker f_1} @>>> \RightFp @>{f_1}>> \LeftFp @>>> 0,
\end{CD}
$$
inducing
$$
\begin{CD}
@>>> H_i({\ker f_1}) @>>> H_i(\RightFp) @>{f_1}>> H_{i-1}(\LeftFp) @>>> ...
\end{CD}
$$ Note that with the earlier grading conventions, $f_1$ must decrease
the grading by one.  Of course,  $\ker f_1$ is a
finite-dimensional graded vector space, so the above gives the
following relation for all sufficiently large $n$: 
\begin{equation}
\label{eq:EulerfOne}
\chi({\ker f_1})=\chi(H_{\leq n}(\RightFp)) + \chi(H_{\leq n-1}(\LeftFp)).
\end{equation}
But from Proposition~\ref{prop:TuraevTorsion} applies in the present
case, to identify $\chi({\ker f_1})=\tau(\spinc_0)$. Note that the
proof of the that proposition does not really require that $\spinc$ be
negative; it suffices to consider the case where $\spinc+\alpha_1^*$,
$\spinc+\beta_j^*$ and $\spinc+\alpha_1^*+\beta_j^*$ are negative, and
$c_1(\spinc)$ is torsion. Combining this result,
Equation~\eqref{eq:LeftRightEuler}, and Equation~\eqref{eq:EulerfOne},
we obtain that: $$\chi(\HFp_{\leq n}(Y,\spinc_0))= -\tau(Y,\spinc_0)+
(-1)^n \rk H_{n}(\LeftFp,\spinc_0).$$ Suppose that $b_1(Y)=1$. Then,
(according to Theorem~\ref{thm:HFinfGen}) for all sufficiently large
$n$, $\rk H_n(\LeftFp,\spinc_0)=1$ if $n$ is even and $0$ when $n$ is
odd. Similarly, when $b_1(Y)=2$, 
we have $$\rk H_n(\LeftFp,\spinc_0)=
\rk \HFinf_n(Y)/2 = 1.$$
\end{proof}

\subsection{On the role of $n_z$}

The ``triviality'' of $\HFinf(Y)$ -- its dependence on the homological
information of $Y$ alone -- underscores the importance of the
quantity $n_z$ in the construction of interesting Floer-homological
invariants. 

Another manifestation of this is the following.  When $Y$ is an
integral homology three-sphere, we needed the base-point to define
$\Z$-grading between intersection points. However, there is still a
$\Zmod{2}$ graded-theory $CF'(Y)$, which is freely generated by the
transverse intersection points of $\Ta\cap\Tb$, and $\Zmod{2}$-graded
by the local intersection number between $\Ta$ and $\Tb$. The map
$$\partial \x = \sum_{\y}\sum_{\{\phi\in\pi_2(\x,\y)|\Mas(\phi)=1\}}
\left(\#\UnparModFlow(\phi)\right)\y$$
gives a well-defined boundary map, and in fact, we can consider the
homology group 
$$\HF'(Y)=H_*(\CF'(Y),\partial).$$ 

However, it is a consequence of Theorem~\ref{thm:HFinfZHS} that 
$$\HF'_*(Y)\cong \Z\oplus 0.$$ To see this, note that as a
$\Z/2$-graded chain complex, $\CFinfty(Y)$ is
naturally a (finitely generated, free) module over the 
ring of Laurent polynomials
$\Z[U,U^{-1}]$. Moreover, its quotient by the action of $U$ and $U^{-1}$ is
the complex $\CF'(Y)$ defined above. More algebraically, we have that
$$\CF'(Y)=\CFinf(Y)\otimes_{\Z[U,U^{-1}]}\Z,$$
where the homomorphism $\Z[U,U^{-1}] \longrightarrow \Z$ sends $U$ to
$1$. Theorem~\ref{thm:HFinfZHS} says that $\HFinfty(Y)$ is a free
$\Z[U,U^{-1}]$-module of rank one.
The claim about $\HF'_*(Y)$ then follows immediately from
the universal coefficients theorem spectral sequence 
(see, for instance~\cite{EilenbergMacLane}).

\section{Applications}
\label{sec:Applications}

In this section, we prove the remaining results 
(Theorems~\ref{thm:Complexity} and \ref{thm:BoundIntPtsSurg})
claimed in the introduction. 

\subsection{Complexity of three-manifolds}

The theorems in the introduction dealing with fractional surgeries are
proved using surgery exact sequences with twisted theories
(Theorems~\ref{thm:ExactFrac} and \ref{thm:ExactFraca}). 
Consequently, we will need the following
analogue of Theorem~\ref{thm:EulerOne}
for the twisted theory:

\begin{lemma}
\label{lemma:TwistedEuler}
Let $Y_0$ be a homology $S^1\times S^2$, and choose a coefficient
system corresponding to a representation
$$H^1(Y_0;\Z)\longrightarrow \Zmod{n}.$$
Then, for each non-torsion $\SpinC$ structure over $Y_0$, we have that
$$\chi(\uHFp(Y_0,\Zmod{n};\spinc))=
n\cm \chi(\HFp(Y_0,\spinc))=
- n\cm \tau_t(Y_0,\spinc)$$  (where on the left we are
still taking the rank as a $\Z$-module, and $t$ here is the component of 
$H^2(Y;\Z)-0$ containing $c_1(\spinc)$).
Similarly, for a torsion $\SpinC$ structure $\spinc_0$, we have that
$$\chi(\HFp_{\leq 2n+1}(Y_0,\spinc_0;\Zmod{n})
= -n\cm \tau(Y_0,\spinc_0).$$
\end{lemma}

\begin{proof}
The proof proceeds exactly as in the proof of
Theorem~\ref{thm:EulerOne} (with the sign pinned down in
Proposition~\ref{prop:PreciseChi}, and 
Theorem~\ref{thm:TruncEuler} in the case where the $\SpinC$ structure
is torsion), together with the observation that now
$\chi(\Ker f_1)$ multiplies by $n$. 
\end{proof}

We will also need the following result, which follows along the lines
of Section~\ref{sec:HFinfty}.

\begin{lemma}
\label{lemma:TwistedHFinf}
Suppose that $Y_0$ is a homology $S^1\times S^2$, and choose a
coefficient system corresponding to a map $H^1(Y_0;\Z)\cong \Z
\longrightarrow
\Zmod{n}$ which maps generators to generators. Then, if $\spinc_0$ is a torsion $\SpinC$ structure, 
then ${\underline{\HF}^\infty_i}(Y_0,\spinc_0,\Zmod{n})\cong \Z$ in all degrees. 
\end{lemma}

\begin{proof}
We still have the long exact sequence
$$
\begin{CD}
... @>>> \uHFinf(Y_0,\spinc_0,\Zmod{n})@>>> 
H_*(\uRightFinf,\Zmod{n}) @>{\delta}>> H_*(\uLeftFinf,\Zmod{n}) @>>> ...
\end{CD}
$$ We place a reference point $p$ at the intersection of $\gamma$ (the
perturbing curve) with $\alpha_1$.  It is clear that
$H_*(\uLeftFinf,\Zmod{n})\cong
H_*(\LeftFinf)\otimes_{\Z}\Z[\Zmod{n}]$.  Moreover, the coboundary
splits as $\delta=\delta_1 - \zeta \delta_2$, where $\zeta$ is is a
primitive $n^{th}$ root of unity, and $\delta_1$ and $\delta_2$
are the maps obtained from the $\delta_1$ and $\delta_2$ using $\Z$
coefficients, by a base-change to $\Zmod{n}$. In particular, both
$\delta_1$ and $\delta_2$ are isomorphisms
(Lemmas~\ref{lemma:DeltaOne} and \ref{lemma:DeltaTwo}).  Thus, in view
of Theorem~\ref{thm:HFinfGen} (indeed, we're using here the special
cases from Subsection~\ref{subsec:HFinfSmallY}), we have exactness for
$$
0 \longrightarrow {\underline{\HFinf_i}}(Y_0,\spinc_0,\Zmod{n}) 
\longrightarrow
\Z[\Zmod{n}] 
\stackrel{1-\zeta}{\longrightarrow} 
\Z[\Zmod{n}] 
\longrightarrow
\underline{\HFinf_{i-1}}(Y_0,\spinc_0,\Zmod{n}) 
\longrightarrow 0
$$
\end{proof}

We can now prove Theorem~\ref{thm:Complexity}.

\vskip.2cm
\noindent{\bf{Proof of Theorem~\ref{thm:Complexity}.}}
This is an application of the $U$-equivariant exact sequence of
Theorem~\ref{thm:ExactFrac}, which gives:
$$\begin{CD}
... @>{F_1}>>  \uHFp(Y_0;\Zmod{n})@>{F_2}>> \HFp(Y_{1/n})@>{F_3}>> \HFp(Y)
@>>> ...,
\end{CD}.
$$ 

Now, we claim that for all sufficiently large $d$, the map induced
by $F_2$ $$\Image U^d \uHFp(Y_0,\Zmod{n}) \longrightarrow \Image U^d
\HFp(Y_{1/n})$$ is surjective. 
It suffices to consider the $\spinc_0$-summand of $\uHFp(Y_0,\Zmod{n})$, where
$\spinc_0$ is the torsion $\SpinC$ structure. There, $F_2$ has a
natural $\Z$-graded lift. For one parity, the corresponding
$\HFinfty(Y_{1/n})$ vanishes (so the claim is obvious). For the other parity, 
in sufficiently high degree $k$, the image of $F_1$ is trivial, so, with the help of 
Lemma~\ref{lemma:TwistedHFinf}, our exact sequence reads:
$$
\begin{CD}
0 @>>> {\underline{\HFp_k}}(Y_0,\spinc_0;\Zmod{n})\cong 
{\underline{\HFinf_k}}(Y_0,\spinc_0;\Zmod{n})\cong \Z @>{F_2}>> \HFp_{k} (Y_{1/n})\cong \Z.
\end{CD}
$$
Since $\HFinf(Y)$ has no torsion, it easily follows that $F_2$ must surject onto the 
generator in $\HFp_k(Y_{1/n})$. 

From this observation, together with the $U$-equivariant exact sequence,
it follows that the map
$$ \begin{CD}
\frac{\HFp(Y)}{U^d \HFp(Y)} @>>>
\frac{\uHFp(Y_0,\Zmod{n})}{U^d \uHFp(Y_0,\Zmod{n})} @>>>
\frac{\HFp(Y_{1/n})}{U^d \HFp(Y_{1/n})} .
\end{CD}
$$
is exact in the middle, and hence that
\begin{equation}
\label{eq:BoundN}
\rk\left(\uHFred(Y_0,\Zmod{n})\right)
\leq 
\rk\left(\uHFred(Y)\right)+
\rk\left(\uHFred(Y_1)\right).
\end{equation}
(Here, as in the case where $b_1=0$, 
$\uHFred(Y_0,\Zmod{n})$ is defined to be the quotient of $\uHFp(Y_0,\Zmod{n})$ by 
the image of $\uHFinf(Y_0,\Zmod{n})$.)

Now, observe that if $\spinc\neq \spinc_0$,
$\HFp(Y_0,\spinc;\Zmod{n})$ is finitely generated, so that for
sufficiently large $d$,
\begin{equation}
\label{eq:ReduceDoesNothing}
\uHFred(Y_0,\spinc;\Zmod{n})=\frac{\uHFp(Y_0,\spinc;\Zmod{n})}{U^d\uHFp(Y_0,\spinc,\Zmod{n})}
=\uHFp(Y_0,\spinc;\Zmod{n}).
\end{equation}

For $\spinc=\spinc_0$, we observe that
\begin{equation}
\label{eq:BoundTrunc}
\max(0,-\chi(\HFp_{\leq 2n+1}(Y_0,\spinc_0;\Zmod{n})))
\leq \rk \HFp_{\leq 0}(Y_0,\spinc_0;\Zmod{n}).
\end{equation}
The reason for this is that for all sufficiently large $n$, we have 
\begin{eqnarray*}
\lefteqn{\chi(\uHFp_{\leq 2n+1}(Y_0,\spinc_0;\Zmod{n}))
=} \\ 
&& \chi\left(\uHFred(Y_0,\spinc_0;\Zmod{n})\right) 
+ \chi\left(\uHFp_{\leq 2n+1}(Y_0,\spinc_0;\Zmod{n})\cap \Image\uHFinf(Y_0,\spinc_0;\Zmod{n})
\right).
\end{eqnarray*}
The second term above is negative: owing to the algebraic structure of
$\HFinf(Y_0,\spinc_0;\Zmod{n})$ (the even-dimensional generators are the images of
the odd-dimensional ones under an isomorphism), there
are more odd-dimensional
than even-dimensional 
generators coming from $U^d \uHFp(Y_0,\spinc_0;\Zmod{n})$ in $\uHFp_{\leq 2n+1}f(Y_0,\spinc_0;\Zmod{n})$. 

The theorem is obtained by combining Inequality~\eqref{eq:BoundN},
Equation~\eqref{eq:ReduceDoesNothing}, Inequality~\eqref{eq:BoundTrunc}, and 
Lemma~\ref{lemma:TwistedEuler}.
\qed

\subsection{Gradient trajectories}

We turn to the bounds on the simultaneous trajectory number of an
integral homology three-sphere discussed in the introduction. First,
we dispatch with Theorem~\ref{thm:BoundIntPts} of from the introduction:

\vskip.2cm
\noindent{\bf{Proof of Theorem~\ref{thm:BoundIntPts}.}}
This is clear: if $(\Sigma,\alphas,\betas,z)$ is a pointed Heegaard
diagram for $Y$, where the $\alpha_i$ meet the $\beta_j$ in general
position, the intersection corresponding chain complex $\CFa(Y)$ is
freely generated by intersection points $\Ta\cap\Tb$, and its rank
is bounded below by the rank of its homology.
\qed
\vskip.2cm

We turn to Theorem~\ref{thm:BoundIntPtsSurg}.

\vskip.2cm
\noindent{\bf{Proof of Theorem~\ref{thm:BoundIntPtsSurg}.}}
As a first step, observe that, since $$\chi(\uHFp(Y_0,\spinc_0\pm
iH;\Zmod{n}))= \pm n\cm \MT_i(K),$$ it follows that the rank of
$\HFp(Y_0,\Zmod{n},\spinc)$ is non-zero for at least $2k$ distinct
non-torsion $\SpinC$ structures; thus the rank of
$\uHFa(Y_0,\spinc,\Zmod{n})$ is also non-zero in these $\SpinC$
structures (c.f. Proposition~\ref{prop:NonVanishHFa}). Moreover, from
Lemma~\ref{lemma:TwistedHFinf}, it follows that the rank of
$\uHFp(Y_0,\Zmod{n},\spinc_0)$ is non-zero, and hence so is the rank
of $\uHFa(Y_0,\spinc_0,\Zmod{n})$.  Now, since for all $\SpinC$
structures, $$\chi(\uHFa(Y_0,\spinc,\Zmod{n}))=0$$ (again, using the
twisted analogue of Proposition~\ref{prop:EulerHFa}), the rank of
$\HFa(Y_0,\Zmod{n})$ is at least $4k+2$.  The result then follows from
the exact sequence of Theorem~\ref{thm:ExactFraca}, together with 
Theorem~\ref{thm:BoundIntPts}.
\qed

\commentable{
\bibliographystyle{plain}
\bibliography{biblio}
}
\end{document}